\documentclass[11pt]{amsart}
\usepackage{amsmath,amsxtra,amssymb,amsthm,amsfonts}

\def\theequation{\thesection.\arabic{equation}}

\numberwithin{equation}{section}
\newtheorem{lemma}{Lemma}[section]
\newtheorem{prop}[lemma]{Proposition}
\newtheorem{theorem}[lemma]{Theorem}
\newtheorem{cor}[lemma]{Corollary}

\newtheorem{rem}[lemma]{Remark}

\newcommand{\re}{\begin{rem}\rm}
  \newcommand{\mar}{\end{rem}}
\newtheorem{exam}[lemma]{Example}

\newtheorem{defi}[lemma]{Definition}

\newcommand{\kla}{\left ( }
\newcommand{\mer}{\right ) }

\renewcommand{\for}{\begin{eqnarray*}}

\newcommand{\mel}{\end{eqnarray*}}

\newcommand{\kl}{\pl \le \pl}
\newcommand{\gl}{\pl \ge \pl}

\newcommand{\lel}{\pl = \pl}

\newcommand{\ez}{{\mathbb E}}

\newcommand{\rz}{{\mathbb R}}
\newcommand{\zz}{{\mathbb Z}}

\newcommand{\cz}{{\mathbb C}}

\newcommand{\ten}{\otimes}

\DeclareMathOperator{\dom}{dom}

\DeclareMathOperator{\rg}{rg}

\newcommand{\p}{\hspace{.05cm}}
\newcommand{\pl}{\hspace{.1cm}}
\newcommand{\pll}{\hspace{.3cm}}

\newcommand{\qd}{\end{proof}\vspace{0.5ex}}

\newcommand{\Om}{\Omega}
\newcommand{\om}{\omega}

\newcommand{\al}{\alpha}
\newcommand{\si}{\sigma}
\newcommand{\Si}{\Sigma}

\newcommand{\La}{\Lambda}
\newcommand{\la}{\lambda}
\newcommand{\eps}{\varepsilon}

\newcommand{\A}{{\mathcal A}}

\newcommand{\M}{{\mathcal M}}

\newcommand{\N}{{\mathcal N}}

\newcommand{\U}{{\mathcal U}}

\newcommand{\xspace}{\hbox{\kern-2.5pt}}
\newcommand{\xyspace}{\hbox{\kern-1.1pt}}

\newcommand\tnorm[1]{\left\vert\xspace\left\vert\xspace\left\vert\mskip2mu
#1\mskip2mu \right\vert\xspace\right\vert\xspace\right\vert}

\newcommand{\supn}{{\sup_n}^+}

\newcommand{\pr}{{\mathrm P\mathrm r}}

\allowdisplaybreaks

\newcommand{\ttt}{{\bf t}}

\begin{document}
\title[Noncommutative Riesz transforms]{\bf \Large Noncommutative Riesz transforms --\\ a probabilistic approach}

\author[M. Junge]{M. Junge}
\address{Department of Mathematics\\
University of Illinois, Urbana, IL 61801, USA} \email[Marius
Junge]{junge@math.uiuc.edu}

\author[T. Mei]{T. Mei}
\address{Department of Mathematics\\
University of Illinois, Urbana, IL 61801, USA} \email[Tao
Mei]{mei@math.uiuc.edu}

\begin{abstract} For $2\le p<\infty$ we show the lower estimates
 \[ \|A^{\frac 12}x\|_p \kl c(p)\max\{\pl \|\Gamma(x,x)^{\frac{1}{2}}\|_p,\pl \|\Gamma(x^*,x^*)^{\frac{1}{2}}\|_p\}  \]
for the Riesz transform associated to a  semigroup $(T_t)$ of
completely positive maps on a von Neumann algebra with negative
generator $T_t=e^{-tA}$, and gradient form
\[ 2\Gamma(x,y)\lel Ax^*y+x^*Ay-A(x^*y)\pl .\]
As additional hypothesis we assume that $\Gamma^2\gl 0$ and the
existence of a Markov dilation for $(T_t)$. We give applications to
quantum metric spaces and show the equivalence of semigroup Hardy
norms and martingale Hardy norms derived from the Markov dilation.
In the limiting case we obtain a viable definition of BMO spaces for
general semigroups of completely positive maps which can be used as
an endpoint for interpolation. For torsion free ordered groups we
construct a connection between Riesz transforms and the Hilbert
transform induced by the order.
\end{abstract}

\thanks{ The first author is partially supported by the
National Science Foundation Foundation DMS 05-56120. The second
author was partially supported by a Young Investigator Award of the
N.S.F supported summer workshop in Texas A\&M university 2007.}

\maketitle

{\bf Introduction:}

\noindent Riesz transforms provide important examples in classical
harmonic analysis and have been studied extensively in the
literature in many different aspects. The aim of this paper is to
continue the work of P. A. Meyer, Bakry, Emery, Ledoux, Varopoulos
and many others on probabilistic aspects of the theory of Riesz
transforms, however in the noncommutative setting. The importance of
analyzing semigroups of completely positive maps on von Neumann
algebras has been impressively demonstrated by the recent work of
Popa \cite{Po-an}, Peterson \cite{Pe}, Popa and Ozawa \cite{op} and
also occurs in the work of Shlyahktenko/Connes \cite{CS} on Betti
numbers for von Neumann algebras. A common thread in this analysis
is to adapt some differential geometric concepts in the setting of
von Neumann algebras.

It was discovered by P.A. Meyer that the general theory of
semigroups provides an appropriate framework to formulate Riesz
transforms which relate the norm of different derivatives in the
classical setting. To be more precise we consider a family $(T_t)$
of contractive completely positive maps on a finite von Neumann
algebra $N$ with normal faithful trace $\tau$ such that
 \[ \tau(T_tx)\le \tau(x) \]
holds for positive $x$ and $t>0$. In the classical setting this is
certainly satisfied for a semigroup of measure preserving positive
maps on $N=L_{\infty}(\Om,\Si,\mu)$. Then the maps $T_t$ act on all
$L_p$ spaces $L_p(N,\tau)$ and in particular on the Hilbert space
$L_2=L_2(N,\tau)$. Let $A$ be the negative generator of
$T_t=e^{-tA}$. For nice elements the gradient form
 \[ 2\Gamma(x,x) \lel A(x^*)y+x^*A(y)-A(x^*y) \]
is defined. In all our examples we will assume that the $T_t$'s are
selfadjoint, i.e. $\tau(T_txy)=\tau(xT_ty)$ and then $A$ is indeed a
positive (unbounded) operator. Under this circumstances we can
formulate P.A. Meyer's problem: Is it true that
 \begin{equation}\label{Mey0}
  \|\Gamma(x,x)^{1/2}\|_p \sim _{c(p)} \|A^{1/2}x\|_p
  \end{equation}
holds for all reasonable elements $x$?

Let us illustrate this question by considering the Laplace operator
$A(f)=-\Delta(f)$, where $\Delta=\sum_i \frac{\partial^2}{\partial
x_i^2}$. Then it is easily verified that
 \[ \Gamma(f,h) \lel \sum_{i=1}^n
  \overline{\frac{\partial f}{\partial x_i}}\frac{\partial f}{\partial x_i}  \lel (\nabla f,\nabla h) \pl
 .\]
In this context P.A. Meyer's inequality reads as follows
 \[ \| |\nabla f|  \|_p \sim_{c(p)} \| |\Delta|^{1/2}(f)\|_p \pl
 .\]
In dimension $n=1$ this follows easily from the continuity of the
Hilbert transform. In higher dimension this are the first examples
of singular integrals and we refer to Stein's work
\cite{Stein,Stein2} for credentials and further information. Indeed,
P.A. Meyer was strongly motivated to provide a probabilistic
approach to Stein's work on Riesz transforms and succeeded in
showing his estimate for the Ornstein-Uhlenbeck semigroup, even in
the infinite dimensional case \cite{Me1,Me2}. In some sense Bakry
\cite{Ba0,Ba1,Ba2,Ba3,Ba4,Ba4} continued Meyer's line of research
and showed \eqref{Mey0} for many diffusion semigroups satisfying the
$\Gamma^2\gl 0$ condition. In the context of semigroups given by the
Laplace-Betrami operator on a manifold the positivity of $\Gamma^2$
is equivalent to the positivity of the Ricci curvature. We recall
that
 \[ 2\Gamma^2(x,y) \lel \Gamma(Ax,y)+\Gamma(x,Ay)-A\Gamma(x,y) \pl
 .\]
More generally the higher order gradients are defined as
 \[ 2\Gamma^{k+1}(x,y) \lel \Gamma^k(Ax,y)+\Gamma^k(x,Ay)-A\Gamma^k(x,y) \pl
 .\]
The Bochner identities for manifolds show that
 \[ \Gamma^2(f,f) \lel {\rm Ric}(df,df)+\|\nabla df\|_{HS}^2 \pl
 .\]
Here $\nabla$ is the second covariant derivative and $HS$ stands for
the Hilbert Schmidt of the corresponding matrix of second
derivatives.

In the noncommutative setting the notion of diffusion process is not
(yet) well-defined. It is however clear that Meyer's approach
requires the semigroup to have a \emph{Markov dilation}. This means
that there exists a family of homomorphisms $\pi_s:N\to \tilde{N}$
such that $\pi_s(N)$ is contained in a filtration $\tilde{N}_t$
with conditional expectation $E_t$ such that
 \[ E_t(\pi_s(x))\lel \pi_t(T_{s-t}x) \]
holds for $t<s$. Our main result is one half of P.A. Meyer's
inequality for $p\gl 2$.

\begin{theorem}\label{main} Let $(T_t)$ be a semigroup of
completely positive selfadjoint  maps with Markov dilation and
$\Gamma^2\gl 0$. Let $2\le p<\infty$ then
 \[ \|A^{1/2}x\|_p \kl c(p) \max\{ \|\Gamma(x,x)^{\frac 12}\|_p,\|\Gamma(x^*,x^*)^{1/2}\|_p\}
 \pl .\]
\end{theorem}

The assumptions of Theorem \ref{main} are satisfied for Fourier
multipliers on discrete groups, sometimes also called Herz-Schur
multipliers. Indeed, let $G$ be a discrete group and $VN(G)$ the
group von Neumann algebra given by the left regular representation
$\la:G\to B(\ell_2(G))$, $\la(g)\delta_h=\delta_{gh}$. It is well
known that $\tau(\sum_g a_g \la(g))=a_e$ extends to a normal
faithful trace on $VN(G)$. A Fourier multiplier is given by a
semigroup $(\phi_t)$ of positive definite functions and the normal
extensions of
 \[ T_t(\la(g))\lel \phi_t(g)\la(g) \pl .\]
Following the recent work of Ricard, we know that $T_t$ has a Markov
dilation provided that $\phi_t(1)=1$, $\phi_t$ is real valued and
$\phi_t(g)=\phi_t(g^{-1})$, i.e. in the selfadjoint case. Let us
note that then, according to Schoenberg's theorem, there exists a
conditionally negative function $\psi:G\to [0,\infty)$. We obtain an
immediate application to quantum metric spaces. Recall that a triple
$(\A,B,\|\pl\|_{Lip})$ is a quantum metric space, if $\A\subset B$
is a dense, not necessarily closed, sub$^*$-algebra of a
$C^*$-algebra $B$ and $\tnorm{\pl}_{Lip}$ is a norm on $\A$ such
that
 \[ \tnorm{ab}_{Lip}\kl \tnorm{a}_{Lip}\|b\|_{B}+\|a\|_{B}\tnorm{b}_{Lip} \pl
 \]
and that
 \[ d_{Lip}(\phi,\psi)\lel \sup_{|\hspace{-0.03cm}\|a\|\hspace{-0.02cm}|_{Lip}\le 1}
 |\phi(a)-\phi(b)| \]
induces the weak$^*$ topology on the state space of $B$.
Equivalently, the natural quotient  map $\iota:
(\A,\tnorm{\pl}_{Lip})\to B/\cz1$ is a compact operator (see
\cite{OR}). Theorem \ref{main} allows us to show the compactness
condition for discrete groups with rapid decay. Let us recall that a
finitely generated group has rapid decay if
 \[ \|x\|_{\infty} \kl C(s) k^s \|x\|_2 \]
holds for some $s$ and every  $x=\sum_{|g|=k} a_g \al(g)$ supported
on the words of length $k$. This notion is independent of the choice
of the generators.

\begin{cor}\label{hyp} Let $G$ be a finitely generated
discrete group with rapid decay and  $\psi:G\to \cz$  be a
conditionally negative function such that $\psi(1)=0$,
$\psi(g)=\psi(g^{-1})$, and
 \[ \inf_{|g|=k }  |\psi(g)| \gl c_{\al} k^{\al} \]
holds for some $\al>0$. Then $\cz[G]$ equipped with the norm
 \[ \tnorm x  \lel \|\Gamma_{A}(x,x)\|^{\frac12} \]
defines a quantum metric space $(\cz[G],C_{red}(G),\tnorm\pl)$.
\end{cor}

In view of the examples from Riemann manifolds it is clear that
$\tnorm\pl$ is the ``correct'' Lipschitz norm.

Compared to Bakry's work in the commutative case the assumption
``diffusion semigroup'' is deleted in our Theorem \ref{main} (but we
have to keep the assumption of a nice algebra of invariant
functions). In section 2 we provide a first proof of Theorem
\ref{main} requiring some extra regularity assumptions (easily
verified for Fourier multipliers). In section 3 we follow
Meyer/Bakry's footsteps. Unfortunately, this also requires to
develop the theory of $H_p$ spaces for noncommutative continuous
filtrations, see \cite{JK}. A key ingredient in  Meyer's martingale
approach is the use of a stopped  brownian motion $P_{B_s}$ as
indices for the associated Poisson semigroup. More explicitly, let
$\pi_t$ be a Markov dilation, $a>0$, $B_t$ a Brownian motion with
$B_0=a$ almost everywhere  and $\ttt_a(\om)=\inf\{t: B_t(\om)=0\}$
the hitting time for the boundary. Then Meyer's approach consists in
investigating the martingale
 \[ \rho_a(x) \lel \pi_{\ttt_a}(x) \pl .\]
The $H_p^c$ norm of the resulting martingale decomposes in a time
and a space component. Without going into details let us mention
that we are able to ``compare''  the martingale $H_p^c$ norm in
Meyer's model and show that they there are (almost) equivalent to
the Hardy norms for semigroups investigated in the joint work C. Le
Merdy and Q. Xu, see \cite{JLX}. This explicit connection between
martingale and semigroup $H_p$ norms seems to be new and extends to
the $H_p$ norms given by the reversed Markov filtration. 

Moreover, we are able to define spaces of bounded mean oscillation
for quite general semigroups as follows
 \[\|x\|_{BMO_c(T)}\lel   \sup  \|T_t|x|^2-|T_tx|^2\|_{\infty}^{1/2} \pl .\]
This norm is closely related to Garsia's BMO norm for the Poisson
semigroup on the circle. The space $BMO$ is defined as the
completion of $N$ with respect to the norm $\|x\|_{BMO}\lel
\max\{\|x\|_{BMO_c},\|x^*\|_{BMO_c}\}$. The good news is that the
BMO space serves as an endpoint for interpolation.

\begin{theorem}\label{intpp} Let $(T_t)$ be a semigroup  with Markov dilation
and $\Gamma^2\gl 0$. Then
 \[ [BMO,L_p(N)]_{\frac{1}{q}} \lel L_{pq}(N) \pl .\]
\end{theorem}

For the associated Poisson semigroup $P_t=e^{-tA^{1/2}}$ the
connection between Garsia's BMO norm and Meyer's martingale approach
is very explicit
 \[ \|x\|_{BMO_c(P)} \lel \|\rho_a(x)\|_{BMO_c} \]
holds for every $a$. This allows us to reduce the interpolation
theorem to previous results on martingales. It is also a good
indication that we have found an appropriate BMO spaces  for
semigroups. In section 4 we confirm this observation by showing that
many natural candidates for BMO norms are equivalent. Section 6 is
devoted to applications of Theorem \ref{main} and Theorem
\ref{intpp}. The applications towards quantum metric spaces are
prepared in Section 5. Our last applications concerns torsion free
ordered groups  which admit a filtration of normal divisors
$G=G_0\supset G_1\supset G_2 \supset \cdots $ such that
 \[ \bigcap_k G_k \lel \{e\} \quad ,\quad G_k/G_{k+1} \lel \zz \pl
 .\]
This holds for example for free groups in $n$ generators. Using the
extension $id\ten P_t^{\zz}$ of the classical Poisson group on
$G\times \zz$ we are able to reduce the boundedness of the Hilbert
transform for ordered groups to estimates for Riesz transforms
associated with $P_t$. This gives a link between the $H_p$-theory
related to sub-diagonal von Neumann algebras and the $H_p$-theory
for semigroups, see Section 6 for more details. In the text the
absolute constant $c(p)$ may differ from line to line.

\section{Preliminaries and notation}
We will use standard notation in the theory of operator algebras
which can be found in \cite{Tak,TK-II,TK-III}, \cite{kar-I,kar-II}
or \cite{sz,strat}. Modular theory does not play an important role
in this paper, because in most of our applications we are interested
in von Neumann algebras with a finite trace. At any rate, using the
Haagerup reduction method (see \cite{JHX}) it usually suffices to
consider the finite case. As a standard reference for noncommutative
$L_p$-spaces we refer to \cite{PX-hand} and the references therein.
For basic properties of  the space of $\tau$-measurable operators
and noncommutative integration we refer to \cite{Ne}. We will also
use operator space terminology, in particular the notion of
completely bounded maps, see the books by Effros-Ruan \cite{ER},
Pisier \cite{psop} or Paulsen \cite{Pa}. We allow for a slight
deviation in the notion of \emph{completely bounded} maps $T:X\to
Y$, where $X\subset L_p(N)$, $Y\subset L_p(M)$ are subspaces of a
noncommutative $L_p$ space. Indeed, we use
 \[ \|T\|_{cb} \lel \sup_{\M} \|id\ten T:L_p(\M;X)\to L_p(\M;Y) \]
where the supremum is taken over all von Neumann algebras $\M$ and
the space
 \[ L_p(\M;X)\subset L_p(\M\ten N) \]
is the completion of the tensor product $L_p(\M)\ten X$ with respect
to the induced norm. In the usual definition of the cb-norm, the
supremum is only taken over $\M=K(\ell_2)$, the compact operators on
$\ell_2$. If Connes' embedding conjecture were true the two norms
coincide. Our policy in general is to prove the estimates with
respect to the stronger norm. Indeed, as so often in martingale
theory these estimates are automatic, i.e. follow because $T$ and
$id\ten T$ satisfy the same assumptions.

We will frequently use standard tools from noncommutative
probability, in particular Doob's inequality
 \[ \| \supn E_n(x)\|_p \kl d_p \|x\|_p \]
for $1<p\le \infty$ and the dual Doob inequality
 \[ \|\sum_n E_n(x_n)\|_p \kl c_p \|\sum_n x_n\|_p \]
which holds for $1\le p<\infty$ and positive $x_n$ with the constant
$c_p=d_{p'}$. Here $(E_n)$ is a sequence of normal conditional
expectations onto an increasing or decreasing sequence of von
Neumann subalgebras $(N_n)$ of a given von Neumann $N$. In the
finite setting these conditional expectations are unique and hence
commute, i.e. $E_nE_m=E_mE_n=E_{\min(n,m)}$ (increasing filtration)
or $E_nE_m=E_mE_n=E_{\max(n,m)}$ (decreasing filtration). We will
always assume this commutation relation. The notation $\supn$ is
taken from \cite{JD} and \cite{JX}. In the noncommutative setting
the pointwise supremum can not be defined directly. However, for
selfadjoint operators $x_n$ we have an order analogue
 \[ \|\supn x_n\|_p \lel \inf_{-a\le x_n\le a} \|a\|_p \pl .\]
In full generality we use Pisier's definition
 \[ \|\supn x_n\|_p \lel \inf_{x_n=ay_nb} \|a\|_{2p} \sup_n
 \|y_n\|_{\infty} \|b\|_{2p} \pl .\]
The same definition holds for arbitrary index sets. We also need
some basic facts about $H_p$-spaces for martingales. Let us recall
some definitions from \cite{JX2}. As usual the martingale difference
are denoted by $d_kx=E_kx-E_{k-1}x$.
 \begin{align*}
 \|x\|_{H_p^c} &= \|(\sum_k d_kx^*d_kx^{})^{\frac12}\|_p \pl ,
 \|x\|_{h_p^c}  \lel \|(\sum_k E_{k-1}(d_kx^*d_kx^{}))^{\frac12}\|_p \pl ,
 \|x\|_{h_p^d} \lel  (\sum_k \|d_kx\|_p^p)^{\frac1p} \pl .
 \end{align*}
The row versions are  given by $\|x\|_{H_p^r}=\|x^*\|_{H_p^c}$,
$\|x\|_{h_p^r}=\|x^*\|_{h_p^c}$. In the noncommutative theory the
definition of the $H_p$-spaces is as follows
 \[ H_p \lel \begin{cases} H_p^c\cap H_p^r &\mbox{if } p\gl 2\\
  H_p^c+H_p^r &\mbox{if } 1\le p \le 2\end{cases} \pl .\]
The Burkholder Gundy inequality reads as
 \begin{equation}\label{BG}
  L_p(N) \lel    H_p  \quad 1<p<\infty \pl .
  \end{equation}
The Burkholder inequalities can be formulated as
\begin{equation}\label{Bk}
 L_p(N) \lel  \begin{cases}
 h_p^c\cap h_p^r\cap h_p^d & p\gl 2 \pl, \\
 h_p^c+h_p^d+h_p^r &1<p\le 2 \pl .
 \end{cases} \end{equation}
All the equalities hold  with equivalent norms. Since this will be
needed in the paper we want to show $H_p=h_p^d+h_p^c$ for $1\le p\le
2$. This requires us to use the dual norms
 \[ \|x\|_{L_p^cMO} \lel \|\supn
 E_n(|x-E_{n-1}(x)|^2)\|_{p/2}^{1/2}\pl ,\pl
 \|x\|_{L_p^cmo} \lel \|\supn
 E_n(|x-E_{n}(x)|^2)\|_{p/2}^{1/2} \pl .\]
Extending the Fefferman Stein duality $(\overline{H_1^c})^*=BMO_c$
from \cite{PX} it was shown in \cite{JX} that
 \[ \overline{H_p^c}^* \lel L_{p'}^cMO \quad 1\le p<2 \pl .\]
Here $\overline{X}^*$  refers to the anti-linear duality $\langle
x,y\rangle=tr(x^*y)$. The following observation is probably of
independent interest.
\begin{lemma}\label{h11} Let $1\le p<2$ and $(N_k)_{k\gl 1}$ be a
martingale filtration
 \begin{enumerate}
 \item[i)] Let $L_p^{c,cond}=\{(x_k)_k\pl:\pl x_k\in N_k\}\subset
 L_p(N;\ell_2^c)$. Then the antilinear dual of  $L_p^{c,cond}$ is
isomorphic to the space $L_{p'}^{c,cond}MO$ of sequences $(x_k)$
with $x_k\in L_{p'}(N_k)$ such that
  \[ \|(x_k)\|_{L_{p'}^{c,cond}MO}\lel \|\sup_n E_n(\sum_{k\gl n}
  x_k^*x_k)\|_{p'/2}^{1/2} \pl .\]
 \item[ii)] Let $h_p^{c,\circ}$ be the subspace of
 $h_p^c$ of elements with $d_1=0$. The anti-linear dual of
 $h_p^{c,\circ}$ is $L_{p'}^cmo$.
 \item[iii)] $H_p^c\lel h_p^d+h_p^c$.
 \end{enumerate}
\end{lemma}
\begin{proof} In \cite{JD} it is shown that
 \begin{eqnarray}\label{fss}
  |\sum_k tr(x_k^*y_k)|
 \kl  \sqrt{2} \pl  \|(\sum_k E_k(x_k^*x_k))^{\frac12}\|_p \pl
 \|\sup_n E_n(\sum_{k\gl n} y_k^*y_k)\|_{p'/2}^{\frac12} \pl .
 \end{eqnarray}
In particular, we have a continuous inclusion
$L_{p'}^{c,cond}MO\subset (\overline{L_{p}^{c,cond}})^*$. For the
converse we note that $L_p^{c,cond}$ is a subspace of
$L_p(N,\ell_2^c)$. Hence a linear functional
$f:\overline{L_{p}^{c,cond}}\to \cz$  of norm one is given by a
sequence $(z_k)\subset L_{p'}(N,\ell_2^c)$ such that
 \[ f(x)\lel \sum_k tr(x_k^*z_k) \pl .\]
We define $y_k=E_k(z_k)$ and deduce from Doob's inequality for
$p'/2>1$ that
 \begin{align*}
  \|\sup_n E_n(\sum_{k\gl n}  y_k^*y_k)\|_{p'/2}
  &\le \|\sup_n E_n(\sum_{k\gl n}  E_k(z_k^*z_k)) \|_{p'/2}
  \kl  \|\sup_n E_n(\sum_k z_k^*z_k) \|_{p'/2}\\
 &\le  d_{p'/2} \|\sum_k z_k^*z_k\|_{p'/2} \pl .
  \end{align*}
For the proof of ii) we assume $d_1(x)=0$ or $d_1(y)=0$ and note
that according to \eqref{fss} we have
  \begin{align*}
  |tr(x^*y)| &=  |\sum_{k\gl 2} tr(d_k(x)^*d_k(y))|\\
  & \kl \sqrt{2} \pl \|(\sum_{k\gl 2}
  E_{k-1}(d_k(x)^*d_k(x)))^{\frac12}\|_p\pl
 \|\sup_{n\gl 1} E_n(\sum_{k-1\gl n} d_k(y)^*d_k(y))\|_{p'/2}^{\frac12}
 \pl .
 \end{align*}
Let us recall that $L_p^cmo$ consists of martingales  with
$d_1(y)=0$. It has been proved in \cite{JD}  that there are linear
maps $u_{k}:N\to C\bar{\ten}(N_{k})$ (the space of weakly converging
columns with values in $N_k$) such that
 \[ u_k(x)^*u_k(x)\lel E_k(x^*x) \pl .\]
Moreover, $u_k$ is a $N_k$ right module map with complemented range
(see \cite{JD} and for the non-separable case \cite{JS}). Then
$u:h_p^{c,\circ}\to L_{p}^{c,cond}$ given by $u(x)\lel
(u_{k-1}(d_k(x))_{k\gl 2}$ is an isometric isomorphism. Hence an
antilinear functional $f:\overline{h_p^{c,\circ}}\to \cz$ is given
by a sequence $z_k\in L_{p'}(N_{k-1},\ell_2^c)$ such that
 \[ f(u(x))\lel \sum_k tr(u_{k-1}(d_k(x))^*z_k) \]
and $\|\sup_n E_n(\sum_{k-1\gl n}z_k^*z_k)\|_{p'/2}\le c_{p}\|f\|$.
Since the range of $u_{k-1}$ is complemented, we may use the
projection and find $y_k$ such that $u_{k-1}(y_k)=Pz_k$ and
 \begin{align*}
   &\|\sup_n E_n(\sum_{k-1\gl n} y_k^*y_k) \|_{p'/2}
   \lel \|\sup_n E_n(\sum_{k-1\gl n} E_{k-1}(y_k^*y_k))
   \|_{p'/2}\\
   &= \|\sup_n E_n(\sum_{k-1\gl n}  u_{k-1}(y_k)^*u_{k-1}(y_k))\|_{p'/2}
 \kl  \|\sup_n E_n(\sum_{k-1\gl n} z_k^*z_k)\|_{p'/2} \kl c_p
 \|f\| \pl .
 \end{align*}
We define $y=\sum_k d_k(y_k)$. Using  the triangle inequality and
$E_nE_k=E_n$  we get that
 \begin{align*}
 &\|\sup_n E_n(\sum_{k-1\gl n} d_k(y)^*d_k(y))
 \|_{p'/2}^{\frac12}\\
 &\le \|\sup_n E_n(\sum_{k-1\gl n} E_k(y_k)^*E_k(y)) \|_{p'/2}^{\frac12}
 + \|\sup_n E_n(\sum_{k-1\gl n} E_{k-1}(y_k)^*E_{k-1}(y)) \|_{p'/2}^{\frac12}\\
 &\kl 2 \pl  \|\sup_n E_n(\sum_{k-1\gl n} y_k^*y_k) \|_{p'/2}
 \kl 2c_p \pl \|f\|
 \pl .
 \end{align*}
For the proof of iii) we recall that $h_p^d+h_p^c\subset H_p^c$. We
claim that the unit ball $B_{H_p^c}$ is contained in the norm
closure of $\overline{C(B_{h_p^d}+B_{h_p^c})}^{\|\pl\|_{H_p^c}}$ for
some constant $C>0$. If not there exists $x\in H_p^c$  and a
continuous linear functional $y$ such that $tr(y^*x)=1$ and
 \[ |tr(y^*x')|\kl \frac{1}{C} \quad \mbox{ for all } \quad x'\in CB_{h_p^d}\cup
 CB_{h_p^c}\pl .\]
We know that $y=\sum_n d_n$ satisfies
 \begin{align*}
  &\|y\|_{(\overline{H}_p^c)^*} \kl \sqrt{2}
   \|y\|_{L_p^cMO}^2
    \kl \sqrt{2} \pl \|\sup_n d_n^*d_n\|_{p'/2}
  + \sqrt{2} \pl \|\sup_n E_n(\sum_{k>n} d_k^*d_k)\|_{p'/2} \\
  &\le \sqrt{2} \kla\sum_n \|d_n\|_{p'}^{p'}\mer^{\frac{1}{p'}}+ \sqrt{2} \|\sup_n E_n(\sum_{k>n} d_k^*d_k)\|_{p'/2}
  \\
  &\kl 2\sqrt{2} \pl \|y\|_{\overline{h_p^d}^*} + \sqrt{2}c_p
  \|y\|_{(\overline{h_p^c})^*}
 \kl \frac{\sqrt{2}(2+c_p)}{C} \pl .
  \end{align*}
Since $\|x\|_{H_p^c}\le 1=tr(y^*x)$ we reach a contradiction for
$C>\sqrt{2}(2+c_p)$. An approximation argument allows us to replace
the norm closure of $C(B_{h_p^d}+B_{h_p^c})$ by
$2C(B_{h_p^d}+B_{h_p^c})$. \qd

Let us briefly prove  the martingale inequalities for potentials in
the noncommutative setting. We recall a classical martingale
inequality from \cite{JD} which is derived from \eqref{fss}. Let
$(N_k)$ be a (discrete) increasing filtration and $a_k\in N$ be
positive elements. For $p\gl 2$ we have
 \begin{equation}\label{ddoob}
 \|\sum_k E_k(a_k)\|_{\frac p2}
 \kl 2c_p^2 \pl \|\sup_m \sum_{k\gl m} E_m(a_k) \|_{\frac p2} \pl .
 \end{equation}
Here $c_p$ is the constant in Stein's inequality. Let $(z_k)_{k\le
n}$ be a finite submartingale, i.e.
 \[ E_k(z_{k+1})\gl z_k \pl .\]
The corresponding positive increasing part of $z$ is defined as
 \[ a_k \lel \langle z\rangle_k \lel \sum_{j<k} E_{j}(z_{j+1}-z_j) \pl .\]
Clearly, we obtain a martingale
 \[ m_k \lel z_k-a_k  \pl .\]
Indeed, $m_k-m_{k-1}=z_k-z_{k-1}-E_{k-1}(z_k-z_{k-1})$ is a
martingale difference sequence.  In the language of potentials, we
have
 \[ E_j(z_{j+1}-z_j)\lel E_j(m_{j+1}-m_j)+E_j(a_{j+1}-a_j) \lel E_j(a_{j+1}-a_j)
 \pl \]
Moreover, we note that $a_{j+1}-a_j$ is still positive. Hence for
$p>1$, we deduce from \eqref{ddoob} and Doob's inequality
 \begin{align*}
 &\|a_n\|_p \lel  \|\sum_{j=1}^n E_j(z_{j+1}-z_j)\|_p
 \lel   \|\sum_{j=1}^n E_j(a_{j+1}-a_j)\|_p \\
 &\le 2
 c_{2p}^2 \|{\sup_m}^+ \sum_{m\kl j\le n}^n E_m(a_{j+1}-a_j)\|_p\kl
 2c_{2p}^2 \|{\sup_m}^+ E_m(z_n-z_m)\|_p \\
 &\le
 2c_{2p}^2 \|{\sup_m}^+ E_m(z_n)\|_p + \|\sup_m z_m \|_p \kl
 2c_{2p}^2 d_p \|z_n \|_p + \|\sup_m z_m \|_p  \pl .
 \end{align*}
If in addition $z_m\gl 0$ for all $m$, we may ignore the second term
and obtain the following result.

\begin{lemma}\label{potential} Let $z_k=a_k+m_k$ be a positive submartingale with increasing
part $(a_k)$ and martingale part $(m_k)$. Then
\[   \|a_n \|_p \kl c(p)\pl  \| z_n\|_p \]
holds for $1<p<\infty$ and  some universal constant $c(p)$.
\end{lemma}

The $H_p$ theory for continuous filtrations $(N_t)_{t\gl 0}\subset
N$ has only been considered recently (see \cite{JK}). We will always
assume that $\bigcap_{s>t}N_s=N_t$. It is well-known that the theory
of $H_p$-norms  is closely related to stochastic integrals. However,
given how nicely the $H_p$-theory translates in the noncommutative
setting, we should not expect surprises. There are two candidates
for the $H_p^c$-norm on a finite interval $[0,T]$
 \[ \|x\|_{H_p^c} \lel \lim_{\si,\U}
 \|\sum_{j=0}^{|\si|-1} |E_{s_{j+1}}x-E_{s_j}x|^2\|_{p/2}^{1/2}
  \quad, \quad \|x\|_{\hat{H}_p^c}\lel
 \| \lim_{\si,\U} \sum_{j=0}^{|\si|-1}
 |E_{s_j+1}x-E_{s_j}x|^2\|_{p/2}^{1/2} \pl .\]
Here $\si=\{0=s_0,...,s_n=T\}$ is a partition $|\si|=n$ and $\U$ is
an ultrafilter refining the natural order of inclusion on the set of
all partitions. In the second term we take the weak$^*$-limit (at
least for $p\gl 2$). It was shown that in \cite{JK} that the two
norms are equivalent and, up to a constant $c_p$, independent of the
choice of $\U$. The main tool in this argument is the observation
that $H_p^c=h_p^d\cap h_p^c$ for $p\gl 2$, where
 \[ \|x\|_{h_p^d} \lel \lim_{\si,\U} (\sum_{j=0}^{n-1}
 \|E_{s_{j+1}}x-E_{s_j}x\|_p^p)^{\frac1p}
 \quad ,\quad
 \|x\|_{h_p^c} \lel \lim_{\si,\U} \|\sum_{j=0}^{n-1}E_{s_{j-1}}
 |E_{s_{j+1}}x-E_{s_j}x|^2\|_{p/2}^{1/2} \pl .\]
Again it is shown that the limit can be taken inside. This gives the
norm $\|\cdot\|_{\hat h_p^c}$ and conditioned bracket
 \[ \langle x,x\rangle_{T} \lel \lim_{\si,\U}
 \sum_{j=0}^{n-1}E_{s_{j-1}} |E_{s_{j+1}}x-E_{s_j}x|^2 \pl \]
and
 \[ \|x\|_{h_p^c} \sim _{c(p)} \| \langle x,x\rangle_T
 \|_{p/2}^{1/2} \lel \|x\|_{\hat{h}_p^c}  \pl .\]
In the continuous context the Burkholder inequalities reads as
follows
 \[ L_p(N) \lel \begin{cases} \hat{h}_p^c\cap \hat{h}_p^r\cap h_p^d &\mbox{if } p\gl 2\\
  \hat{h}_{p}^c+ \hat{h}_{p}^r+ h_{p}^d &\mbox{if } 1<p\le 2
  \end{cases} \pl ,\]
where $h_p^r,\hat h_p^r$ are the corresponding row spaces. The exact
form of the Feffermann-Stein duality for $p=1$ is not yet explored.
For $2\geq p>1$ we used the definition
$\hat{h}_{p}^c=\overline{\hat{h}_{p'}^c}^*$  and refer to \cite{JK}
for more information. A martingale $x$ is said to have \emph{a.u.
continuous path} if for every $T>0$, every $\eps>0$ there exists a
projection $e$ with $\tau(1-e)<\eps$ such that the function
$f_e:[0,T]\to N$ given by
 \[ f_e(t) \lel x_te \in N \]
is continuous. For a martingale with {a.u.\!\!}  continuous path we
have ${\rm var}_p(x)=\|x\|_{h_p^d}=0$ for all $2<p<\infty$.  We
recall from \cite{JK} that the condition ${\rm var}_p(x)=0$ implies
that
 \[ \lim_{\si,\U} \|\sum_{j=0}^{n-1} |E_{s_j+1}x-E_{s_j}x|^2-
 \sum_{j=0}^{n-1} E_{s_j}(|E_{s_j+1}x-E_{s_j}x|^2)\|_{p/2}
 \lel 0 \]
for all $p>2$ and the norm convergence of
  \[ L_{p/2}-\lim_{\si}  \sum_{j=0}^{n-1} E_{s_j}(|E_{s_j+1}x-E_{s_j}x|^2)
  \lel \langle x,x\rangle_T \pl .\]
This implies that for  martingales with ${\rm var}_p(x)=0$ the
equality (see \cite{JK})
 \[ \|x\|_{H_p^c([0,T])} \lel \|\langle x,x\rangle_T\|_{p/2}^{1/2} \]
holds without constants. Also  we have
 \begin{eqnarray}\label{brackket}
  \lim_{\si,\U} \|{\sup_j}^+
  E_{s_j}(|E_Tx-E_{s_{j-1}}x|^2)\|_{p/2}
  \lel \| {\sup_s}^+ E_{s}(\langle x,x\rangle_T-\langle
  x,x\rangle_{s}) \|_{p/2}  \pl .
  \end{eqnarray}
The correct definition of the norm in  $L_p^cmo$ for continuous
filtration is
 \[ \|x\|_{L_p^cmo} \lel \sup_{T} \| {\sup_s}^+ E_{s}(\langle x,x\rangle_T-\langle
  x,x\rangle_{s}) \|_{p/2}^{1/2}\]
for $2\le p\le \infty$.

It is shown in \cite{JX3} that the ergodic averages
$M_tx=\frac{1}{t}\int_0^t T_sxds$ satisfy a maximal inequality
 \begin{eqnarray}\label{Tdp}
 \|\sup_t M_t(x)\|_p \kl d_p \|x\|_p \quad 1<p\le \infty
 \end{eqnarray}
and the dual inequality
 \begin{eqnarray}\label{Tddp}
 \|\sum_k M_{t_k}(x_k)\|_p \kl c_p \|\sum_k x_k\|_p \quad 1\le
 p<\infty\end{eqnarray}
for positive $x_k$. We will make extensive use of the Poisson
semigroup
 \begin{eqnarray}\label{ps-form}
  P_s \lel \frac 1{2\sqrt{\pi }}\int_0^\infty se^{-\frac{s^2}{4u}}u^{-\frac
 32}T_udu.
 \end{eqnarray}
It has been shown in \cite{JX3} that $P_s$ is a positive average of
the $M_t$ and hence \eqref{Tdp} and \eqref{Tddp} hold for $P_t$
instead of $M_t$. If we assume in addition that the $T_t$'s are
selfadjoint, i.e.
 \begin{eqnarray}\label{sa}
 \tau(T_txy)\lel \tau(xT_ty)
 \end{eqnarray}
holds for all $x,y\in N$, then  \eqref{Tdp} and \eqref{Tddp} holds
for the $T_t$'s instead of the $M_t$, see \cite{JX2}. We will impose
a \emph{standard assumption}, namely the existence of
weak$^*$-dense, not necessarily  closed $^*$-algebra $\A\subset N$
such that
 \[ T_t(\A)\subset \A \quad \mbox{and} \quad  A(\A)\subset \A \]
holds for $t>0$ and the negative generator $A$ of the semigroup
$T_t=e^{-tA}$. Whenever we talk about $P_t$ we will also assume that
$P_t(\A)\subset \A$. This assumption is satisfied for our main
examples, Fourier multipliers, but somewhat problematic for certain
manifolds. We are convinced that with some extra effort our results
still hold for these manifolds and most of the interesting examples,
but the arguments are far more transparent with this additional
assumption. We will make crucial use of the following inequality
from \cite{Mei}.

\begin{prop}\label{monot} Let $x\in \A$ be positive  and $0<t<s$. Then
 \[ P_sx\kl  \frac{s}{t} \pl  P_tx  \pl .\]
\end{prop}

\begin{proof} We use \eqref{ps-form} and $e^{-\frac{s^2}{4u}}\le
e^{-\frac{t^2}{4u}}$ for all $u$. This yields the assertion
\begin{align*}
 \frac{P_sx}s &=\frac 1{2\sqrt{\pi }}\int_0^\infty e^{-\frac{s^2}{4u}}u^{-\frac 32}T_u(x)du  \kl  \frac 1{2\sqrt{\pi }}\int_0^\infty
 e^{-\frac{t^2}{4u}}u^{-\frac 32}T_u(x)du \lel
 \frac{P_tx}t \pl . \qedhere
\end{align*}
 \qd

We will use $H_p(T), H_p(P)$ to denote the Hardy spaces associated
with semigroup $(T_t)$ and subordinated Poisson semigroup $(P_t)$
respectively. See \cite{JLX} for more details.
\section{Reversed martingale filtration}

We will assume that $(T_t)$ is a semigroup of completely positive
maps on a tracial von Neumann algebra $N$. A Markov dilation for
$T_t$ is given by a family $\pi_s:N\to \M$ of trace preserving
$^*$-homomorphism  with the following properties
 \begin{enumerate}
 \item[i)] Let $M_{s]}$ be the von Neumann algebra generated by
 the $\pi_v(x)$, $x\in N$, $v\le s$. Then
  \[ E_{s]}(\pi_t(x))\lel \pi_s(T_{t-s}x) \]
holds for $s<t$ and $x\in N$.
 \item[ii)] Let $M_{[s}$ be the von Neumann algebra generated by
 the $\pi_v(x)$, $x\in N$, $v\gl s$. Then
  \[ E_{[s}(\pi_t(x))\lel \pi_s(T_{s-t}x) \]
holds for $t<s$ and $x\in N$.
\end{enumerate}
This definition is adapted for selfadjoint $(T_t)$. We recall that
$\Gamma^2\gl 0$ is equivalent to
 \[ \Gamma(T_tx,T_tx)\kl T_t\Gamma(x,x) \pl .\]
Whenever we invoke the associated Poisson semigroup
$P_tx=e^{-tA^{1/2}}x$ we assume in addition that $P_t(\A)\subset
\A$. It is easy to see that $\Gamma^2\gl 0$ also implies
 \[ \Gamma(P_tx,P_tx)\kl P_t\Gamma(x,x) \pl .\]

\begin{lemma}\label{ctpath} Let $(\pi_s)$ be a Markov dilation.
\begin{enumerate}
\item[i)]   Then $m(x)$ defined by
 \[ m_s(x)\lel \pi_s(T_sx) \]
is a martingale with respect to the reversed filtration $M_{[s}$.
\item[ii)] Assume $\Gamma^2\gl 0$ or $\Gamma(T_rx,T_rx)$ is
uniformly bounded in $L_2$ for $x\in \A$. Then for $x\in N$ the
martingale $m_s(x)$ has continuous path. If moreover,
$f(r,s)=T_{s}\Gamma(T_rx,T_rx)$ is continuous in $L_2$, or
$\Gamma^3\gl 0$, or $\Gamma^2(T_rx,T_rx)$ is locally bounded in
$L_2$, then
 \[ \langle m(x),m(y)\rangle \lel 2\int_0^{\infty} \pi_t(\Gamma(T_tx,T_tx)) dt \pl . \]
\end{enumerate}
\end{lemma}

\begin{proof} The first assertion is obvious, because
$E_{[s}(\pi_0(x))=\pi_s(T_sx)$. For the second we assume that $x\in
\A$ be selfadjoint. Let $s<t$. Then we deduce from the Cauchy
Schwarz inequality in the form $|\tau(abab)|\le \tau(abba)$ that
 \begin{align*}
 &\|m_s(x)-m_t(x)\|_4^4
  \lel  \|
  \pi_s(|T_sx|^2)+\pi_t(|T_tx|^2)-\pi_s(T_sx)\pi_t(T_tx)-\pi_t(T_tx)\pi_s(T_sx)\|_2^2\\
 &= \tau((\pi_s(|T_sx|^2)+\pi_t(|T_tx|^2))^2)+\tau((\pi_s(T_sx)\pi_t(T_tx)+
 \pi_t(T_tx)\pi_s(T_sx))^2)\\
 &\quad -2\tau((\pi_s(|T_sx|^2)+\pi_t(|T_tx|^2))
 (\pi_s(T_sx)\pi_t(T_tx)+\pi_t(T_tx)\pi_s(T_sx)) \\
 &=\tau((T_sx)^4)+\tau((T_tx)^4)+2\tau(\pi_s((T_sx)^2)\pi_t((T_tx)^2))\\
 &\quad
 +2\tau(\pi_s(T_sx)\pi_t(T_tx)\pi_s(T_sx)\pi_t(T_tx))+2\tau(\pi_s((T_sx)^2)
 \pi_t((T_tx)^2))\\
 &\quad  -2\bigg(\tau(\pi_s((T_sx)^3)\pi_t(T_tx))+
 \tau(\pi_t((T_tx)^3)\pi_s(T_sx))\bigg) \\
 &\kl \tau((T_sx)^4) +\tau((T_tx)^4)
 +6\tau(\pi_s((T_sx)^2)\pi_t((T_tx)^2)) \\
 &\quad  - 4\big(\tau(\pi_s((T_sx)^3)\pi_s(T_{t-s}T_tx))+\tau((T_tx)^4)\big)\\
  &=\tau( (T_sx)^4) +6\tau(T_{t-s}(T_sx)^2(T_tx)^2)
  -4\tau((T_sx)^3T_{2t-s}x) -3\tau((T_tx)^4) \\
  &= 3\bigg(\tau((T_tx)^4)-\tau((T_sx)^4))\bigg) +
  6\bigg(\tau(T_{t-s}(T_sx)^2(T_tx)^2) -\tau((T_tx)^4)\bigg)\\
  &\quad -4\bigg(\tau((T_sx)^3T_{2t-s}x)-\tau((T_sx)^4)\bigg) \pl.
  \end{align*}
Now, it suffices to consider the terms separately. Let
$d=T_tx-T_sx$. Then we use $|\tau(abab)|\le \tau(abba)$ again and
find
 \begin{align*}
 &|\tau( (T_sx+d)^4-(T_sx)^4)| \kl
 4|\tau(d^3T_s(x))+4|\tau(dT_s(x)^3|+ 2|\tau(dT_sxdT_sx)|+4\tau(d^2(T_sx)^2)\\
 &\kl 4|\tau(d^3T_sx)|+4  |\tau(dT_s(x)^3)|+6 \tau(d^2(T_sx)^2) \\
 &\kl 12 \|d\|_2 \max\{ \|(T_tx-T_sx)^2T_sx\|_2,\|(T_sx)^3\|_2,
 \|(T_tx-T_sx)(T_sx)\|_2\} \pl .
 \end{align*}
For $x$ in the domain of $A$, we know that $f(t)=T_tx$ is
differentiable and hence
 \begin{equation}\label{2ess}
  \|T_tx-T_sx\|_2 \lel \|\int_{s}^t T_rAx\|_2\kl (t-s)\|Ax\|_2
 \pl \end{equation}
Equation \eqref{2ess} also allows us to estimate the last term.
Indeed,  by Cauchy-Schwarz
 \[ |\tau((T_sx)^3T_{2t-s}x)-\tau((T_sx)^4)|
 \kl 2(t-s) \|AT_sx\|_2 \|T_sx\|_6^3
 \kl 2(t-s) \|Ax\|_2 \|x\|_N^3 \pl .\]
For the middle term we consider the function
 \[ f(r)\lel T_{t-r}(T_rx^*T_ry) \pl .\]
Due to our assumption this function is differentiable and
 \[ f'(r) \lel
 T_{t-r}A(T_rx^*T_ry)-T_{t-r}(AT_rx^*T_ry)-T_{t-r}(T_rx^*AT_{r}y)
 \lel -2T_{t-r}\Gamma(T_rx,T_ry) \pl .\]
This implies
 \begin{equation}\label{wies}
 T_{t-s}(|T_sx|^2)-|T_tx|^2
 \lel \int_s^t 2T_{t-r}\Gamma(T_rx,T_rx) dr  \pl.
 \end{equation}
If $\Gamma^2\gl 0$ we have
 \[ \|T_{t-s}(|T_sx|^2)-|T_tx|^2\|_2 \kl 2(t-s) \|\Gamma(x,x)\|
 \pl .\]
A similar estimate holds if just assume  $\sup_r
\|\Gamma(T_rx,T_rx)\|_2\kl C$, because $T_{t-r}$ is a contraction on
$L_2(N)$. This implies
 \[ \|m_t-m_s\|_4^4 \kl (t-s) (40\pl
 \|x\|_N^3\|Ax\|_2+12\pl \|\Gamma(x,x)\|\|x\|_N^2) \]
for all $x\in \A$. The noncommutative version of Kolmogorov's
theorem is proved in \cite{Golu}. Thus $m_t(x)$ has continuous path.
Due to Doob's inequality the class of martingales with continuous
path is closed in $L_p(N)$. Since $\A$ is assumed to be weakly dense
and hence norm dense in $L_p(N)$, we deduce the assertion for all
$x\in N$. For the last formula we observe for $x\in \A$ that
 \begin{align*}
  E_{[t}(|\pi_s(T_sx)-\pi_t(T_tx)|^2)&=
  E_{[t}(\pi_s(T_s|x|^2))-\pi_t(|T_tx|^2)
   \lel \pi_t( T_{t-s}|T_sx|^2)-|T_tx|^2) \pl .
  \end{align*}
We deduce from \eqref{wies} that
 \begin{equation}\label{kff}
  E_{[t}(\pi_s(T_s|x|^2)-\pi_t(|T_tx|^2))
 \lel \pi_t(\int_{s}^t 2T_{t-r}\Gamma(T_rx,T_rx) dr)
 \lel E_{[t} \int_s^t \pi_r(2\Gamma(T_rx,T_rx)) dr  \pl .
 \end{equation}
This implies for the  limit
 \begin{align*}
  &\langle m(x),m(x)\rangle_s- \langle m(x),m(x)\rangle_t
   \lel \lim_{|\pi|\to 0} \sum_{j} E_{[s_{j+1}}(
 |\pi_{s_j}(T_{s_j})-\pi_{s_{j+1}}(T_{s_{j+1}})|^2)\\
  &= \lim_{|\pi|\to 0} \sum_{j}
  E_{[s_{j+1}} \int_{s_j}^{s_{j+1}}
   2 \pi_r(\Gamma(T_rx,T_rx)) dr  \lel
   \int_s^t 2\pi_r(\Gamma(T_rx,T_rx))dr \pl .
  \end{align*}
Here the limit is taken over refining partitions. Since $m(x)$ has
continuous path, we know that the brackets coincide. In order to
remove the extra term $ E_{s_{j+1}}$ from the integral in the last
inequality, we first note that for $r<s$
 \begin{align}
 \|\pi_r(a)-E_{s}\pi_r(a)\|_2^2
 &=\|\pi_r(a)-\pi_s(T_{s-r}a)\|_2^2 \lel \tau(a^2-(T_{s-r}a)^2) \nonumber \\
 &\le \|T_{s-r}a-a\|_2(\|a\|_2+\|T_{s-r}a\|_2)
 \kl |s-r| \pl  \|Aa\|_2\pl (\|a\|_2+\|T_{s-r}a\|_2) \label{l222}
 \pl .
 \end{align}
Applying this to $a=\Gamma(T_rx,T_rx)$, we are done assuming the
continuity of $f(r,s)=T_s\Gamma(T_rx,T_rx)$ in $L_2$. In fact, we
deduce from the Cauchy Schwarz inequality that
 \begin{align*}
  \|A\Gamma(T_rx,T_rx)\|_2&\le 2 \pl \|\Gamma^2(T_rx,T_rx)\|_2
  +\|\Gamma(AT_rx,T_rx)\|_2+\|\Gamma(T_rx,AT_rx)\|_2\\
  &\le 2 \pl  \|\Gamma^2(T_rx,T_rx)\|_2+ 2\pl
  \|\Gamma(T_rAx,T_rAx)\|_2^{\frac12 }
  \|\Gamma(T_rx,T_rx)\|_2^{\frac12}\\
  &\le  2 \pl \|\Gamma^2(T_rx,T_rx)\|_2 +
  2 \pl \|\Gamma(Ax,Ax)\|_2^{\frac12}\|\Gamma(x,x)\|_2^{\frac12}\pl .
  \end{align*}
If $\Gamma^3\gl 0$, we have $\|\Gamma^2(T_rx,T_rx)\|\le
\|\Gamma^3(x,x)\|_2$ and $\Gamma^3(x,x)\in N$ for $x\in \A$.
Assuming only the boundedness of $\Gamma^2(T_rx,T_rx)$, still allows
us to obtain the assertion by refining the partition. \qd

\begin{rem}{\rm Under the condition ii) we see moreover that
 \[ \langle m(x),m(x)\rangle_s-\langle m(x),m(x)\rangle_t
 \kl 2\pl  |s-t| \pl  \|\Gamma(x,x)\|_{\infty} \]
holds for $s<t$. By approximation this implies that the bracket is
absolutely continuous with respect to the Lebesgue measure. In the
commutative case this is enough to show that the martingale can be
obtained as a stochastic integral against the brownian motion. We
refer to \cite{JK} for similar  applications of this regularity
condition.}
\end{rem}

\begin{prop}\label{shift} Let $2\le p<\infty$ and $\Gamma^2\gl 0$. Let $x\in \N$ be a selfadjoint mean $0$ element. Then
 \[ c_{p}^{-1}\|x\|_{p} \kl \|(\int_0^{\infty} \Gamma(T_tx,T_tx) dt)^{\frac12 }\|_p
 \kl c_p \|x\|_p \pl .\]
Moreover, for every $x\in N$ there exists a martingale $m^2(x)$ such
that
 \[ \|m^2(x)\|_{H_p^c} \kl c(p) \|\int_0^{\infty}
 \Gamma(T_sx,T_sx) ds \|_{\frac p2}^{\frac12} \quad \mbox{and} \quad
  \tau(\pi_0(y)^*m^2(x)) \lel \frac{1}{3} \tau(y^*(I-\pr)x) \pl .\]
Here $\pr$ is the projection onto the kernel of $A$. If in addition
$f(r,s)=T_r\Gamma(T_sx,T_sx)$ is $L_2$-continuous for $x\in \A$,
then
 \[ \langle m^2(x),m^2(x)\rangle \lel \int_0^{\infty}
 \pi_s(\Gamma(T_{2s}x,T_{2s}x)) ds  \pl. \]
\end{prop}

\begin{proof} Let $s_0$ be fixed and $\si=\{0,...,s_0\}$ a
partition. We define the martingale differences
 \[ d_j \lel
 (\pi_{s_j}(T_{s_j+s_{j+1}}x)-\pi_{s_{j+1}}(T_{2s_{j+1}})) \pl. \]
Indeed,
$E_{s_{j+1}}\pi_{s_j}(T_{s_j+s_{j+1}}x)=\pi_{s_{j+1}}(T_{2s_{j+1}}x)$
shows that
 \[ m_{\si} \lel \sum_{j=1}^{s_0} d_j
 \pl .\]
is a martingale with respect to the discrete filtration
$(N_{[s_j})$. Following Lemma \ref{ctpath} (in particular
\eqref{kff}), we may also calculate the conditioned bracket
 \[ \sum_{j} E_{[s_{j+1}}(|d_j|^2)
 \lel \sum_j E_{s_{j+1}}\int_{s_j}^{s_{j+1}}
 \pi_r(\Gamma(T_{r+s_{j+1}}x,T_{r+s_{j+1}}x)) dr \pl .\]
On the other hand, we deduce again from Lemma \eqref{ctpath} that
 \begin{align*}
  \|d_j\|_4^4 &\lel  \|m_{s_j}(T_{s_{j+1}}x)-m_{s_{j+1}}(T_{s_{j+1}}x)
 \|_4^4 \\
 & \kl C(t-s)  (\|T_{s_{j+1}}x\|_N^3 \|AT_{s_{j+1}}x\|_2 +
 \|T_{s_{j+1}}x\|_N^2\|\Gamma(T_{s_{j+1}}x,T_{s_{j+1}}x)\|)
   \pl. \end{align*}
This allows us to define the weak$^*$ limit
 \[ m_{s_0} \lel \lim_{\si} m_{\si} \]
as a martingale in $L_4(N)$ with continuous path such that
 \begin{align*}
  &\lim_{\si} \|m_{\si}\|_{H_p^c(\si)}^2
  \kl c_{p/2} \lim_{\si} \|\sum_j \int_{s_j}^{s_{j+1}}
  \pi_r(\Gamma(T_rT_{s_{j+1}}x,T_rT_{s_{j+1}}x)) \|_{p/2}\\
  &\le c_{p/2} \lim_{\si} \|\sum_j \int_{s_j}^{s_{j+1}} \!\!\!
  \pi_r(T_r\Gamma(T_{s_{j+1}}x,T_{s_{j+1}}x)) \|_{p/2}\\
  &\lel
 c_{p/2} \lim_{\si} \|\sum_j \int_{s_j}^{s_{j+1}} \!\!\!
  E_{[r}\pi_0(\Gamma(T_{s_{j+1}}x,T_{s_{j+1}}x)) \|_{p/2} \\
  &\le c_{p/2}^2 \lim_{\si} \|\sum_j
  \int_{s_j}^{s_{j+1}}\Gamma(T_{s_{j+1}}x,T_{s_{j+1}}x) dr
  \|_{p/2}
  \lel c_{p/2}^2 \|\int_{0}^{s_0} \Gamma(T_sx,T_sx) ds\|_{p/2} \pl
  .
  \end{align*}
With the results from \cite{JK}, this implies that $y_{s_0}\in
H_p^c$ for all $p<\infty$. Using the same estimate for $x^*$, we
deduce that $m_{s_0}\in L_p$ for all $p<\infty$. We may take another
weak$^*$ limit to define $m^2(x)=\lim_{s_0\to \infty} m_{s_0}$ which
has continuous path because of the uniform estimate for the norm
$\|x\|_{h_4^d}$. Moreover, the proof of
 \begin{equation} \label{brry}
  \langle m^2(x),m^2(x)\rangle \lel 2\int_0^{\infty}
 \pi_s(\Gamma(T_{2s}x,T_{2s}x)) ds \pl .
 \end{equation}
is the same as in Lemma \ref{ctpath}. Even without knowing exactly
what the bracket looks like, it is easy to complete the estimate for
selfadjoint $x$. We assume that $x$ has mean $0$, i.e. $\pr(x)=x$.
Let $y\in \A$ such that $\|y\|_{p'}=1$ and
 \[ \|x\|_p \kl (1+\delta) |\tau(y^*x)| \pl .\]
Let us consider
 \begin{align*}
  &\tau( \pi_0(y)^*m_{\si})
   \lel \sum_j  \tau\big((\pi_{s_{j}}(T_{s_{j}}(y^*))-
  \pi_{s_{j+1}}(T_{s_{j+1}}(y^*)))d_j\big)\\
  &= 2\sum_j
  \int_{s_j}^{s_{j+1}}\tau\bigg(\pi_r\big(\Gamma(T_ry,T_{s_{j+1}+r}x)\big)\bigg)
  dr
  \lel 2
   \sum_j
  \int_{s_j}^{s_{j+1}}\tau\bigg(\Gamma(T_ry,T_{s_{j+1}+r}x)\bigg) dr \pl .
  \end{align*}
Thus in the limit we obtain
 \begin{align*}
 \lim_{s_0}\lim_{\si} \tau( \pi_0(y)^*m_{\si})
  \lel 2\int_0^{\infty} \tau(\Gamma(T_ry,T_{2r}x)) dr
  \lel 2 \int_0^{\infty} \tau(y^*AT_{3r}x) dr \pl .
  \end{align*}
Using the spectral resolution for $A=\int_0^{\infty} \la dE(\la)$
and $d\nu_x(\la)=(x,dE(\la)x)$,  we see that
 \begin{align*}
 \int_0^{\infty} (x,AT_{3r}x) dr
 &=\int_0^{\infty} \int_0^{\infty} \la e^{-3r\la} dr d\nu_x(\la)
  \lel  \frac{1}{3} (E_{>0}x,E_{>0}x) \pl .
 \end{align*}
Thus, $\int_0^{\infty} AT_{3r}x=\frac{1}{3}E_{>0}x$, where
$E_{>0}=I-\pr$ is the orthogonal projection onto the complement of
of  the kernel of $A$. Since for selfadjoint $x$ the martingale
$m^{2}(x)$ is also selfadjoint, we deduce
 \begin{align*}
  &\|x\|_p \kl 3(1+\delta) \lim_{s_0}\lim_{\si}
  |\tau(\pi_0(y^*)m_{\si})| \lel 3(1+\delta) |\tau(\pi_0(y)^*m^2(x))| \\
  &\le c_p 3(1+\delta) \|y\|_{p'} \|m_2(x)\|_{H_p^c} \kl
   c_pc_{p/2}3(1+\delta)  \lim_{s_0}\lim_{\si}
 \|\sum_j \int_{s_j}^{s_{j+1}}  \Gamma(T_{s_{j+1}}x,T_{s_{j+1}}x)
 dr \|_{p/2}^{1/2} \\
 &=  c_pc_{p/2}3(1+\delta)  \pl \|(\int_0^{\infty} \Gamma(T_sx,T_sx) ds
 )^{\frac12 } \|_p \pl .
 \end{align*}
For the upper estimate we refer to \cite{JR1} which applies due to
$H^{\infty}$-calculus. \qd

\begin{lemma}\label{cond} Let $(T_t)$ and $\A$ be as above. Assume that there
is a further  von Neumann algebra $M$, a sequence $(u_j):\A\to \M$
such that
 \[ \Gamma(x,x) \lel \sum_j u_j(x)^*u_j(x) \]
and semigroup $\hat{T}_t$ with cb-$H^{\infty}$ calculus such that
 \[ (u_j(T_tx)) \lel (\hat{T}_tu_j(x)) \pl. \]
Then
 \[ \|\int_0^{\infty}\Gamma(A^{1/2}T_tx,A^{1/2}T_tx) dt\|_{p/2}
 \kl c_p^2 \pl \|\Gamma(x,x)\|_{p/2} \]
holds for all $2\le p<\infty$ and all mean $0$ elements $x\in \A$.
\end{lemma}

\begin{proof} Since $id_{B(\ell_2)}\ten \hat{T}_t$ satisfies
$H^{\infty}$-calculus, we deduce from \cite{JLX} that
 \begin{align}
 &\|\int_0^{\infty} \Gamma(A^{1/2}T_tx,A^{1/2}T_tx) dt \|_{p/2}^{1/2}
  \lel \| (\int_0^{\infty}  |\hat{A}^{1/2}\hat{T}_t(\sum_j e_{j,1} \ten u_j(x))|^2 dt)^{1/2}  \|_p  \label{pppp} \\
 &\le c_p \|\sum_j e_{j,1} \ten u_j(x)\|_{L_p(B(\ell_2)\ten M)}
 \lel
  c_p \| \Gamma(x,x)^{\frac12}\|_p \pl . \nonumber
 \end{align}
Here $\hat{A}$ is the generator of $\hat{T}_t$. Note also that
$u_j(T_tx)=\hat{T}_tu_j(x)$ implies $u_j(P_tx)\lel \hat{P}_tu_j(x)$
according to \eqref{ps-form} and hence by differentiation
$\hat{A}^{1/2}u_j(T_tx)=u_j(A^{1/2}T_tx)$. This justifies the first
equation in \eqref{pppp} and completes the proof. \qd

Let us recall the notation
$\Gamma_{\al_1,\al_2}=(\Gamma_{A^{\al_1}})_{A^{\al_2}}$ for iterated
gradients, see \cite{JR1}.

\begin{prop}\label{pps} Let $(T_t)$ be a semigroup with a Markov dilation and
$\Gamma^2\gl 0$. Let $(P_t)$ be the associated Poisson semigroup
satisfying
 \[ A^{1/2}\A\subset \A \quad, \quad P_s(\A)\subset \A \]
and that $f(s,t)=P_s\Gamma_{1/2,1}(P_tx,P_tx)$ is $L_2$-continuous
for $x\in \A$. Then
 \[ \|(\int_0^{\infty} \Gamma(A^{1/2}P_sx,A^{1/2}P_sx) sds)^{\frac12} \|_{p}
 \kl c_p \|\Gamma(x,x)^{\frac12}\|_p \pl .\]
\end{prop}

\begin{proof} A glance at \eqref{ps-form} shows a Markov
dilation for $T_t$ implies that $P_t$ is factorable. According to
\cite{Cl}, we deduce that $P_t$ also has a Markov dilation. Let us
denote this family of maps again with $\pi_s$. We consider the
submartingale
 \[ y_s\lel \pi_s(\Gamma(P_sx,P_sx)) \pl .\]
Indeed, for $s<t$ we deduce from $\Gamma^2\gl 0$ that
 \[ E_{[t}(y_s)
 \lel \pi_t(P_{t-s}\Gamma(P_sx,P_sx))
 \gl \pi_t( \Gamma(P_tx,P_tx)) \pl. \]
As in Lemma \ref{ctpath} we  consider $f(r)=
P_{t-r}\Gamma(P_rx,P_rx)$ and obtain
  \begin{align*} f'(r) &=
 A^{1/2}P_{t-r}\Gamma(P_rx,P_rx)-P_{t-r}\Gamma(A^{1/2}P_rx,P_rx)-
 P_{t-r}\Gamma(P_rx,A^{1/2}P_rx)\\
 &= -2P_{t-r}\Gamma_{1/2,1}(P_rx,P_rx) \pl .
 \end{align*}
This implies
 \[ E_{[t}(y_s-y_t)
 \lel 2E_{[t}\int_{s}^t \pi_r(\Gamma_{1/2,1}(P_rx,P_rx)) dr \pl
 .\]
We apply \eqref{l222} for $P_s$ and obtain
 \[ \|\pi_r(y)-\pi_t(P_{t-r}y)\|_2^2
 \kl \|y\|_2 \|P_r(y)-P_t(y)\|_2 \pl .\]
We deduce that
 \begin{align*}
 & \|\sum_{j} \int_{s_j}^{s_{j+1}}
  \pi_r(\Gamma_{1/2,1}(P_rx,P_rx))-E_{s_{j+1}}(\pi_r(\Gamma_{1/2,1}(P_rx,P_rx))
  dr \| \\
  &\le \sum_{j} \int_{s_{j+1}}^{s_j}
  \|\Gamma_{1/2,1}(P_rx,P_rx)\|_2^{1/2}
  \|(P_r-P_{s_{j+1}})\Gamma_{1/2,1}(P_rx,P_rx)\|_2^{1/2} dr \pl .
  \end{align*}
Thus uniform continuity of $P_t\Gamma_{1/2,1}(P_rx,P_rx)$ implies
that the limit converges to $0$ as long as the mesh size converges
to $0$. Therefore we obtain
 \[ \lim_{|\si|\to 0} \sum_{j} E_{s_{j+1}}(y_{s_j}-y_{s_{j+1}})
 \lel 2\int_{s}^t \pi_r(\Gamma_{1/2,1}(P_rx,P_rx)) dr \pl .\]
Now, we apply the inequality for potentials Lemma \ref{potential}
and find
 \[ \|\sum_{j} E_{s_{j+1}}(y_{s_j}-y_{s_{j+1}})\|_{p/2} \kl c_{p/2}
 \|y_0\|_{p/2} \pl .\]
Indeed, since we are working with a reversed martingale,
$y_0=\pi_0(\Gamma(x,x))$ is the endpoint. Passing to the limit, we
deduce that
 \[ \|\int_0^{\infty} \pi_s(\Gamma_{1/2,1}(P_sx,P_sx)) ds \|_{p/2}
 \kl c_{p/2} \|y_0\|_{p/2} \lel c_p \|\Gamma(x,x)\|_{p/2} \pl . \]
Conditioning on $\pi_0$ yields
 \[ \|\int_0^{\infty} P_s\Gamma_{1/2,1}(P_sx,P_sx) ds\|_{p/2}
 \kl c_{p/2} \|\Gamma(x,x)\|_{p/2} \pl .\]
Now, we recall from \cite{JR1} that
 \begin{align*}
 & \Gamma_{1/2,1}(y,y)
 \lel \int_{0}^{\infty} P_t\Gamma(A^{1/2}P_ty,A^{1/2}P_ty) dt
 + \int_{0}^{\infty} \Gamma^2(P_ty,P_ty) dt\\
 & \gl \int_{0}^{\infty} P_t\Gamma(A^{1/2}P_ty,A^{1/2}P_ty) dt \pl
 . \end{align*}
Here we use $\Gamma^2\gl 0$. Therefore, we obtain
 \begin{align*}
  &\|\int_0^{\infty} \Gamma(A^{1/2}P_sx,A^{1/2}P_sx) sds \|_{p/2}
  \lel 2 \pl \|\int_0^{\infty} \Gamma(A^{1/2}P_{2s}x,A^{1/2}P_{2s}x) 2sds
  \|_{p/2}\\
  &\le 2\pl \|\int_0^{\infty} P_s\Gamma(A^{1/2}P_{s}x,A^{1/2}P_{s}x) 2sds
  \|_{p/2} \\
  &=
   2 \pl \|\int_0^{\infty}\int_0^{\infty}
  P_{s+t}\Gamma(A^{1/2}P_{s+t}x,A^{1/2}P_{s+t}x) ds dt \|_{p/2} \\
  &\kl 2 \pl \|\int_0^{\infty} P_s\Gamma_{1/2,1}(P_sx,P_sx)) ds
  \|_{p/2}\kl
 2 c_{p/2}  \pl \|\Gamma(x,x)\|_{p/2} \pl . \qedhere
  \end{align*}
 \qd

Our next task is to replace $P_s$ by $T_s$ following a well-known
path in \cite{JLX}.

\begin{lemma}\label{cond3} Let $2\le p<\infty$, $\tilde{\Gamma}$ be a completely positive  form on $\bar{\A}\times \A$,
and  $(T_t)$ be a semigroup of selfadjoint maps with selfadjoint
generator such that
  \begin{equation}\label{assanag2}
 \|(\sum_k \tilde{\Gamma}(T_{z_k}x_k,T_{z_k}x_k))^{1/2} \|_p
  \kl c(p,\theta) \|(\sum_k \tilde{\Gamma}(x_k,x_k))^{1/2}\|_p \pl\end{equation}
for all $z_k$ with $0\le Arg(z_k)\le \theta$, where $0<\theta< \pi$.
Moreover, assume that $A^{1/2}L_2(N)$ is dense in $(I-\pr)L_2(N)$
with respect to $\|x\|_{\tilde{\Gamma}}=\tau(\tilde{\Gamma}(x,x))$.
Then
 \begin{equation}\label{anag2}
  \|(\int_0^{\infty} \tilde{\Gamma}(T_sx,T_sx) ds)^{\frac12}\|_p
 \kl c_0 c(p,\theta) \|(\int_0^{\infty}\tilde{\Gamma}(P_sx,P_sx)
 sds)^{\frac12}\|_p \pl \end{equation}
holds for all $x\in \A$ with $\pr(x)=0$.
\end{lemma}

\begin{proof} We introduce the space $L_p(L_2^c(\tilde{\Gamma}))$ as the the
closure of continuous functions such that
 \[ \|f\|_{L_p(L_2^c(\tilde{\Gamma}))}
 \lel  \| (\int_0^{\infty} \tilde{\Gamma}(f(s),f(s)) \pl \frac{ds}{s})^{1/2}
 \|_p  \pl .\]
As in \cite[Corollary 4.9]{JLX} our assumption implies that the
family $z(z-a)^{-1}$ is Col-bounded on $L_p(L_2^c(\Gamma))$ for the
same angle. Then we may apply \cite[Theorem 4.14]{JLX} and deduce
that $T_{\Phi}$ with the kernel $\Phi(s,t)=F_2(sA)F_1(tA)$ is
bounded on $L_p(L_2(\tilde{\Gamma}))$. We may choose
$F_2(z)=z^{1/2}e^{-z}$ and $F_1(z)=ze^{-z}$. Let us assume that
$x=A^{1/2}y$. Let us define $f(t)=\sqrt{At}P_{\sqrt{t}}y$ Using a
change of variable, we deduce that
 \begin{align*}
  \int_0^{\infty} \tilde{\Gamma}(f(t),f(t)) \pl \frac{dt}{t}
  &= \int_0^{\infty} \tilde{\Gamma}(P_{\sqrt{t}}x,P_{\sqrt{t}}x) t
 \pl   \frac{dt}{t} \lel \frac12  \int_0^{\infty}
 \tilde{\Gamma}(P_{s}x,P_{s}x) \pl s^2
  \frac{ds}{s}  \pl .
  \end{align*}
In order to apply $T_{\Phi}$ we have to calculate
 \begin{align*}
  \int_0^{\infty} F_1(tA)f(t) \pl \frac{dt}{t}
  &= \int_0^{\infty} tAT_t(t^{1/2}A^{1/2}P_{\sqrt{t}}y) \pl \frac{dt}{t}
  \lel \int_0^{\infty} t^{\frac32}A^{\frac32} T_tP_{\sqrt{t}}y  \pl
  \frac{dt}{t} \pl .
  \end{align*}
Let $dE_{\la}$ be the spectral measure for $A$ and $d\mu_{y_1,y}$
the induced probability measure for elements $y_1,y\in L_2(N)$. Then
 \begin{align*}
 \tau(y_1^*\int_0^{\infty} t^{\frac32}A^{\frac32} T_tP_{\sqrt{t}}y  \pl
  \frac{dt}{t})
  &= \int_0^{\infty} \int_0^{\infty} t^{\frac32}\la^{\frac32}
  e^{-t\la}e^{-\sqrt{t}\sqrt{\la}} \frac{dt}{t}
  d\mu_{y_1,y}(\la)\\
 &= \kla \int_0^{\infty} t^{\frac32}
  e^{-t}e^{-\sqrt{t}} \frac{dt}{t} \mer \tau(y_1^*y_2) \pl .
  \end{align*}
Let us denote by $c$ the constant given by the converging integral.
Since $y_1$ is arbitrary we deduce by approximation that
 \[ T_{\Phi}(f)(s) \lel c A^{1/2}T_s(y) \lel c T_s(x) \pl \]
holds for all $x\in L_2(L_2(\tilde{\Gamma}))$ with $\pr(x)=0$. Let
us note that the assumption that $A$ is selfadjoint is not
necessary. In the sectorial case we refer to \cite[Lemma 6.5]{JLX}
and the argument in \cite[Theorem 6.7]{JLX}. The assertion follows
from the boundedness of $T_{\Phi}$ on $L_p(L_2(\tilde{\Gamma}))$.
\qd

\begin{rem}\label{hppp} {\rm For a semigroup of (selfadjoint)  completely
positive maps and the canonical form  $\tilde{\Gamma}(x,y)=x^*y$, we
deduce that
 \[ \|x\|_{H_p^c(T)} \kl c(p)  \|x\|_{H_p^c(P)} \pl \]
without assuming $H^{\infty}$-calculus. The reverse inequality is
shown in \cite{JR1} and hence the equivalence of different semigroup
$H_p$-norms holds without using a Markov dilation.}
\end{rem}

\begin{rem}\label{cond2}{\rm 
Assume the assumption of Lemma \ref{cond} is satisfied. Then
\eqref{assanag2} holds for $\tilde{\Gamma}=\Gamma$,
$\theta<\pi(\frac12-\frac1p)$ and \eqref{anag2} holds also. }
\end{rem}

\begin{proof} For a selfadjoint semigroup the results in \cite{JX} imply
 \[ \|\sum_k |\hat{T}_{t_k}y_k|^2 \|_{p/2}
 \kl \|\sum_k \hat{T}_{t_k}|y_k|^2 \|_{p/2}
 \kl c_p \|\sum_k |y_k|^2 \|_{p/2} \pl .\]
For $p=2$ we have
 \[ \|\sum_k |\hat{T}_{z_k}y_k|^2\|_1\lel \sum_k \|\hat{T}_{z_k}y_k\|_2^2
 \kl \sum_k \|y_k\|_2^2  \pl \]
whenever $Re(z_k)\gl 0$. Then the assertion is a standard
application of Stein's theorem on interpolation of analytic families
applied to $y_k=u(x_k)$ and yields  \eqref{assanag2}. Moreover, we
may directly apply the argument in Lemma \ref{cond3} for $\hat{T}_t$
and $(y_j)$. Note that $x$ orthogonal to the kernel of $A$ implies
that $(y_j)$ is orthogonal to the kernel of $\hat{A}$. Thus we
obtain \eqref{anag2} without using the extra density assumption.\qd

The following argument is based on a continuous version of a result
of Stein. 

\begin{theorem}\label{stein3} Let   $1<p<\infty$.
\begin{enumerate}
\item[i)] Let $m \lel \int_0^{\infty} dm_s$ be a martingale with
continuous path in $H_p$. Then
 \[ \|\kla \int_0^{\infty} |\frac{1}{r}\int_0^r s\pl dm_s|^2 \frac{dr}{r}
 \mer^{1/2}\|_p
 \kl c_p \|m\|_{H_p^c} \pl .\]
\item[ii)] Let $(T_t)$ be semigroup with a martingale dilation  and
 $x\in \A$. Then
 \[  \| x\|_{H_p^c(T)} \kl c_p \|\pi_0(x)\|_{H_p^c} \pl .\]
\item[iii)] Let in addition $2\le p<\infty$ and $\Gamma^2\gl 0$. Then
 \[ c_p^{-1} \|x\|_{H_p^c(T)} \kl \|\kla \int_0^{\infty}
 \Gamma(T_sx,T_sx)ds \mer^{\frac12} \|_p \kl C_p \|x\|_{H_p^c(T)} \pl .\]
\end{enumerate}
\end{theorem}

\begin{proof} The bulk of the argument is due to Stein, the noncommutative part of the argument is contained in
\cite[Proposition 10.8]{JLX} where it is proved that for a
martingale difference sequence $(d_j)$,
 \[ E_k\lel \sum_{j=0}^k d_j \quad, \quad  \La_m(x)\lel \frac{1}{m+1}\sum_{k=0}^m E_k(x) \]
and $\Delta_m(x)=\La_m(x)-\La_{m-1}(x)$ one has
 \begin{align}
 \|(\sqrt{m}\Delta_m(x))\|_{L_p(\ell_2^c)}
 \kl C_p  \| (d_j)_j\|_{L_p(\ell_2^c)}
 + C_p \|(\sum_{j=2^{k}+1}^{2^{k+1}}d_j)_k \|_{L_p(\ell_2^c)} \pl .
 \label{stein}
 \end{align}
Here $C_p$ is the norm of Stein's projection in $L_p$. It is also
important to note that this argument is true for decreasing or
increasing martingale differences. Moreover, let $(m_t)$ be a
continuous martingale so that
 \[ \lim_{|\si| \to 0} \| (\sum_{t_j} |m_{t_{j+1}}-m_{t_j}|^2)^{1/2} \|_p
 \lel  \|m\|_{H_p^c} \pl .\]
Here the limit is taken over partitions of a fixed interval
$[\al,\beta]$ and the mesh size of the partitions converges to $0$.
Then we note that the right hand side of \eqref{stein} is controlled
by two partitions and hence
 \[ \lim_{|\pi| \to 0}
 \|(\sqrt{l}\Delta_l(m))\|_{L_p(\ell_2^c)}
 \kl C_p \|m\|_{H_p^c} \pl .\]
Let us now fix $0<\al<\beta<\infty$ and assume that
 \[ E_k(x) \lel \int_{\al}^{\al+\frac{k}{n}} dx_s  \pl \]
holds in terms of stochastic integrals. This implies
  \begin{align*}
 &\La_l(m)-\La_{l-1}(m)
 \lel \frac{1}{l+1} E_l(m)+ \sum_{k=1}^{l-1}
 E_k(m)(\frac{1}{l+1}-\frac{1}{l}) \\
 &= \frac{1}{l+1} \int_{\al+\frac{l-1}{n}}^{\al+\frac{l}{n}} dm_s
 + \frac{1}{l+1} \int_{\al}^{\al+\frac{l-1}{n}}
 (1-(\frac{l-\lceil n(s-\al) \rceil}{l})) dm_s \\
 &= \frac{1}{l+1} \int_{\al+\frac{l-1}{n}}^{\al+\frac{l}{n}} dm_s
 + \frac{1}{l(l+1)} \int_{\al}^{\al+\frac{l-1}{n}}
  \lceil n(s-\al) \rceil  dm_s  \pl .
 \end{align*}
Here $\lceil x \rceil$ is the smallest integer $\gl x$. The first
part is easy to control and hence
 \[  \| \sum_{l=2}^{\frac{\al-\beta}{n}}
 \frac{1}{l-1} |\int_{\al}^{\al+\frac{l-1}{n}}
  \frac{\lceil n(s-\al)  \rceil}{l-1} dm_s|^2 \|_{p/2}^{1/2}
  \kl c_p \|x\|_{H_p^c} \pl .\]
Passing to the limit (see \cite[]{JLX} for a similar reasoning),  we
deduce that
 \begin{align*}
  \|\kla \int_{\al}^{\beta} |\int_{\al}^r \frac{(s-\al)}{r} dm_s|^2
 \frac{dr}{r} \mer^{\frac12} \|_p
 \kl c_p \|x\|_{H_p^c} \pl
  \end{align*}
provided the square function  is Riemann integrable. Finally, we may
sent $\al$ and $\beta$ to the boundary and obtain
 \[   \|\kla \int_{0}^{\infty} |\frac{1}{r} \int_{0}^r sdm_s|^2
  \frac{dr}{r} \mer^{\frac12} \|_p
 \kl c_p \|m\|_{H_p^c}  \pl .\]
This completes the proof of i). We apply this inequality first in
the particular case where $m=\pi_0(x)=\int_0^{\infty} dm_s(x)$ is
the reversed martingale decomposition. We also add the conditional
expectation $E_0$. Let $\pi_0(y)=\int_0^{\infty} dm_s(y)$ an other
element. Then we deduce from the calculus of brackets that
 \begin{align*}
 &\tau(y^*\int_{0}^r \frac{s}{r} dm_s(x))
 \lel 2\int_0^r \tau( \pi_s(\Gamma(T_sy,T_sx)\frac{s}{r})) ds \\
 &= 2\int_0^r \tau( y^*AT_{2s}x) \frac{s}{r})) ds
 \lel 2 \tau(y^* \int_0^{r}AT_{2s}x) \frac{s}{r} ds
 \lel \tau(y^* \int_0^1 AT_{2sr}x 2sr \pl ds)
 \pl.
 \end{align*}
Let us define
 \[ f(z) \lel \int_0^1 (2zs) e^{-2sz} ds
 \lel \frac{1-e^{-2z}}{2z}-e^{-2z} \pl.  \]
Note that $f(0)=0$ and vanishes at $\infty$.  Then we see that
 \[ f(rA)x \lel \int_0^1 2rsAT_{2rs}x ds \pl .\]
Using the equivalence of different square functions \cite{JLX}, we
deduce that
 \[ \|x\|_{H_p^c(T)} \kl c_p  \pl \|\pi_0(x)\|_{H_p^c} \pl .\]
For our last assertion we consider $p\gl 2$ and the martingale
$m=m^2(x)$. For fixed $r$, we deduce from the fact that the function
$f(s)=s$ is adapted that
 \[ \int_{0}^{t}  s \pl dm^2_s(x)  \]
is a martingale and hence
 \[ \langle \int_{0}^{\infty}  dm_s(y),\int_0^{r} dm^2_s(x) \p s\rangle
 \lel 2 \int_0^r \pi_s(\Gamma(T_sy,T_{2s}x))s ds \pl . \]
Here we use the projection on $\pi_0(N)$ from Lemma \ref{shift} and
the calculus of brackets for stochastic integrals. This implies
 \[ E_0(\frac{1}{r} \int_0^r s  dm^2_s(x))
  \lel  \frac{2}{r} \int_0^r AT_{3s}x \pl s ds
 \lel \frac{2}{3}  \int_0^1 AT_{3sr}x  \pl (3sr) ds \pl .\]
Thus replacing $f$ by $\tilde{f}(z)=\int_0^1 (3sz) e^{-3sz} ds$, we
deduce again with the equivalence of different square functions that
  \[ \| x\|_{H_p^c(T)} \kl c_p \|m^2(x)\|_{H_p^c}
  \kl c_p^2 \pl \|\kla \int_0^{\infty} \Gamma(T_sx,T_sx) ds \mer^{\frac12}
 \|_p \pl .\]
The last estimate is of course taken from Proposition \ref{shift}.
Assuming $\Gamma^2\gl 0$, we can refer to \cite{JR1} for the
estimate
 \begin{align*}
   \|\kla \int_0^{\infty} \Gamma(T_sx,T_sx) ds \mer^{\frac12}
  \|_p &\le  c_p \|x\|_{H_p^c(T)} \pl . \qedhere \end{align*}
\qd

\begin{theorem} Let $(T_t)$ be a semigroup satisfying
$\Gamma^2\gl 0$. Assume that
 \begin{enumerate}
 \item[i)]  The assumption of Lemma \ref{cond} is satisfied \emph{or}
 \item[ii)] Condition \eqref{anag2} is satisfied for $\widetilde \Gamma=\Gamma$ and $f(t,s)=P_s\Gamma_{1/2,1}(P_tx,P_tx)$ is
$L_2$-continuous for $x\in \A$.
\end{enumerate}
Then
  \[  \|A^{\frac12}x\|_{H_p^c(T)} \kl c_p\pl  \|\Gamma(x,x)^{\frac12}\|_p
  \]
and
 \[ \|A^{\frac12}x\|_p\kl c_p \max\{\|\Gamma(x,x)^{\frac12}\|_p
  ,\|\Gamma(x^*,x^*)^{\frac12}\|_p\}\pl .\]
holds for all mean $0$ elements $x\in \A$.
\end{theorem}

\begin{proof} For the first assertion we combine Theorem
\ref{stein3}iii)  with  Lemma \ref{cond} or Remark \ref{cond2}. This
allows us to apply Proposition \eqref{pps}. For the second
assertion, we refer to \cite{JLX} for the fact that a Markov
dilation implies $H^{\infty}$-calculus and hence
 \[ \|A^{1/2}x\|_p\sim_{c_p} \|A^{1/2}x\|_{H_p(T)}
 \lel \max\{\|A^{1/2}x\|_{H_p^c(T)},\|A^{1/2}x^*\|_{H_p^c(T)}\} \pl .\]
This immediately implies the second assertion. \qd

\begin{rem}{\rm In a forthcoming publication we will show that the
assumption in Lemma \ref{cond} is satisfied for Fourier multipliers
on discrete groups. Motivated by  the recent work of Ricard we will
also construct a Markov dilation satisfying the conditions i) and
ii) at the beginning of this section.}
\end{rem}

\section{The probabilistic model}

The probabilistic approach to Littlewood-Paley theory goes back to
the work of P.A. Meyer and has found many applications. Instead of
adding a time component to the manifold as in Stein's approach, the
probabilistic approach adds an additional brownian motion to the
picture. As in the previous section we assume that $(T_t)$ is a
semigroup of completely positive maps and that $\pi_s$ is a Markov
dilation satisfying i). (The reversed condition is no longer
necessary). Let us keep the notation $\pr_0$ for the projection on
the kernel of $A$ and start with a simple observation, well-known in
the commutative case.

\begin{lemma}\label{Mey} Let $-A$ be the negative generator of $T_t=e^{-tA}$ and
$x\in \dom(A)$ with $\pr_0(x)=0$. Then
\[ m_t(x)\lel \pi_t(x)+\int_0^{t} \pi_s(Ax) ds \]
is a martingale. \end{lemma}

\begin{proof} For $p>0$ we calculate
 \begin{align*}
 &E_s(\int_0^{\infty} e^{-pt}\pi_t((p+A)x) dt) \lel
  \int_0^{s}e^{-pt} \pi_t((p+A)x) dt + \int_{s}^{\infty}
 (p+A)e^{-pt}\pi_s(T_{t-s}(p+A)x)dt \\
 &= \int_0^{s}e^{-pt} \pi_t((p+A)x) dt +e^{-ps}
 \pi_s(\int_0^{\infty}e^{-t(p+A)}(p+A)x dt) \pl .
 \end{align*}
A change of variables shows that
$\int_0^{\infty}e^{-t(p+\la)}(p+\la) dt=1$ holds for every $\la\in
\rz$. Thus (arguing in $L_2$ if necessary), we see that
$\int_0^{\infty}e^{-t(p+A)}(p+A)x dt=x$ for all $x$ with
$\pr_0(x)=0$. Hence
 \[ m_{s,p}\lel e^{-ps}\pi_s(x)+\int_0^se^{-pt} \pi_t(px+Ax)dt  \]
is a martingale for  all $p>0$. Sending $p\to 0$ implies that
$m_s(x)$ is a martingale. \qd

\begin{lemma}\label{gamma-M} Let $-A$ be the generator of $T_t=e^{-tA}$ and
$\Gamma$ the associated gradient form. Assume that the filtration
$M_s$ is continuous. Let $x,y\in \A$. Then
\[ \langle  m(x),m(y)\rangle_t \lel 2\int_0^t \pi_s(\Gamma(x,x))ds
 \pl .\]
\end{lemma}

\begin{proof} Let us recall that for adapted process $(a_s)$ and $(b_s)$ we have
 \[ \langle a,b\rangle_s \lel \lim_{\pi,\U} \sum_j
 E_{s_j}((a_{s_{j+1}}-a_{s_j})^*(b_{s_{j+1}}-b_{s_j})) \pl .\]
The limit is taken in the weak sense. In particular, the bracket is
bilinear and vanishes on martingales, because then
$E_{s_j}(m_{s_{j+1}}-m_{s_j})=0$. It is best to start with
 \begin{align*}
 \pi_t(x^*x)-\pi_s(x^*x)&= m_t(x^*x)-m_s(x^*x)- \int_s^t
 \pi_r(A(x^*x)) dr \pl .
 \end{align*}
We use the notation $\pi_t(x)=m_t(x)+a_t(x)$ where $a_t(x)=\int_0^t
\pi_r(-Ax)dr$ and $m_t(x)$ is the martingale from the previous Lemma
\ref{Mey}. Then we find
 \begin{align*}
 &\pi_t(x^*x)-\pi_s(x^*x) \lel
 \pi_t(x)^*\pi_t(x)-\pi_s(x)^*\pi_s(x)\\
 &= (\pi_t(x)-\pi_s(x)+\pi_s(x))^*(\pi_t(x)-\pi_s(x)+\pi_s(x))\\
 &= (\pi_t(x)-\pi_s(x))^*\pi_s(x)+\pi_s(x)^*(\pi_t(x)-\pi_s(x))+
 (\pi_t(x)-\pi_s(x))^*(\pi_t(x)-\pi_s(x)) \\
 &= (m_t(x)-m_s(x))^*\pi_s(x)+(a_t(x)-a_s(x))\pi_s(x) \\
 &\pll +
 \pi_s(x)^*(m_t(x)-m_s(x))+\pi_s(x)^*(a_t(x)-a_s(x))\\
 &\pl +
 (m_t(x)-m_s(x))^*(m_t(x)-m_s(x))+(m_t(x)-m_s(x))^*(a_t(x)-a_s(x))\\
 &\pl +
(a_t(x)-a_s(x)) (m_t(x)-m_s(x))+(a_t(x)-a_s(x))^*(a_t(x)-a_s(x)) \pl
.
\end{align*}
After applying $E_{s}$ the first and the third term disappear. Then
we observe that
 \[ \|(a_t(x)-a_s(x))^*(a_t(x)-a_s(x))\|\kl (t-s)^2 \sup_r
 \|\pi_r(Ax)\| \]
and hence for bounded $Ax$ this terms disappears when the mesh size
goes to $0$. Since we assume that the filtration is continuous, we
know that $m_s(x)$ is norm continuous in $L_{2p}$, $p<\infty$. Thus
by uniform continuity we find
 \[ \lim_{\delta\to 0} \sup_{|s_{j+1}-s_j|<\delta} \sum_j
 \|(a_{s_{j+1}}(x^*)-a_{s_j}(x^*))(m_{s_{j+1}}(x)-m_{s_j}(x))\|_p \lel
 0 \pl .\]
Note that the $L_p$ continuity of $m_s$ implies the $L_p$ continuity
of $\pi_s$. Therefore we obtain in the limit (in the Riemann sense)
 \[ \langle m(x),m(x)\rangle_t \lel -\int_0^t
 \pi_r(A(x^*x)) dr+\int_0^t \pi_r(Ax^*)\pi_r(x)dr-
 \int_0^t \pi_r(x^*)\pi_r(Ax)dr \pl .\]
Here we use \eqref{l222} for $a=\Gamma(x,x)$. By assumption
$A\Gamma(x,x)\in \A$ for $x\in \A$ and hence we
$\pi_r(a)-E_{s_j}\pi_r(a)$ goes do $0$ uniformly in $|r-s_{j+1}|$.
By polarization the formula is true for all $x,y$. \qd

The main ingredient in the probabilistic approach towards Riesz
transforms  is to use L\'{e}vy's stopping time argument for the
Brownian motion (see however \cite{GuCR}, \cite{GuVa} for more
compact notation). Let $(B_t)$ be a classical brownian motion with
generator $ds$ (instead of the usual $\frac12 ds$) such that $B_0=a$
holds with probability $1$. Then we consider the stopping time ${\bf
t}_a=\inf\{t: B_t(\om)=0\}$. Instead of $\A$ we consider now the
tensor product $\A(B)\ten \A$ where $\A(B)$ is the algebra of
polynomials in the variables $B_t$. The new generator  is
 \[ \hat{A}\lel -\frac{d^2}{dt^2}+ A \pl .\]
This leads to
 \[ \hat{\Gamma}(x,y)\lel \Gamma(x,y)+\frac{d}{dt}x^*\frac{d}{dt}y
 \pl .\]
The Markov dilation is given by $\hat{\pi}_t(f\ten x)\lel
f(B_t)\ten\pi_t(x)$. Indeed, let $\hat{M}_s\lel M_s^B\ten M_s$ be
the von Neumann algebra given by the tensor product of the Brownian
motion observed until time $s$ and the von Neumann algebra $M_s$
given by the Markov dilation. Then
 \begin{align*}
  \hat{E}_s(f(B_t)\ten \pi_t(x)) &=
  E^B_s(f(B_t))E_s(\pi_t(x))
  \lel  T^B_{t-s}f(B_s) \pl     \pi_s(T_{t-s}(x)) \\
  &=  \hat{\pi}_s((T^B_{t-s}\ten
  T_{t-s})(f\ten x))
  \pl .
  \end{align*}
Here we use the fact that the brownian motion is the Markov process
for the generator $D^2=\frac{d^2}{dt^2}$ with corresponding
semigroup $T^B_t=e^{tD^2}$. For an element $x\in \A$ we use the
notation $Px\in L_{\infty}(\rz_+;N)$ given by the function
 \[ Px(t)\lel P_t(x) \pl .\]
We will also write $P'x$ for the function $\frac{d}{dt}P_tx$.
Harmonicity now leads to a well-known martingale property (again due
to Meyer).

\begin{prop}\label{ccbrack}  Let ${\bf t}_a$ be the stopping time as above. Then
 \[ n_t(x)\lel \hat{\pi}_{{\bf t}_a\wedge t}(Px) \]
is a martingale with bracket
 \[ \langle n(x) ,n(x)\rangle_t \lel 2\int_0^{{\bf t}_a\wedge t}
 \hat{\pi}_s(\hat{\Gamma}(Px,Px)) ds \pl .\]
\end{prop}

\begin{proof} We consider $y=\hat{\pi}_{{\bf t}_a}(Px)\lel \pi_{{\bf t}_a}(x)$. Let
us calculate the conditional expectation $\hat{E}_s$:
\begin{align*}
 \hat{E}_s(\pi_{\ttt_a}(x))
 &= 1_{\ttt_a\le s} \hat{E}_s(\pi_{\ttt_a}(x))+
 1_{\ttt_a>s}\hat{E}_s(\pi_{\ttt_a}(x)) \lel
  1_{\ttt_a\le s}\pi_{\ttt_a}(x) +
  1_{\ttt_a>s}E_s^B(\pi_s(T_{\ttt_a-s}(x)))
 \end{align*}
Now, we fix an $\om \in \Om$ such that $B_s(\om)=b$ and $s<\ttt_a$.
This means $b>0$. Then $\ttt_a-s$ is exactly stopping time until
$B_t-B_s$ hits $0$. Let us recall that (see \cite[page=25]{IMK})
 \[ \ez e^{-\la\ttt_b} \lel e^{-\sqrt{\la}b} \pl .\]
Using the spectral resolution $A=\int \la dE(\la)$ we get
 \[ (x,\ez T_{\ttt_b}(y))
 \lel \ez \int_0^{\infty} e^{-\la \ttt_b}d\nu_{x,y}(\la)
 \lel \int_0^{\infty} e^{-b\sqrt{\la}} d\nu_{x,y}(\la)
  \lel (x,P_b(x)) \pl .\]
By continuity,
 \begin{equation}\label{expev}
 \ez T_{\ttt_b}(y) \lel P_by
 \end{equation}
holds $L_p$'s. Hence we find
 \[ \ez(\pi_s(T_{\ttt_a-s}(x))|B_s=b) \lel \pi_s(P_{B_s}(x))
  \lel \hat{\pi}_s(Px)  \pl .\]
This proves the first assertion. For the second, we recall that
 \[ m_t(Px)\lel \hat{\pi}_t(Px)+\int_0^t \hat{\pi}_s(\hat{A}Px) ds \]
is a martingale and according to Lemma \ref{gamma-M} we have
 \[ \langle m(x),m(x)\rangle_t \lel 2\int_0^t
 \hat{\pi}_s\hat{\Gamma}(Px,Px) ds \pl .\]
Indeed, we may approximate $Px$ by function of the form $\sum_j
f_j\ten x_j$ in the graph norm of $\hat{A}$ such that $f_j(s)\lel 0$
for $s<0$ and the apply Lemma \ref{gamma-M}. So that we read
$P_sx=1_{s>0}P_sx$. However, we have
 \[ \frac{d^2}{dt^2}(P_t(x))\lel AP_tx \pl
 .\]
and hence $\hat{A}Px=0$. (This might no longer be true for the
approximations but it holds in the limit). Thus
$m=(\hat{\pi}_t(Px)_t)$ is martingale such that
 \[ \langle m,m\rangle_t \lel 2\int_0^t
 \hat{\pi}_s\hat{\Gamma}(Px,Px) ds \pl .\]
Thus by conditioning on the stopping time $\ttt_a$ we still have
  \begin{align*} \langle n,n\rangle_{t\wedge \ttt_a} &=  2\int_0^{t\wedge
  \ttt_a}
 \hat{\pi}_s\hat{\Gamma}(Px,Px) ds \pl . \qedhere
 \end{align*}
\qd

\begin{rem} {\rm Using the stopping $\ttt_a$ we can explicitly construct a
reversed Markov dilation for $P_s$ once we have one for $T_t$.
Indeed, let $B$ be a brownian motion such that ${\rm
Prob}(B_0=\infty)=1$ and $\ttt_b=\inf\{t :B_t(\om)=b\}$. Such a
random variable can be constructed using the limit for $a\to \infty$
of the finite brownian motions above. Then the random variable
 \[ \pi_b(x) \lel  \pi_{\ttt_b}(x) \]
satisfies
 \[ E_{c}(\pi_{b}(x))(\om)
 \lel \ez( \pi_c(T_{\ttt_{c-b}}(x))| \ttt_c(\om)=c)(\om)
 \lel \pi_{b}(P_{c-b}x)(\om)  \]
for all $c>b$. Here we use the random filtration $M_{\ttt_c]}$.
 }\end{rem}

In the following we fix the notation $\rho_ax=\pi_{\ttt_a}x$ for the
induced trace preserving $^*$-homo-morphism. For $\kappa>0$ we
follow a similar idea as in section 1 and construct martingales
$\rho_{a}^{\kappa}(x)$ such that
 \[ \langle \rho_{a}^{\kappa}(x),\rho_{a}^{\kappa}(x)\rangle_t
 \lel \int_0^{t\wedge \ttt_a}
 \hat{\pi}_s(\hat{\Gamma}(P^{\kappa}x,P^{\kappa}x)) ds  \pl .\]
Here $P^{\kappa}_s(x)=P_{\kappa s}x$. Indeed, we fix a partition
$\si=\{t_0,...,t_n\}$ and define
 \[ m_{\si} \lel \sum_{j=1}^n
 \hat{E}_{t_{j+1}}\pi_{\rho_a}(P^{\kappa-1}_{B_{t_j}}(x))-
 \hat{E}_{t_{j}}\pi_{\rho_a}(P^{\kappa-1}_{B_{t_j}}(x))
 \pl .\]
According to Lemma \ref{ccbrack} we obtain
 \[ \langle m_{\si},m_{\si}\rangle_t
 \lel 2 \sum_{t_{j+1}\le t} \int_{t_j}^{t_{j+1}}
 \pi_s\hat{\Gamma}(P^{\kappa}(B_{t_j}x),P^{\kappa}(B_{t_j}x))ds
 +\int_{t_{j_t}}^t
 \pi_s\hat{\Gamma}(P^{\kappa}(B_{t_{j_t}}x),P^{\kappa}(B_{t_{j_t}}x))ds
 \pl .\]
Passing to weak$^*$-limit we obtain $\rho_a^{\kappa}(x)$.

\begin{lemma}\label{Bkryest}
Let $\om$ and $t>0$  such that $t<\ttt_{a}(\om)$ and
$b=B_{t}(\om)>0$. Then
 \begin{align*}
 &\bigg(\hat{E}_t\langle \rho_a^{\kappa}x, \rho_a^{\kappa} x\rangle_{\infty}-
 \langle \rho_a^{\kappa}x, \rho_a^{\kappa}x\rangle_{t}\bigg)(\om)
  \kl c(\kappa) \hat{E}_{t}\bigg(\rho_a\big(\int_0^{\infty}
  P_s(\hat{\Gamma}(P_sx,P_sx)
  \min(s,b) ds\big)\bigg)(\om)
 \end{align*}
holds for $\kappa>1$. For $\kappa=1$ and $0<\beta<1$
 \begin{align*}
 &\bigg(\hat{E}_t\langle \rho_ax, \rho_ax\rangle_{\infty}-
 \langle \rho_ax, \rho_ax\rangle_{t}\bigg) (\om)
  \kl\frac {c}{\beta(1-\beta)^3} \pi_{t\wedge \ttt_a}\big(\int_0^{\infty}
  P_{\beta b+s }\hat{\Gamma}(P_sx,P_sx)
  \min(s,b) ds\big)(\om).
 \end{align*}
\end{lemma}
\begin{proof} Our starting point is
 \[ \langle \rho_{a}^{\kappa}(x) ,\rho_{a}^{\kappa}(x)\rangle_t \lel \int_0^{{\bf t}_a\wedge t}
  \hat{\pi}_s(\hat{\Gamma}(P^{\kappa}x,P^{\kappa}x)) ds  \pl .\]
For the little bmo norm this implies
 \[ \hat{E}_t\langle \rho_a^{\kappa}x, \rho_a^{\kappa} x\rangle_{\infty}-
 \langle \rho_a^{\kappa} x, \rho_a^{\kappa}x\rangle_{t} \lel \hat E_t\int_{t\wedge \ttt_a}^{\ttt_a} \hat{\pi}_s(\hat{\Gamma}(P^{\kappa}x,P^{\kappa}x))
 ds\pl .\]
Thus for $t>\ttt_a(\om)$ we have $0$.  Let us assume $t<\ttt_a(\om)$
and $b=B_t(\om)>0$. Then we observe that
 \begin{align*}
  &\ez(E_t\int_{t}^{\ttt_a} \pi_s(\hat{\Gamma}(P^{\kappa}_{B_s}x,P_{B_s}^{\kappa}x))
  ds|B_t=b)(\om)
   \lel \pi_t\bigg(
    \ez\int_0^{\ttt_b} T_s(\hat{\Gamma}(P^{\kappa}_{\tilde{B}_s}x,P^{\kappa}_{\tilde{B}_s}x)) ds \bigg)\pl
  .
  \end{align*}
On the right hand side we used the  notation $\tilde{B}_s$ for a
Brownian motion starting at $b$ and $\ttt_b$ is the stopping time at
$0$. Let us fix $y=P_bx$. We use a well-known formula for local
times (see \cite[Formula (11)]{Ba1})
  \begin{equation}\label{11} \ez \int_0^{\ttt_b} f(t,\tilde{B}_t) dt
 \lel \frac12 \int_0^{\infty} \int_{|b-s|}^{b+s} \int_0^{\infty}f(r,s)d\mu_t(r)
 dt ds \end{equation}
where $\int_0^{\infty} e^{-\la r}d\mu_t(r)=e^{-t\sqrt{\la}}$. This
implies
 \begin{align}
 &\ez\int_0^{\ttt_b} T_s(\hat\Gamma(P^{\kappa}_{\tilde{B}_s}x,P^{\kappa}_{\tilde{B}_s}x)) ds
  \lel \frac{1}{2} \int_0^{\infty} \int_{|b-s|}^{b+s}
 P_t\hat{\Gamma}(P_{\kappa s} x,P_{\kappa s}x) dt ds \label{here}
 \end{align}
For $\kappa>1$, let  $\al=\frac {\kappa+1}{2\kappa}$. Then we
observe with $\Gamma^2\gl 0$ and monotonicity from Proposition
\ref{monot} that
 \begin{align*}
 &\int_{|b-s|}^{b+s}
 P_t\hat{\Gamma}(P_{\kappa s} x,P_{\kappa s}x) dt
 \kl \int_{|b-s|}^{b+s}
 P_{t+\al\kappa s} \hat{\Gamma}(P_{(1-\al)\kappa s} x,P_{(1-\al)\kappa s}x)
 dt\\
 &\le \kla \int_{|b-s|}^{b+s} (t+\al\kappa s) dt \mer
 \frac{P_{|b-s|+\al\kappa s}(\hat{\Gamma}(P_{(1-\al)\kappa s} x,P_{(1-\al)\kappa
 s}x))}{|b-s|+\al\kappa s} \\
 &= \frac{2bs+\al \kappa s(b+s-|b-s|)}{|b-s|+\al\kappa s}
 P_{|b-s|+\al\kappa s}(\hat{\Gamma}(P_{(1-\al)\kappa s} x,P_{(1-\al)\kappa
 s}x)) \pl .
 \end{align*}
Note that $b+s-|b-s|=2\min(b,s)$.  For $s\gl b$ we use monotonicity
again and get
 \begin{align*}
 &\frac{2(1+\al \kappa) bs}{|b-s|+\al\kappa s}
 P_{s-b+\al\kappa s}\hat{\Gamma}(P_{(1-\al)\kappa s} x,P_{(1-\al)\kappa
 s}x)\\
 &\le \frac{2(1+\al \kappa) bs}{b+(\al\kappa-1)s}
 P_{b+(\al \kappa-1)s}\hat{\Gamma}(P_{(1-\al)\kappa s}
 x,P_{(1-\al)\kappa s}x) \\
 &\le
 b\frac{2(1+\al \kappa)}{(\al \kappa-1)}
 P_{b+(\al \kappa-1)s}\hat{\Gamma}(P_{(1-\al)\kappa s}
 x,P_{(1-\al)\kappa s}x)) \pl.
\end{align*}
Note that for $b>s$ we have $|b-s|+\al\kappa s=b+(\al\kappa-1)s$ and
hence
\begin{align*}
 & \frac{1}{2} \int_0^{\infty} \int_{|b-s|}^{b+s}
 P_t\hat{\Gamma}(P_{\kappa s} x,P_{\kappa s}x) dt ds\\
 &
 \kl \max((1+\al \kappa),\frac{1+\al \kappa}{\al \kappa -1})
 \int_0^{\infty} P_{b+(\al \kappa-1)s}\hat{\Gamma}(P_{(1-\al)\kappa
 s }x, P_{(1-\al)\kappa s}x) \min(b,s) ds \\
 &\lel  2\max(\frac{\kappa+3}2,\frac{3+ \kappa}{\kappa -1})
 \int_0^{\infty} P_{b+s}\hat{\Gamma}(P_{s}x, P_sx) \min(b,\frac{2s}{\kappa-1})
 \frac{ds}{\kappa-1} \pl.
 \end{align*}
We deduce the first assertion. For $\kappa=1,\beta<1$, let
$\al=\frac{\beta+1}{2}$ and $\gamma=\frac{1-\beta}{2}$. Then, for
$b>s$, $b-s+\al s\gl \beta b+\gamma s$ and hence
 \begin{align*}
 &  \frac{2bs+\al \kappa s(b+s-|b-s|)}{|b-s|+\al s}
  P_{|b-s|+\al s}(\hat{\Gamma}(P_{(1-\al)s} x,P_{(1-\al)
  s}x)) \\
  & \kl \frac{2s(b+\al s)}{\beta b+\gamma s}
 P_{\beta b+\gamma s}(\hat{\Gamma}(P_{(1-\al)s} x,P_{(1-\al)
 s}x)) \kl \frac {4s}\beta
 P_{\beta b+\gamma s}(\hat{\Gamma}(P_{\gamma s} x,P_{\gamma
 s}x))
 \pl .
 \end{align*}
For $s\gl b$ we have $s-b+\al s\gl \beta b+\gamma s$ and hence
 \begin{align*}
 &  \frac{2(1+\al) bs }{s-b + \al s}
  P_{s-b +\al s}(\hat{\Gamma}(P_{(1-\al)s} x,P_{(1-\al)
  s}x)) \\
 &\kl   \frac{2(1+\al) bs }{\beta b + \gamma s}
  P_{\beta b +\gamma s}(\hat{\Gamma}(P_{(1-\al)s} x,P_{(1-\al)
  s}x))
  \kl 2b\frac{1+\al}{\gamma}
  P_{\beta b +\gamma s}(\hat{\Gamma}(P_{(1-\al)s} x,P_{(1-\al)  s}x))
  \pl \\
  &\kl \frac{8b}{1-\beta}
  P_{\beta b +\gamma s}(\hat{\Gamma}(P_{\gamma s} x,P_{\gamma s}x))
  \pl.
\end{align*}
We deduce the assertion from a change of variables which leads to
$c(1-\beta)^{-3}$ or $c(1-\beta)^{-2}\beta^{-1}$. \qd

The lower estimate in the following result is well-known in the
commutative theory.

\begin{theorem}\label{hpc} Let $x\in \A$ and $2\le p<\infty$ and $\Gamma^2\gl 0$ and $\kappa \gl 1$. Then
  \begin{align*}
  &\|(\int_0^{\infty} P_s\hat{\Gamma}(P_{\kappa s}x,P_{\kappa s}x)\min(s,a) ds)^{1/2}
  \|_p
 \kl  \|\rho_a^{\kappa} x \|_{h_p^c} \pl.
  \end{align*}
For   $\kappa>1$.
 \begin{align*}  \|\rho_a^{\kappa} x \|_{h_p^c}
    \kl
   c_p(\kappa) \pl \|(\int_0^{\infty} P_s\hat{\Gamma}(P_sx,P_sx) s
  ds)^{1/2}\|_p \pl .
  \end{align*}
\end{theorem}

\begin{proof} For the lower
estimate, we calculate the conditional expectation of the square
function onto $\rho_a(N)$. Indeed, let $y\in N$ then
 \begin{align*}
  &\ez \tau(\rho_a(y)^*\int_0^{\ttt_a}\hat{\pi}_s(\hat{\Gamma}(P^{\kappa}x,P^{\kappa}x)) ds)
   \lel  \ez \int_0^{\ttt_a}
  \tau(\hat{E}_{s}(\rho_a(y^*))\hat{\pi}_s(\hat{\Gamma}(P^{\kappa}x,P^{\kappa}x)))ds \\
 &= \ez \int_0^{\ttt_a}
  \tau(\pi_s(P_{B_s}(y^*) \pi_s(\hat{\Gamma}(P_{\kappa B_s}x,P_{\kappa
  B_s}x))) ds
  \lel \ez \int_0^{\ttt_a}
   \tau(y^*P_{B_s}\hat{\Gamma}(P_{\kappa B_s}x,P_{\kappa B_s}x))ds \\
  &= \int_0^{\infty} \tau(y^*P_s\hat{\Gamma}(P_{\kappa s}x,P_{\kappa s}x)) \min(a,s) ds
  \pl .
  \end{align*}
For the upper estimate we note that the $L_p^cMO$ is given by
 \[ \|\rho_a^{\kappa} x\|_{L_p^cmo}
 \lel \|\sup_t \hat{E}_t\langle \rho_a^{\kappa}x, \rho_a^{\kappa}x\rangle_{\infty}-
 \langle \rho_a^{\kappa} x, \rho_a^{\kappa} x\rangle_{t} \|_{p/2}^{1/2} \pl .\]
We may assume $z=\int_0^{\infty} P_s\hat{\Gamma}(P_sx,P_sx) s ds\in
L_{p/2}(N)$. By Doob's inequality we find a $y\in L_{p/2}$ such that
\[ \pi_{\ttt_a\wedge t}(P_{B_t}z) \lel  \hat{E}_t(\rho_a(z)) \kl y
 \]
for all $t\gl 0$ and
\[ \|y\|_{p/2} \kl c_{p/2} \|\rho_a(z)\|_{p/2} \lel c_{p/2} \|z\|_{p/2} \pl. \]
With Lemma \ref{Bkryest} we deduce that
 \[ \hat{E}_t\langle \rho_a^{\kappa}x, \rho_a^{\kappa}x\rangle_{\infty}-
 \langle \rho_a^{\kappa} x, \rho_a^{\kappa} x\rangle_{t}
 \kl c(\kappa) \hat{\pi}_{t\wedge \ttt_a}(Pz)
 \kl c(\kappa)  y   \]
for all $t\gl 0$. This implies the upper estimate. \qd

\begin{rem}{\rm In the semi-commutative case where $P_t=P_t^{\rz^n}\ten id$ is the
Poisson semigroup on $\rz^n$ we have the estimate
 \[ \frac{P_{\beta t}}{\beta t} \kl \beta^{-n} \frac{P_{t}}{t} \pl
\]
which follows from the explicit representation as a convolution
kernel. Choosing $\beta=(1-\frac{1}{n})$ in Lemma \ref{Bkryest}, we
obtain a polynomial estimate
 \[ \|\rho_a(x)\|_{h_p^c} \kl c n^3 \|x\|_{H_p^c} \pl. \]
}\end{rem}

\begin{cor}\label{sqq1} Let $2\le p<\infty$ and $\Gamma^2\gl 0$. Then
 \[ \|(\int_0^{\infty} P_s\hat{\Gamma}_s(P_sx,P_s) sds)^{\frac
 12}\|_p
 \sim_{c(p,\kappa)} \lim_{a\to \infty} \|\rho_a^{\kappa} (x)\|_{h_p^c} \pl .\]
holds for all $\kappa>1$.
\end{cor}

\begin{lemma}\label{ctpp2} Let $x\in \A$ such that
 \[  \sup_s  \|A((P_sx)^2)\|_2  \pl <\pl \infty  \pl. \]
Then the martingale $x_t=\hat{E}_t(\rho_a(x))$ has continuous path.
Moreover, if $\Gamma^2\gl 0$, then every martingale $\rho_a(x)$ with
$x\in \A$ has continuous path.
\end{lemma}

\begin{proof} Let us assume $x$ selfadjoint (for convenience).
 We  follow Lemma \ref{ctpath} and observe that
the $4$-norm satisfies
 \begin{align*}
  \|x_t-x_s\|_4^4 &\lel
  \tau(x_t^4)+\tau(x_s^4)-4\tau(x_t^3x_s)-4\tau(x_s^3x_t)+4\tau(x_t^2x_s^2)+2\tau(x_tx_sx_tx_s)\\
  &\kl
  \tau(x_t^4)-3\tau(x_s^4)-4\tau(x_t^3x_s)+6\tau(x_t^2x_s^2)\\
  &=
   \tau(x_t^4)-\tau(x_s^4) -4(\tau(x_t^3x_s)-\tau(x_s^4))+
  6(\tau(x_t^2x_s^2)-\tau(x_s^4))
  \pl .
  \end{align*}
We note that
 \[ \tau(x_t^4) \lel \ez \tau(\pi_{\ttt_a\wedge t}(P_{B_{\ttt_a\wedge t}}(x))^4))
 \lel \ez \tau((P_{B_{\ttt_a\wedge t}}x)^4) \pl .\]
We use the It\^{o} formula for $f(s)=(P_{s}x)^4$ and obtain
 \begin{align*}
 &(P_{B_{\ttt_a\wedge t}}x)^4 \lel
  (P_{B_{\ttt_a\wedge s}}x)^4\\
 &+ \int_{\ttt_a\wedge s}^{\ttt_a\wedge
 t} (P'_{B_r}x(P_{B_r}x)^3+P_{B_r}xP'_{B_r}x(P_{B_r}x)^2
 +(P_{B_r}x)^2P'_{B_r}x P_{B_r}x+(P_{B_r}x)^3P'_{B_r}x) dB_r +\\
 & \int_{\ttt_a\wedge s}^{\ttt_a\wedge  t}\!\!\!
 \big(P''_{B_r}x(P_{B_r}x)^3+ P'_{B_r}xP'_{B_r}x(P_{B_r}x)^2
 +P'_{B_r}xP_{B_r}x P'_{B_r}xP_{B_r}x+P'_{B_r}x(P_{B_r}x)^2
 P'_{B_r}x\big)
 dr+..
 \end{align*}
Indeed, every term in the second line yields four terms for the
second  derivative in the next line. Taking the expectation it
suffices to estimate the terms with the second derivative (using the
unusual normalization $dr$ instead of $\frac12dr$). Thus a uniform
bound for $A^{1/2}x$ in the $4$ norm implies that
 \[ |\tau(x_t^4)-\tau(x_s^4)|\kl C \ez |\ttt_a\wedge t-\ttt_a\wedge
 s|\kl C |t-s|  \max\{\|Ax\|_4^2\|x\|_4^2+\|A^2x\|_4\|x\|_4^3\} \pl .\]
For the second term we observe that
 \[ \tau(x_t^3x_s) \lel \ez \tau\big(\pi_{\ttt_a\wedge t}((P_{B_{\ttt_a\wedge
  t}}x)^3)\pi_{\ttt_a\wedge s}(P_{B_{\ttt_a\wedge s}}x)\big)
 \lel \ez  \tau\big((P_{B_{\ttt_a\wedge t}}x)^3T_{\ttt_a\wedge t-\ttt_a\wedge
 s}P_{B_{\ttt_a \wedge s}}x\big)  \pl. \]
We have to invoke the It\^{o} formula for
 \begin{align*}
 &(P_{B_{\ttt_a\wedge t}}x)^3 \lel
  (P_{B_{\ttt_a\wedge s}}x)^3\\
 &+ \int_{\ttt_a\wedge s}^{\ttt_a\wedge
 t} (P'_{B_r}x(P_{B_r}x)^2+P_{B_r}xP'_{B_r}xP_{B_r}x
 +(P_{B_r}x)^2P'_{B_r}x) dB_r \\
 & +\int_{\ttt_a\wedge s}^{\ttt_a\wedge  t}\!\!\!
 \big(P''_{B_r}x(P_{B_r}x)^2+ P'_{B_r}xP'_{B_r}xP_{B_r}x
 +P'_{B_r}xP_{B_r}x P'_{B_r}x \big)
 dr+..
 \end{align*}
Here we use $P_s''$ for $\frac{d^2}{ds^2}P_s$. The first term
vanishes again due to the martingale property. Then we use
 \[ \|T_rP_{B_{\ttt_a\wedge s}}x-P_{B_{\ttt_a\wedge s}}x\|_4
 \kl \|T_rx-x\|_4\lel \|\int_0^r T_sAx\|_4\kl r\|Ax\|_4 \]
in
 \begin{align*}
 &\ez  \tau\big((P_{B_{\ttt_a\wedge t}}x)^3T_{\ttt_a\wedge t-\ttt_a\wedge
 s}P_{B_{\ttt_a \wedge s}}x\big) - \ez  \tau\big((P_{B_{\ttt_a\wedge s}}x)^4\big)
 \\
  &\lel \ez  \tau\big((P_{B_{\ttt_a\wedge t}}x)^3P_{B_{\ttt_a \wedge
  s}}x\big)- \ez  \tau\big((P_{B_{\ttt_a\wedge s}}x)^4\big)
 + \tau\big((P_{B_{\ttt_a\wedge t}}x)^3(T_{\ttt_a\wedge t-\ttt_a\wedge
 s}-I)P_{B_{\ttt_a \wedge s}}x\big) \pl .
 \end{align*}
Applying It\^{o}'s formula for $P_{B_{\ttt_a\wedge t}}x$ we find an
estimate of the form
 \[ |\tau(x_t^3x_s)-\tau(x_s^4)|\kl  |t-s| (\|x\|_4^3 \|A^2x\|_4+\|x\|_4^3\|Ax\|_4) \pl \]
for $A^2x,Ax\in L_4$. For the last term we have
 \[\tau(x_t^2x_s^2)\lel \ez \tau((P_{B_{\ttt_a\wedge t}}x)^2T_{\ttt_a\wedge
 t-\ttt_a\wedge s}(P_{B_{\ttt_a\wedge s}}x)^2) \pl.  \]
The It\^{o} formula for $(P_{B_{\ttt_a\wedge t}}x)^2$ is simpler
than above. According to our assumption
 \[ \|(T_{\ttt_a\wedge t-\ttt_a\wedge s}-I)((P_{B_{\ttt_a\wedge s}}x)^2)\|_2
 \kl C |\ttt_a\wedge t-\ttt_a\wedge s|\kl C |t-s| \pl .\]
Therefore the martingale satisfies the assumption of the
noncommutative Kolmogorov theorem due to \cite{Golu}. Finally, let
us assume that $\Gamma^2\gl 0$. Then we use $\hat{A}Px=0$ and
conlcude for selfadjoint $x$ that
 \begin{align*}
  A(P_sxP_sx)&= \hat{A}(P_sxP_sx) +\frac{d}{ds}(P_sx)^2
  \lel 2\hat{\Gamma}(P_sx,P_sx)+P_sxP_s'x+P_s'xP_sx \\
  &= 2\Gamma(P_sx,P_sx) + 2P_s'xP_s'x +P_sxP_s'x+P_s'xP_sx \pl .
  \end{align*}
This implies with $0\le \Gamma(P_sx,P_sx)\le P_s\Gamma(x,x)$ and
H\"older's inequality that
 \begin{align*}
  \|A(P_sxP_sx)\|_2
  &\le 2\|\Gamma(x,x)\|_2+2\|A^{1/2}x\|_4^2+2\|x\|_4\|A^{1/2}x\|_4
  \pl.
  \end{align*}
Thus for $x\in \A$, we have $\Gamma(x,x)\in \A$ and hence a uniform
estimate in $s$. \qd



The main advantage of the probabilistic model is that it allows to
consider time and space derivatives simultaneously. Let us recall
that in the space $h_p^c$, $1\le p<\infty$, we have an orthogonal
projection $P^{br}$ on the space of martingales
 \[ h_p^{br}\lel\{ \int x_sdB_s\pl:\pl (x_s) \mbox{ adapted} \} \]
Of course, we have to read $\int x_sdB_s$ as a stochastic integral
approximated by $\sum_{s_j} x_{s_j}(B_{s_{j+1}}-B_{s_j})$. We refer
to the classical literature for approximation of the stopped process
 \[ (\int x_sdB_s)_{\ttt_a} \lel \int_0^{\ttt_a} x_s dB_s \]
which remains  in $h_p^{br}$. Let us consider a martingale $z_t\in
L_{\infty}(\Om)\bar{\ten}N$. Then the brownian projection $z_t$ is
the unique martingale $b_t\in h_p^{br}$ such that
 \[ \langle b,\int x_sdB_s\rangle_t \lel \langle z, \int
 x_sdB_s\rangle_t \]
holds for every adapted process $x$. Let us consider for example the
simple tensor $z=f\ten y$. We may assume
 \[ E_s^{B}(f) \lel \int_0^s g(r)dB_r \]
Let $z_s=\hat{E}_s(z)$. Note that
 \begin{align}
 z_{s+h}-z_s &=    E_{s+h}^B(f)y_{s+h}- E_s^B(f)\ten y_s \nonumber \\
 &= (E_{s+h}^B(f)-E_s^B(f))y_s +
 E_s^B(f)(y_{s+h}-y_s)+
 (E_{s+h}^B(f)-E_s^B(f))(y_{s+h}-y_s) \pl . \label{dddpo}
 \end{align}
Thus for $m=x_s(B_{s+h}-B_s)$ we find
 \begin{align*}
 &\hat{E}_s((m_{s+h}^*-m_s^*)(z_{s+h}-z_s))\lel
    x_s^*\hat{E}_s((B_{s+h}-B_s)(z_{s+h}-z_s))\\
 &= x_s^*\hat{E}_s((B_{s+h}-B_s)(E_{s+h}^B(f)-E_s^B(f))y_s + x_s^*\hat{E}_s((B_{s+h}-B_s)E_s^B(f)(y_{s+h}-y_s))\\
 &\quad +
 x_s^*\hat{E}_s((B_{s+h}-B_s)(E_{s+h}^B(f)-E_s^B(f))(y_{s+h}-y_s))\\
 &= x_s^*\int_{s}^{s+h} g(r) dr \pl y_s
 \pl .
 \end{align*}
Indeed, for the  two additional terms we use commutativity and
$\hat{E}_s\lel E_s^B\ten E_s$. This yields $0$ in both cases. Thus
in general we find
 \[ \langle \int  a_sdB_s,z\rangle_t \lel \int_0^t g(r) E_r(y) dr
 \quad \mbox{and} \quad  b_t\lel \int_0^t g(r)E_r(y)dB_r \pl .\]
This shows us how to extend the projection $P^{br}$ by linearity.
Since this procedure is less known in the non-commutative context we
shall also show continuity with respect to the $h_p$ norm. We come
back to \eqref{dddpo} and observe as above with the help of
orthogonality that
 \begin{align*}
  &\hat{E}_s((z_{s+h}-z_s)^*(z_{s+h}-z_s))\lel
  \hat{E}_s(((E_{s+h}^B(f)-E_s^B(f))y_s)^*(E_{s+h}^B(f)-E_s^B(f))y_s)\\
  &\pl + \hat{E}_s((E_s^B(f)(y_{s+h}-y_s))^*E_s^B(f)(y_{s+h}-y_s))\\
  &\pl + \hat{E}_s(((E_{s+h}^B(f)-E_s^B(f))(E_{s+h}^B(f)-E_s^B(f))
  (y_{s+h}-y_s))^*(y_{s+h}-y_s)) \pl .
   \end{align*}
For the last term we get
 \begin{align*}
   &\hat{E}_s(((E_{s+h}^B(f)-E_s^B(f))(E_{s+h}^B(f)-E_s^B(f))
  (y_{s+h}-y_s))^*(y_{s+h}-y_s)) \\
 &=  \int_s^{s+h} g(r)dr \pl E_s((y_{s+h}-y_s)^*(y_{s+h}-y_s)) \pl
  .
  \end{align*}
However, the Burkolder inequality  implies that
 \[ \|\sum_j E_{s_j}(|y_{s_{j+1}}-y_{s_j}|^2) \|_{p/2} \kl c(p)
 \|y\|_p^2 \pl .\]
Thus for bounded $g$, the last term vanishes as long as the mesh
size of the partition goes to $0$. This yields
 \[ \langle z,z\rangle_t \lel \langle b,b\rangle_t+ \langle \int_0^{\infty} E_s(f)dy_s,
 \int_0^{\infty} E_s(f)dy_s\rangle_t  \pl .\]
By approximation and linearity we deduce that
 \begin{align}\label{orth}
   \langle P^{br}(z),P^{br}(z)\rangle_t+\langle (I-P^{br})(z),(I-P^{br})(z)\rangle_t
   \lel \langle z,z\rangle_t \pl .
   \end{align}
\begin{lemma}\label{Pcont}  Let $1<p<\infty$. Then $P^{br}$ and $(I-P^{br})$ are bounded, selfadjoint preserving maps on
$L_p(\hat{M})$.
\end{lemma}

\begin{proof} By duality it suffices to consider $2\le p<\infty$.  We see
that on a dense set of martingales of the form $z=\sum_j f_j\ten
y_j$, the images $P^{br}(z)$ have continuous path and satisfy
 \[ \langle P^{br}(z),P^{br}(z)\rangle_t \kl \langle z,z\rangle_t \pl .\]
Thus the Burkholder-inequalities imply that
 \[ \|P^{br}(z)\|_{h_p^c} \kl \|z\|_{h_p^c} \kl c(p) \|z\|_p \pl .\]
Note that that $P^{br}(z^*)=P^{br}(z)^*$. Since $P^{br}(z)$ has
continuous path we deduce from \cite{JK} that
 \[ \|P^{br}(z)\|_{L_p}\kl c_1(p) \pl \|z\|_{h_p} \kl c_1(p)c(p) \|z\|_p
 \pl .\]
The assertion follows by density. Moreover, the jump parts of $z$
are mapped to $(I-P^{br})(z)$.\qd

\begin{lemma}\label{horcomp} Let $x\in N$, then
 \begin{enumerate}
 \item[i)] $\langle P^{br}\rho_a(x),P^{br}\rho_a(x)\rangle_t
 \lel 2
 \int_0^{\ttt_a\wedge t} \hat{\pi}_r(|P'x|^2) dr$;
 \item[ii)] $\langle (Id-P^{br})\rho_a(x),(Id-P^{br})\rho_a(x)\rangle_t \lel
 2\int_0^{\ttt_a\wedge t} \hat{\pi}_r(\Gamma(Px,Px))dr$.
 \end{enumerate}
\end{lemma}

\begin{proof} According to Proposition \ref{ccbrack} and
\eqref{orth}it suffices to show that
 \begin{equation}
 \label{fr}   \langle P^{br}(\pi_{\ttt_a}(x)), P^{br}(\pi_{\ttt_a}(x))\rangle
 \lel 2\int_0^{t\wedge \ttt_a} \hat{\pi}_s(|P'x|^2) ds \pl .
 \end{equation}
We deduce from
 \[ \langle n_{\ttt},m_{\ttt}\rangle_s \lel
 \langle n,m\rangle_{s\wedge \ttt} \]
that $P^{br}$ commutes with stopping times.  Therefore it suffices
to consider the martingale $m_t(x)\lel \hat{\pi}_t(Px)+\int_0^t
\hat{\pi}_s(\hat{A}Px)ds=\hat{\pi}_t(Px)$ and calculate the
component corresponding to the brownian motion.  By approximation it
suffices to consider $n\lel (y\ten f) (B_{s+h}-B_s)$ such that $f$
is a $\Si_s$-measurable bounded function.  Let $F_t:M\to N$ be the
conditional expectation corresponding to the trace preserving map
$\pi_t$ and $\tilde{y}\lel F_{s+h}(y)\in N$. Let $\mu$ be the
spectral measure such that
 \[ \tau(\tilde{y}f(A)x)\lel \int_0^{\infty} f(\la) d\mu(\la) \pl .\]
Then we have
 \begin{align*}
 &\ez \tau((y\ten f)(B_{s+h}-B_s)\pi_{s+h}(P_{B_{s+h}}x))\\
 &= \ez (B_{s+h}-B_s)f \tau_N(F_{s+h}(y)P_{B_{s+h}}x)
 \lel  \int_0^{\infty} \ez (B_{s+h}-B_s)f e^{-\sqrt{\la}B_{s+h}}
 d\mu(\la)\pl .
 \end{align*}
To be more precise, we replace $g(t)=P_tx$ by a function
$h(-\sqrt{A}t)x$ such that $h(z)=e^{z}v_{\eps}(z)$ such that
$v_{\eps}(z)$ vanishes at $0$ and converges to $1$ as $\eps$ goes to
$0$. Using a stopping time $\ttt_a^{\delta}$ which stops the
brownian motion at  $\delta>0$ this calculation can be justified. By
It\^{o}'s formula we have
 \[ e^{-\sqrt{\la}B_{s+h}}\lel
 e^{-\sqrt{\la}B_{s}}-\sqrt{\la} \int_{s}^{s+h}
 e^{-\sqrt{\la}B_r}dB_r+ \la \int_s^{s+h} e^{-\sqrt{\la}B_r}dr \pl
 .\]
This yields
  \begin{align*}
 & \ez \tau((y\ten f)(B_{s+h}-B_s)\pi_{s+h}(P_{B_{s+h}}x))
 \lel   \ez \tau((y\ten f)(B_{s+h}-B_s)\pi_{s+h}(P_{B_s}(x))\\
 &+  \int_{s}^{s+h} \ez \tau((y\ten f)\pi_{s+h}(P'_{B_r}x))\pl 2 dr
  +  \int_s^{s+h} \ez \tau((y\ten f)(B_{s+h}-B_s) \pi_{s+h}(P''_{B_r}x))
  dr\\
 &=\int_{s}^{s+h} \ez \tau((y\ten f)\pi_{s+h}(P'_{B_r}x)) dr\\
 & + \int_s^{s+h} \ez \tau((y\ten f)(B_{s+h}-B_r) \pi_{s+h}(P''_{B_r}x))
  + \int_s^{s+h} \ez \tau((y\ten f)(B_{r}-B_s) \pi_{s+h}(P''_{B_r}x))
  dr \pl .
  \end{align*}
By independence the last two terms are $0$. Note that
 \begin{align*}
 \|\pi_{s+h}(y)-\pi_r(y)\|_2^2 &=
 2(\tau(y^*y)-\tau(T_{s+h-r}(y)^*y)) \kl (s+h-r) \|Ay\|_2\|y\|_2
 \end{align*}
implies
 \[ 1_{B_r>0} \|\pi_{s+h}(P'_{B_r}x)-\pi_r(P'_{B_r}y)\|_2^2
 \kl (s+h-r) \|A^2y\|_2\|Ay\|_2 \pl .\]
Thus by continuity, we obtain
 \[ b_t \lel \int_0^t \pi_r(P'_{B_r}x) dB_r \pl .\]
Hence for the bracket, we deduce \eqref{fr} (with $d\langle
B_r,B_r\rangle =2dr$).\qd

\begin{lemma}\label{erhalt}  Let $x,y\in L_2$. Let $\pr$ be the projection onto $(\ker(A))^{\perp}$. Then
 \[ \lim_{a\to \infty} \tau(\langle (I-P^{br})\rho_a(x),
 (I-P^{br})\rho_a(y)\rangle_{\ttt_a})
 \lel \frac{1}{4} \tau((I-\pr)(x)^*(I-\pr)(y)) \pl .\]
\end{lemma}

\begin{proof} By polarization we have
 \[ \langle
 (Id-P^{br})\pi_{\ttt_a}(x),(Id-P^{br})\pi_{\ttt_a}(y)\rangle_{\infty}
 \lel 2\int_0^{t\wedge \ttt_a} \hat{\pi}_s\Gamma(Px,Py) ds \pl .\]
Thus taking the trace yields
 \begin{align*}
 \tau(\langle
 &(Id-P^{br})\pi_{\ttt_a}(x),(Id-P^{br})\pi_{\ttt_a}(y)\rangle_{\infty})\\
 &= 2\ez \int_0^{\ttt_a} \tau(\hat{\pi}_s\Gamma(Px,Py)) ds
 \lel  2\ez \int_0^{\ttt_a} \tau(P_{B_s}x^*AP_{B_s}y) ds  \pl .
 \end{align*}
Again, we may use polarization and hence it suffices to establish
the result for $x=y$. Let $dE_{\la}$ be the spectral measure of $A$
and $\om_x(T)=(x,Tx)$. We use the well-known formula (see
\cite{Ba1})
  \begin{equation}\label{bkrr}\ez \int_0^{\ttt_a}f(B_s)ds \lel \int_0^{\infty} \min(a,s)
 f(s) ds \pl .\end{equation}
This implies
 \begin{align*}
 & \ez \int_0^{\ttt_a} \tau(P_{B_s}x^*AP_{B_s}x) ds
 \lel  \ez \int_0^{\ttt_a} (x,P_{2B_s}A(x))ds \\
 &= \int_0^{\infty} \ez \int_0^{\ttt_a} e^{-2\sqrt{\la}B_s}\la  ds
 \om_x(dE_{\la}) \lel
   \int_0^{\infty} \int_0^{\infty} \min(s,a)e^{-2\sqrt{\la}s} ds
 \om_x(dE_{\la}) \pl . \end{align*}
For $\lim_{a\to \infty}$ we find
 \[ \la \int_0^{\infty} e^{-2\sqrt{\la}s} sds \lel
  \frac{1}{4}  \int_0^{\infty} e^{-2\sqrt{\la}s} (2\sqrt{\la}s)^2
  \frac{ds}{s} \lel \frac{1}{4}
  \]
for $\la>0$ and equals $0$ else. Thus we get
 \[ \lim_{a\to \infty}
 \tau(\langle
 (Id-P^{br})\pi_{\ttt_a}(x),(Id-P^{br})\pi_{\ttt_a}(x)\rangle_{\infty})
 \lel \frac{1}{2} \langle (I-\pr)(x),(I-\pr)(x)\rangle  \pl .\]
Thus the general formula is established by polarization.\qd

The next Lemma deals with Hardy spaces and follows closely Bakry's
proof.

\begin{lemma}\label{nonsa} Assume that $\Gamma^2\gl 0$ and $(T_t)$ admits
Markov dilation path. Then
 \[ \sup_a \|(I-P^{br})\rho_a(A^{1/2}x)\|_{H_p^{c}}
 \kl c(p) \pl \|\Gamma(x,x)^{\frac12}\|_p \pl .\]
holds for $2<p<\infty$.
\end{lemma}

\begin{proof} Let $x\in \A$. We consider the function
 \[  f(s)\lel \Gamma(P_sx,P_sx) \]
and claim that $y_t=\hat{\pi}_t(f)$ is a submartingale. Indeed, we
know that
 \[ m_t(f)\lel \pi_t(f)+\int_0^t \hat{\pi}_r(\hat{A}f) dr \]
is a martingale. This implies that
 \begin{align*}
  E_s(\hat{\pi}_t(f))\lel E_s(m_t(f))-E_s(\int_0^t \hat{\pi}_r(\hat{A}f) dr)
  \lel m_s(f) - E_s(\int_0^t \hat{\pi}_r(\hat{A}f) dr) \pl .
  \end{align*}
Let us calculate the right hand side:
 \[ \frac{\partial f}{\partial s}(s)
 \lel \Gamma(P_s'x,P_sx)+\Gamma(P_sx,P_s'x)\]
and hence
 \begin{align*} \frac{\partial^2 f}{\partial s^2}(s)
 &=  \Gamma(P_s''x,P_sx)+\Gamma(P_sx,P_s''x)+
 2\Gamma(P_s'x,P_s'x)\\
 &=
 \Gamma(AP_sx,P_sx)+\Gamma(P_sx,AP_sx)+ 2\Gamma(P_s'x,P_s'x)
 \pl .
 \end{align*}
Therefore we obtain
 \begin{equation}\label{subh10}
 -\hat{A}(f) \lel (\frac{\partial^2}{\partial s^2}-A)(f)
 \lel 2\Gamma^2(P_sx,P_sx)+ 2\Gamma(P_s'x,P_s'x) \pl .
 \end{equation}
The same equation will allow us to estimate the increasing part of
the bracket $y_t$ defined as the limit of
 \begin{equation}\label{opo}
  \langle y\rangle_t \lel \lim_{\si} \sum_j E_{t_{j}}(y_{t_{j+1}}-y_{t_j})
  \end{equation}
where the limit is taken along some ultrafilter on partitions of the
interval $[0,t]$. Let $r_t\lel \int_0^t \hat{\pi}_s(\hat{A}f) ds$.
Clearly, the bracket operation vanishes on the martingale part. We
obtain
 \begin{align*}
 E_{t_j}(\int_{t_j}^{t_{j+1}} \hat{\pi}_s(\hat{A}f)ds) &=
 \int_{t_j}^{t_{j+1}}\hat{\pi}_{t_j}(\hat{T}_{s-t_j}\hat{A}f) ds
 \lel \pi_{t_j}(f-\hat{T}_{t_{j+1}}f) \approx -(t_{j+1}-t_j)\pi_{t_j}(\hat{A}f)   \pl .
 \end{align*}
Thus by $L_p$  continuity of $\hat{\pi}_s(\hat{A}f)$ we find
 \[ \langle y\rangle_t \lel \int_0^t \hat{\pi}_s(2\Gamma^2(P_sx,P_sx)+ 2\Gamma(P_s'x,P_s'x))ds \pl .\]
This does not change if we add stopping times, i.e. we have
 \begin{align*}
 \langle y_{\ttt_a}\rangle_t &=  2 \int_0^{t\wedge \ttt_a}
 \hat{\pi}_s(\Gamma^2(P_sx,P_sx)) ds + 2
 \int_0^{t\wedge \ttt_a} \hat{\pi}_s(\Gamma(P'_sx,P_s'x)) ds\\
  &\gl \langle (I-P^{br})(\pi_{\ttt_a}(A^{1/2}x)),
 (I-P^{br})(\pi_{\ttt_a}(A^{1/2}x))\rangle_t \pl .
 \end{align*}
According to Lemma \ref{potential} we find with $p/2>1$ that
 \begin{equation}
 \label{majoriz}
  \|\langle y \rangle_{\ttt_a} \|_{\frac p2}
  \kl c_p \| y_{\ttt_a}\|_{\frac p2} \pl
  .
  \end{equation}
It is time to apply  subharmonicity again. Now in the form
 \[ \Gamma(P_sx,P_sx)\kl P_s(\Gamma(x,x)) \pl .\]
Therefore
 \[ \hat{\pi}_{\ttt_a}(\Gamma(Px,Px))
 \kl \hat{\pi}_{\ttt_a}(P\Gamma(x,x)) \lel \rho_a(\Gamma(x,x))  \pl.\]
Now, we note that $\pi_{\ttt_a}:N\to M$ is a trace preserving
$^*$-homomorphism. This implies
 \[ \|\rho_a(\Gamma(x,x))\|_{\frac p2}\lel
 \|\Gamma(x,x)\|_{\frac p2} \pl .\]
We deduce that
 \begin{align*}
  \|(I-P^{br})\rho_a(A^{1/2} x)\|_{h_p^c} &\kl c(p) \| \langle
 (I-P^{br})\rho_a x, (I-P^{br})\rho_a x\rangle \|_{p/2}^{1/2}
 \kl c'(p) \|\Gamma(x,x)^{\frac12}\|_p \pl .
 \end{align*}
We may replace the $h_p^c$ norm by the $H_p^c$ norm, because
$\rho_a(x)$ and $P^{br}(\rho_a(x))$ have continuous path. Thus
$I-P^{br}(\rho_a(x))$ also has continuous path. \qd

We are now well-prepared for our main result on Riesz transforms.

\begin{theorem}\label{lowd} Let $T_t$ be a completely positive selfadjoint
semigroup with  negative generator $A$, which admits  a Markov
dilation and satisfies $\Gamma^2\gl 0$. Let $2<p<\infty$. Then
 \[ \|A^{\frac12}x\|_p \kl c_p \pl  \max\{\|\Gamma(x,x)^{\frac12}\|_p,\|\Gamma(x^*,x^*)^{\frac12}\|_p\}  \pl
 \]
holds for all  $x$.
\end{theorem}

\begin{proof} Let $\delta>0$. Using the fact that $\pr$ is a contraction,
we may find  $y_0\in L_{p'}(N)$ such that $\pr(y_0)=0$,
$\|y_0\|_{p'}\le 1$  and
 \begin{equation}\label{naa}
  \|A^{1/2}x\|_p\kl 2(1+\delta)  |\tau(y_0^*A^{1/2}x)| \pl .
 \end{equation}
By approximation we may assume that $y_0\in L_2(N)$ and still
satisfies \eqref{naa}. Since
$\rg(A^{1/2})=\ker(A^{1/2})^{\perp}=(\ker(A))^{\perp}$ we may
approximate $y_0$ by $A^{1/2}y\in L_{p'}(N)$ such that
$\|A^{1/2}y\|_{p'}\le 1$ and
 \[ \|A^{1/2}x\|_p \kl 3 |\tau(A^{1/2}y^*A^{1/2}x)| \pl .\]
We fix $a>0$ and according to \cite{JK} we decompose
$\rho_a(A^{1/2}y)\lel m_c+m_r+m_d$ such that
 \[ \|m_c\|_{h_{p'}^c} + \|m_r\|_{h_{p'}^r} + \|m_d\|_{h_{p'}^d}
 \kl c(p') \pl \|\rho_a(A^{1/2}y)\|_{p'} \kl 2 c(p')\pl .\]
Since $\rho_a(A^{1/2}x)$ has continuous path (see Lemma
\ref{ctpp2}), we know that
 \[ \langle m_d^*,(I-P^{br})\rho_a(A^{1/2}x)\rangle \lel 0 \pl .\]
Therefore we obtain from Lemma \ref{nonsa} that
 \begin{align*}
 &|\ez \tau
 ((I-P^{br})\rho_a(A^{1/2}y^*)(I-P^{br})(\rho_a(A^{1/2}x))|\\
 &=  |\ez \tau ( (m_c^*+m_r^*+m_d^*)(I-P^{br})(\rho_a(A^{1/2}x))|\\
 &= |\ez \tau ( ((I-P^{br}(m_c)^*+ (I-P^{br})(m_r^*))(I-P^{br})(\rho_a(A^{1/2}x))|\\
 &= |\ez \tau( \langle m_c,(I-P^{br})(\rho_a(A^{1/2}x)\rangle)|
  + |\overline{ \ez \tau( \langle
 m_r^*,(I-P^{br})(\rho_a(A^{1/2}x^*)\rangle)}|\\
 &\le \|m_c\|_{h_{p'}^c} \|(I-P^{br})(\rho_a(A^{1/2}x))\|_{h_p^c}
 + \|m_r^*\|_{h_{p'}^c}
  \|(I-P^{br})(\rho_a(A^{1/2}x^*))\|_{h_p^c}\\
 &\le c(p')c(p)  \|\rho_a(A^{1/2}y)\|_{p'} (
 \|\Gamma(x,x)^{\frac12}\|_p + \|\Gamma(x^*,x^*)^{\frac12}\|_p)
 \pl .
 \end{align*}
Note that $(I-\pr)(A^{1/2}x)=A^{1/2}x$. Therefore Lemma \ref{erhalt}
shows that
 \begin{align*}
 |\tau(A^{1/2}y^*A^{1/2}x)| &\le  4 \lim_{a\to \infty } |\ez \tau
 ((I-P^{br})\rho_a(A^{1/2}y^*)(I-P^{br})(\rho_a(A^{1/2}x))| \\
 &\le 8 c(p')c(p) (\|\Gamma(x,x)^{\frac12}\|_p + \|\Gamma(x^*,x^*)^{\frac12}\|_p)
 \pl .
 \end{align*}
By our choice of $y$ we deduce the assertion. \qd

As a further application we compare the martingale $H_p$-norms and
the semigroup $H_p$-norms from \cite{JLX}.

\begin{theorem}\label{hpp} Let $2\le p<\infty$ and $\Gamma^2\gl 0$ and $\kappa>1$.
\begin{enumerate}
\item[i)] Then
 \[ \|x\|_{H_p^c} \sim_{c(p)} \lim_a \|P^{br}\rho_a^{\kappa}(x)\|_{h_p^c}
 \sim_{\tilde{c}(p)} \lim_a \|\rho_a^{\kappa}(x)\|_{H_p^c} \pl .\]
\item[ii)]  If moreover, the assumption of Lemma \ref{cond} or
Lemma \ref{cond3} is satisfied, then
 \[ \|x\|_{H_p^c} \sim_{c(p)} \lim_a
 \|(I-P^{br})\rho_a^{\kappa}(x)\|_{h_p^c}\pl .\]
and
 \[ \|A^{1/2}x\|_{H_p^c(P)} \kl c(p) \|\Gamma(x,x)^{\frac12}\|_p \pl
 .\]
\end{enumerate}
\end{theorem}

\begin{proof} Let us start with an easy application of Theorem \ref{hpc}, namely the condition $\Gamma^2\gl 0$ implies
 \begin{align*}
 &\|\int_0^{\infty} \hat{\Gamma}(P_sx,P_sx) s^2\frac{ds}{s} \|_{p/2}^{1/2}
 \lel   (\kappa+1)  \|\int_0^{\infty} \hat{\Gamma}(P_{(\kappa+1)s}x,P_{(\kappa+1)s}x) sds
 \|_{p/2}^{1/2}\\
 &\le (\kappa+1)\pl
  \|\int_0^{\infty} P_s\hat{\Gamma}(P_{\kappa s}x,P_{\kappa  s}x) sds
 \|_{p/2}^{1/2}\\
 &\le   (\kappa+1) \pl \lim_{a} \|\int_0^{\infty} P_s\hat{\Gamma}(P_{\kappa s}x,P_{\kappa s}x) \min(s,a) sds
 \|_{p/2}^{1/2}
 \kl  (\kappa+1) \lim_{a} \|\rho_a^{\kappa} (x)\|_{h_p^c} \pl .
 \end{align*}
The converse is given by Theorem \ref{hpc}. The same argument in
combination with Lemma \ref{horcomp} also shows that
 \[ \|\int_0^{\infty} |P_s'x|^2 sds\|_{p/2}^{1/2}
 \sim_{c(p,\kappa)}\lim_{a} \|P^{br}\rho^{\kappa}_a(x)\|_{h_p^c}  \]
and
 \[ \|\int_0^{\infty} \Gamma(P_sx,P_sx) sds\|_{p/2}^{1/2}
 \sim_{c(p,\kappa)}\lim_{a} \|(I-P^{br})\rho^{\kappa}_a(x)\|_{h_p^c}  \pl. \]
We refer to \cite{JR1} for
 \[ \|\int_0^{\infty} \Gamma(P_sx,P_sx) sds\|_{p/2}^{1/2}
 \kl c(p) \|x\|_{H_p^c(P)} \pl .\]
Assuming the condition of Lemma \ref{cond3} for $\Gamma$ or under
the assumption of Lemma \ref{cond} we have
 \[ \|\int_0^{\infty} \Gamma(T_sx,T_sx) sds\|_{p/2}^{1/2}  \kl c(p) \|\int_0^{\infty} \Gamma(P_sx,P_sx)
 sds\|_{p/2}^{1/2} \pl .\]
Thus Theorem \ref{stein3} iii) yields the missing estimate in ii),
because the $H_p^c(P)$ and $H_p^c(T)$ are comparable, see again
\cite{JR1}. In that situation the last assertion follows from Lemma
\ref{nonsa}. \qd

\section{BMO spaces}

In the recent years the theory of BMO spaces has been extended to
semigroups of positive operators on $\rz^n$ assuming that their
kernels satisfy certain regularity conditions (see \cite{DuLi}). In
this part we compare different candidates for the  BMO-norm. Our
main motivation is the Garsia norm for the Poisson semigroup on the
circle. In full generality we define for a semigroup $T_t$ of
completely positive maps the norm
 \[ \|x\|_{BMO_c(T)} \lel \sup_t \|
 T_t|x|^2-|T_tx|^2\|_{\infty}^{1/2}\pl .\]

\begin{lemma}\label{bmo1}  Let  $(T_t)$ be a semigroup
of completely positive maps on a von Neumann algebra. Then
$BMO_c(T)$ defines a normed space. Moreover, if $(T_t)$ has a
reversed martingale dilation and $\Gamma^2\gl 0$, then the reversed
martingale $(\pi_s(T_s))$ satisfies
 \[ \|\pi_0(x)\|_{bmo_c} \lel \|x\|_{BMO_c(T)} \pl .\]
\end{lemma}

\begin{proof} Let us fix $t>0$ and define the homogenous expression $\|x\|_t=
\|T_t|x|^2-|T_tx|^2\|^{1/2}$. Let $d=T_t(1)$ and $e$ the support
projection of $d$. It is easy to see that $T_t(N)\subset eNe=M$.
Then $\tilde{T}_t(x)=d^{-1/2}T_t(x)d^{-1/2}$ is a well-defined
unital completely positive map $\tilde{T}_t:N\to M$. Let
$N\ten_{T_t}M$ be the Hilbert $C^*$-module over $M$ with $M$ inner
product
 \[ \langle a\ten b,c\ten d\rangle \lel b^*\tilde{T}_t(a^*c)d \pl. \]
Since $\tilde{T}_t$ is unital we obtain $^*$-homomorphism $\pi:N\to
{\mathcal L}(N\ten_{\tilde{T}_t}M)\cong {\rm Mult}(K(\ell_2)\ten M)$
such that
 \[ \tilde{T}_t(x) \lel e_{11}\pi(x)e_{11} \pl .\]
This implies
 \[ T_t(x)\lel d^{1/2}e_{11}\pi(x)e_{11}d^{1/2} \pl .\]
Therefore we get
  \begin{align*}  & T_t(x^*x)-T_t(x^*)T_t(x)
  \lel d^{1/2}e_{11}\pi(x)^*\pi(x) e_{11}d^{1/2}-
 d^{1/2}e_{11}\pi(x)^*e_{11}d^{1/2}d^{1/2}e_{11}\pi(x)
 e_{11}d^{1/2}\\
 &= d^{1/2}e_{11}\pi(x)^*(1-e_{11}de_{11})\pi(x)e_{11}d^{1/2}
 \pl .\end{align*}
This implies that the linear map $u:N\to {\mathcal
L}(N\ten_{\tilde{T}_t}M)$ defined by
 \[ u(x) \lel (1-e_{11}de_{11})^{\frac12} \pi(x)e_{11}d^{1/2} \]
is an isometric embedding of $N$ equipped with the norm $\|\pl
\|_t$. An alternative proof can be derived from \eqref{wies}
 \[ T_t|x|^2-|T_tx|^2 \lel 2 \int_0^t T_{t-s}\Gamma(T_sx,T_sx) ds
 \pl .\]
Thus the GNS construction for the positive form $T_{t-s}\Gamma$
allows us to find linear maps $u_{ts}:N\to C(N)$ such that
 \[ T_t|x|^2-|T_tx|^2 \lel \int_0^t |u_{ts}(x)|^2 ds \pl.\]
This provides an embedding in $L_2^c([0,t])\ten_{\min}C(N)$. Now, we
assume that $T_t$ admits a martingale dilation as in section 1. We
consider the martingale $m_s\lel \pi_s(T_sx)=E_{[s}(\pi_0(x))$ (see
Lemma \ref{ctpath}). The part ii) shows that
 \begin{align*}
   E_{[t}(\langle m,m\rangle_{0}-\langle m,m\rangle_{t})
   &=2E_{[t} \int_0^t \pi_s(\Gamma(T_sx,T_sx))ds
   \lel  2 \pi_t \int_0^t T_{t-s}(\Gamma(T_sx,T_sx))ds\\
   &=  \pi_t(T_t|x|^2-|T_tx|^2) \pl .
   \end{align*}
Taking the supremum over  all $t$, we deduce the assertion. \qd

The BMO norm for the probabilistic model is closely related to the
associated Poisson semigroup.

\begin{prop}\label{bmo2} Let $(T_t)$ be a semigroup with a Markov dilation.
\begin{enumerate}
 \item[i)] Let $x\in N$ and $a>0$. Then
 \[ \|x\|_{BMO_c(P)} \lel \|\rho_a(x)\|_{BMO_c} \pl.\]
 \item[ii)]$
 \frac{1}{90} \|P^{br}\rho_a(x)\|_{bmo_c} \kl \sup_b \|\int_0^{\infty}
  P_{b+s}|P'_{s}|^2 \min(s,b) ds\|^{\frac12}
  \kl 2  \|P^{br}\rho_a(x)\|_{bmo_c}$.
 \item[iii)] Assume $\Gamma^2\gl 0$. Then
  \begin{align*}
   \frac{1}{90} \|(I\!-\!P^{br})\rho_a(x)\|_{bmo_c} &\le \!  \sup_b
   \|\int_0^{\infty}\!\!
  P_{b+s}\Gamma(P_{s}x,P_{s}x) \min(s,b) ds\|^{\frac 12 }
  \le  2 \|(I\!-\!P^{br})\rho_a(x)\|_{bmo_c}\!
   \pl. \end{align*}
\end{enumerate}
\end{prop}

\begin{proof} We recall that $\hat{E}_t(\rho_a(x))=
\hat{\pi}_{\ttt_a\wedge t}(Px)$ and $\rho_a(x)=\pi_{\ttt_a}(x)$.
Hence we get
 \begin{align*}
  \hat{E}_t(|\rho_a(x)|^2)-|\hat{E}_t(\rho_a(x))|^2
  &= \pi_{\ttt_a\wedge t}(P_{B_{\ttt_a\wedge t}}|x|^2-|P_{B_{\ttt_a\wedge
  t}}x|^2) \pl,
  \end{align*}
for $\ttt_a(\om)>t$. Thus in any case we have
 \[ {\rm ess}\sup_{\om} \|\hat{E}_t(|\rho_a(x)|^2)-|\hat{\pi}_{\ttt_a\wedge
 t}(Px)|^2\| \kl \sup_s \|\pi_{\ttt_a\wedge
 t}(P_s|x|^2-|P_sx|^2)\|\kl \|x\|_{BMO_c(P)}^2 \pl .\]
However, for $t=0$ we recall that $B_0(\om)=a$ almost everywhere.
This means $B_t=a+\tilde{B}_t$ where $\tilde{B}_t$ is a centered
brownian motion. Since $\limsup_t |\tilde{B}_t|/\sqrt{2t\log\log
t}=1$, we know that with probability $1$ the process $|\tilde{B}_t|$
exceeds $a$. Thus with probability $1$ the process $B_t$ hits $0$ or
$2a$. Hence with probability $\frac12$ the process hits $2a$ before
it hits $0$. Let us assume that $B_{t(\om)}(\om)=2a$ and
$B_s(\om)>0$ for $0<s<t(\om)$. By starting a new brownian motion at
$t(\om)$, we see with conditional probability $\frac12$ we have
$B_{t'(\om)}=4a$ for some $t(\om)<t'(\om)$ and $B_s(\om)>0$ for all
$t(\om)<s<t'(\om)$. By induction we deduce  that with probability
$2^{-n}$ the process $B_t$ hits $2^na$ before it hits $0$. Thus
given any $b>0$, we may choose $n$ such that $2^{n}a>b$. We see that
with positive probability there exists $t_n(\om)$ such that
$B_{t_n(\om)}=2^na$ and $B_s(\om)>0$ on $[0,t_n(\om)]$ and $B_s$ is
continuous. By continuity there exists $t(\om)\in
[t_n(\om),\ttt_a(\om)]$ such that $B_{t_{\om}}=b$. In particular,
 \[ \|\hat{E}_{t(\om)}(|\rho_a(x)|^2)-|\hat{\pi}_{\ttt_a\wedge
 t}(Px)|^2\| \lel \|\pi_{t(\om)}(P_{B_{t(\om)}}|x|^2-|P_{B_{t(\om)}}x|^2)\|
 \lel \|P_b|x|^2-|P_bx|^2\| \pl . \]
Taking the supremum over all $b$ yields i). For the proof of iii) we
first apply Lemma \ref{horcomp} and then Lemma \ref{Bkryest}. This
immediately yields the first inequality (after a concise review of
the involved constant for $\beta=\frac{2}{3}$). For the upper
estimate of this term, we recall that with positive probability
every value $b$ is hit. Then we start in \eqref{here} for a fixed
$b=B_t(\om)$. We use the monotonicity $\frac{P_{b+s}(z)}{b+s}\kl
\frac{P_t(z)}{t}$ and find
 \begin{align*}
 &\ez\int_0^{\ttt_b} T_s(\Gamma(P_{\tilde{B}_s}x,P_{\tilde{B}_s}x)) ds
 \lel \frac{1}{2} \int_0^{\infty} \int_{|b-s|}^{b+s} P_{t}\Gamma(P_sx,P_sx)
 dt ds \\
 &\gl   \frac{1}{2} \int_0^{\infty}
\frac{P_{b+s}\Gamma(P_sx,P_sx)}{b+s}
 \kla \int_{|b-s|}^{b+s} t dt \mer ds \lel
  \int_0^{\infty} \frac{P_{b+s}\Gamma(P_sx,P_sx)}{b+s} bs \pl
 ds \\
 &\gl \frac{1}{2} \int_0^{\infty} P_{b+s}\Gamma(P_{s}x,P_{s}x) \min(b,s)
 ds \pl .
 \end{align*}
The proof of ii) is similar but we only need $|P_tz|^2\le P_t|z|^2$
instead of $\Gamma^2\gl 0$. \qd

\begin{prop}\label{bmo3} Let $(T_t)$ be a semigroup of completely positive
selfadjoint maps and $(P_t)$ the associated Poisson semigroup. Let
$x\in \A$. Then
\begin{enumerate}
 \item[i)] $P_b|x|^2-|P_bx|^2=\int_0^\infty \int_{\max\{0,v-b\}}^{b}P_{b+2u-s}\hat{\Gamma
 }(P_sx,P_sx)duds $, and, assuming $\Gamma^2\gl 0$,
 \[ \frac14 \int_0^{\infty}  P_{b+s}\hat{\Gamma}(P_{s}x,P_sx)   \min(s,b)
 ds \kl P_b|x|^2-|P_bx|^2\kl 180 \int_0^{\infty} P_{\frac{b}{3}+s}
 \hat{\Gamma}(P_{s}x,P_{s}x) \pl \min(\frac b3,s) ds\pl ;\]
 \item[ii)] $\sup_b \|\int_0^{\infty} P_{b+s}|P'_sx|^2 \min(s,b)
 ds\| \kl 4 \pl \|x\|_{BMO_c(P)}^2$;
 \item[iii)] $\sup_b \|\int_0^{\infty} P_{b+s}\Gamma(P_sx,P_sx)\min(s,b)
  ds\| \kl  4\pl \|x\|_{BMO_c(P)}^2$
provided   $\Gamma^2\gl 0$.
\end{enumerate}
\end{prop}

\begin{proof} For the proof of i) we recall from \eqref{wies}
applied to $P_t$ that
 \[ P_b|x|^2-|P_bx|^2 \lel 2 \int_{0}^b  P_{b-s}\Gamma_{A^{1/2}}(P_sx,P_sx) ds \pl .\]
We recall from \cite{JR1} that
 \[ \Gamma_{A^{1/2}}(y,y) \lel \int_0^{\infty}
 P_t\Gamma(P_ty,P_ty) dt +
 \int_0^{\infty}
 P_t|P_t'y|^2  dt \]
holds for $y\in \A$. Combining these equations we obtain with a
change of variables ($v=s+t$, $u=t$)
 \begin{align}
  &P_b|x|^2-|P_bx|^2 \lel  2 \int_{0}^b \int_0^{\infty}
  P_{b-s+t}\hat{\Gamma}(P_{s+t}x,P_{s+t}x) dtds  \nonumber \\
  &\lel
  2 \int_{0}^{\infty} \int_{\max\{0,v-b\}}^{v}
  P_{b-v+2u}\hat{\Gamma}(P_{v}x,P_{v}x) du dv \label{start}
\end{align}
We apply $\hat{\Gamma}^2\gl$ and monotonicity \ref{monot} and split
the integral
 \begin{align*}
  & 2 \int_{0}^{\infty} \int_{\max\{0,v-b\}}^{v}
  P_{b-v+2u}\hat{\Gamma}(P_{v}x,P_{v}x) du dv\\
  &\gl 2 \int_{0}^{\infty} \int_{\max\{0,v-b\}}^{v}
  P_{b+2u}\hat{\Gamma}(P_{2v}x,P_{2v}x) du dv\\
  &\gl 2 \int_{0}^{\infty} \kla \int_{\max\{0,v-b\}}^v \frac{b+2u}{b+2v}
   du\mer  P_{b+2v}\hat{\Gamma}(P_{2v}x,P_{2v}x)  dv\\
 &\lel   \int_{0}^{\infty}  \frac{2bv+4v^2-2b\max\{0,v-b\}-4\max\{0,v-b\}^2}{2(b+2v)}
   P_{b+2v}\hat{\Gamma}(P_{2v}x,P_{2v}x)  dv\\
  &\gl   \int_{0}^{b}
   P_{b+2v}\hat{\Gamma}(P_{2v}x,P_{2v}x)  vdv
  +  \int_{b}^{\infty} \frac{4bv}{2(b+2v)}
    P_{b+2v}\hat{\Gamma}(P_{2v}x,P_{2v}x)  dv
    \\
  &\gl \frac12  \int_{0}^{b} P_{b+2v} \hat{\Gamma}(P_{2v}x,P_{2v}x)  2vdv
  +\frac{1}{2} b  \int_{b}^{\infty} P_{b+2v}\hat{\Gamma}(P_{2v}x,P_{2v}x)
  dv \\
  &\gl \frac12 \int_0^{\infty} P_{b+2v}\hat{\Gamma}(P_{2v}x,P_{2v}x)  \min(2v,b)
  dv\pl .
 \end{align*}
Without $\Gamma^2\gl 0$ we  only obtain
 \[ P_b|x|^2-|P_bx|^2 \gl
 \frac12 \int_0^{\infty}  P_{b+2v} |P'_{2v}x|^2   \min(2v,b) dv
 \lel \frac14 \int_0^{\infty}  P_{b+v} |P'_{v}x|^2   \min(v,b) dv
 \pl .\]
This yields iii) and iv). To complete the proof of i) we start with
\eqref{start} and the $\Gamma^2$ condition:
 \begin{align*}
  &P_b|x|^2-|P_bx|^2 \lel 2 \int_{0}^{\infty} \int_{\max\{0,v-b\}}^{v}
  P_{b-v+2u}\hat{\Gamma}(P_{v}x,P_{v}x) du dv \\
 &\kl 2 \int_{0}^{\infty} \int_{\max\{0,v-b\}}^{v}
  P_{b-\frac{v}{3}+2u}\hat{\Gamma}(P_{\frac{v}{3}}x,P_{\frac{v}{3}}x) du dv \\
 &\lel  2 \int_0^b \int_{0}^{v}  P_{b-v+2u}\hat{\Gamma}(P_{v}x,P_{v}x) du dv
 + 2 \int_b^{\infty} \int_{v-b}^{v}  P_{b-v+2u}\hat{\Gamma}(P_{v}x,P_{v}x) du dv
 \lel I+II \pl .
 \end{align*}
For $v\gl b$ we have
 \[  \frac{b+v}{3} \kl  b-\frac{v}{3}+2u \kl \frac{5}{3}(b+v) \pl. \]
Thus monotonicity  implies
 \begin{align*}
  II &\kl 2 \int_{0}^{\infty} \int_{v-b}^{v}
  P_{b-\frac{v}{3}+2u}\hat{\Gamma}(P_{\frac{v}{3}}x,P_{\frac{v}{3}} x) du dv
 \kl   10b  \int_{b}^{\infty}
 P_{\frac{b+v}{3}}\hat{\Gamma}(P_{\frac{v}{3}}x,P_{\frac{v}{3}}x) dv \\
 &= 90   \int_{\frac b3}^{\infty}
 P_{\frac{b}{3}+s}\hat{\Gamma}(P_{s}x,P_sx)
 \min(s,\frac{b}{3}) \pl
  ds\pl .
 \end{align*}
In the range $v\le b$ and $0\le u\le v$  we also have
 \[ \frac{b+v}{3}\kl b+2u-\frac{v}{3} \kl \frac{5}{3}(b+v) \pl .\]
Again by monotonicity  and with $s=\frac{v}{3}$ we obtain
 \begin{align*}
 I &\le 10  \int_0^{b}  P_{\frac{b+v}{3}}\hat{\Gamma}(P_{\frac{v}{3}}x,P_{\frac{v}{3}}
 x) \pl v dv
 \lel 90 \int_0^{\frac b3}  P_{\frac{b}{3}+s}\hat{\Gamma}(P_{s}x,P_{s}
 x) \pl s  ds \pl.
 \end{align*}
This yields
 \begin{align*}
  P_b|x|^2-|P_bx|^2 &\le  180 \int_0^{\infty} P_{\frac{b}{3}+s}
 \hat{\Gamma}(P_{s}x,P_{s}x) \pl \min(\frac{b}{3},s) ds \pl .\qedhere
 \end{align*}\qd

We will now study different BMO norms motivated by the expressions
above, namely
 \for
  \|x\|_{BMO_c(\hat{\Gamma})} &=&
 \sup_b \|P_b\int_0^{b} \hat{\Gamma}(P_sx,P_sx)
 sds\|_{\infty}^{\frac12} \pl ,\\
 \|x\|_{BMO_c^*(T)} &=&  \sup_t \|T_t|x-T_tx|^2\|^{1/2}
 \pl .
 \mel
The second norm has been introduced in \cite{Mei}, motivated by the
expression
 \[ \|f\|_{BMO_1} \lel  \sup_z P_z(|f-f(z)|)  \pl .\]
We also noticed that it was studied in the commutative case in
\cite{Duli2}. With respect to $\| \pl \|_{BMO_1}$ it is easy to show
that the conjugation operator is bounded from $L_{\infty}$ to BMO.
Here $f(z)$ gives the value of the harmonic extension in the
interior of the circle (see \cite{Ga}). This means in $f-f(z)$,
$f(z)$ is considered as a constant function. In some sense $x-P_tx$
is similar, but clearly $P_tx$ still is a function, even when $P_tx$
is the Poisson integral of $f$.

\begin{lemma}\label{lem1} Let $(T_t)$ be a semigroup satisfying $\Gamma^2\gl
 0$. Then
\[
 \frac14  \pl ||x||_{BMO_c(\hat \Gamma)}^2\kl  ||\sup_t\int_0^\infty P_{s+t}\hat{\Gamma }%
 (P_sx,P_sx)\min(s,t) ds||_\infty ^2\kl
  32 \pl ||x||_{BMO_c(\hat \Gamma)}^2.
 \]
\end{lemma}

\begin{proof} For the first estimate we note that due to $\Gamma^2\gl
0$ we have
 \begin{align*}
  & \|\int_0^tP_t\hat{\Gamma }(P_vx,P_vx)vdv \|_\infty
   \kl \|\int_0^tP_{\frac v2+t}\hat{\Gamma }(P_{\frac v2}x,P_{\frac v2}x)vdv\|_\infty
   \lel  4 \|\int_0^{\frac t2}P_{s+t}\hat{\Gamma}(P_sx,P_sx)sds \|_\infty
  \\
 &\kl  4 \|\int_0^\infty P_{s+t}\hat{\Gamma }(P_sx,P_sx) \min(s,t) ds \|_\infty .
 \end{align*}
For the other argument we use a dyadic decomposition. Indeed,
according to Proposition \ref{monot}, we have
 \[
 \frac{2^ntP_{s+t}}{s+t}\hat{\Gamma }(P_sx,P_sx)\leq P_{2^nt}\hat{%
 \Gamma }(P_sx,P_sx)
 \]
for $s\geq 2^nt.$ This implies
  \begin{align*}
 & \frac12 \int_0^\infty P_{s+t}\hat{\Gamma }(P_sx,P_sx) \min(s,t) ds  \kl
 \int_0^\infty P_{s+t}\hat{\Gamma }(P_sx,P_sx) \frac{st}{s+t} ds  \\
 &\lel \int_0^{2t}P_t\hat{\Gamma}(P_sx,P_sx)\frac{st}{s+t}ds+\sum_{n=1}^\infty \frac1{2^n}\int_{2^nt}^{2^{n+1}t}\frac{2^ntP_{s+t}}{s+t}\hat{\Gamma }(P_sx,P_sx)sds \\
  &\kl \int_0^{2t}P_t\hat{\Gamma}(P_sx,P_sx)sds+\sum_{n=1}^\infty \frac1{2^n}\int_{2^nt}^{2^{n+1}t}P_{2^nt}\hat{\Gamma }(P_sx,P_sx)sds \\
 &\kl  \int_0^{2t}P_t\hat{\Gamma}(P_sx,P_sx)sds+\sum_{n=1}^\infty \frac1{2^n}\int_0^{2^{n+1}t}P_{2^nt}\hat{\Gamma }(P_sx,P_sx)sds.
 \end{align*}
However, we can replace $2t$ by $t$ using $\Gamma^2\gl 0$ and Lemma
\ref{monot}:
 \begin{align*}
  &\int_0^{2t}P_t\hat{\Gamma }(P_sx,P_sx)sds
  \kl \int_0^{2t}P_{t+\frac s2}  \hat{\Gamma }(P_{\frac s2}x,P_{\frac s2} x)sds
  \lel 4  \int_0^{t} P_{t+v}  \hat{\Gamma }(P_vx,P_vx) vdv\\
  &\le  8 \int_0^{t} P_{t}  \hat{\Gamma }(P_vx,P_vx) vdv
  \pl.\end{align*}
Applying this argument for every $2^{n+1}t$, we deduce the
assertion. \qd

The next observation is true for arbitrary semigroups $(T_t)$.

\begin{prop}\label{gp} Let $(T_t)$ be a semigroup of completely positive
maps. Then
\begin{enumerate}
 \item[i)] $\|T_sx\|_{BMO_c(T)}\kl \|x\|_{BMO_c(T)}$ for all $s>0$
 and $x\in N$;
 \item[ii)] $\|x\|_{BMO_c^*(T)}\kl 2\pl \|x\|_{BMO_c(T)}+\sup_t
 \|T_tx-T_{2t}x\|_{\infty}^{1/2}$ for all $x\in N$.
 \end{enumerate}
\end{prop}

\begin{proof} Let us start with i) and the pointwise estimate
  \[ 0\kl T_{t}|T_sx|^2-|T_{t+s}x|^2\kl T_{t+s}|x|^2-|T_{t+s}x|^2
  \pl .\]
By definition of the $BMO_c(T)$ norm  this implies
 \begin{align*}
  \|T_sx\|_{BMO_c(T)}
 &=  \sup_t \| T_{t}|T_sx|^2-|T_{t+s}x|^2\|^{\frac12}_\infty \kl
  \sup_t \|T_{t+s}|x|^2-|T_{t+s}x|^2||^{\frac12}_\infty \kl
 \|x\|_{BMO_c(T)} \pl .
 \end{align*}
For the proof of ii), we fix  $t>0$ and use the triangle inequality
(see Lemma \ref{bmo1}):
\begin{align*}
 &\|T_{t}|x-T_tx|^2\|_\infty
 \kl \|T_{t}|x-T_tx|^2-|T_t(x-T_tx)|^2\|_\infty +
  \|\p  |T_t(x-T_tx)|^2\|_\infty \\
  &\le \|x-T_tx\|_{BMO_c(T)}^2+\||T_t(x-T_tx)|^2\|_\infty\\
  &\le  2\|x\|_{BMO_c(T)}^2+2\|T_tx\|_{BMO_c(T)}^2+\|T_t(x-T_tx)\|^2_\infty
 \end{align*}
We apply (i) and obtain
\begin{align*}
  \|T_{t}|x-T_tx|^2\|_\infty
 &\le  4\|x\|_{BMO_c(T)}^2+\|T_t(x-T_tx)\|^2_\infty.
 \end{align*}
Taking supremum over $t$ yields the assertion. \qd

Our next goal is to show that the $BMO_c(P)$-norm is in fact larger
than the $BMO_c^*(P)$-norm.

\begin{lemma}\label{aa} Let $a>1$. Then
 \[
   \sup_t \|P_tx-P_{at}x\| \kl  \sqrt{2}(1+{\rm \log }_{\frac 32}a)
  \|x\|_{BMO_c(\hat \Gamma_A)}.
\]
\end{lemma}

\begin{proof} For $t$ fixed, we have the
 \begin{align*}
 &|P_{3t}x-P_{2t}x|^2 \kl
  P_{\frac{3t}2}(|P_{\frac{3t}2}x-P_{\frac t2}x|^2)
  \lel P_{\frac{3t}2}(|\int_{\frac t2}^{\frac{3t}2} P_s'x ds|^2)
 \\
 & \kl   P_{\frac{3t}2}(t\int_{\frac t2}^{\frac{3t}2}|P'_sx|^2ds)
  \kl  2P_{\frac{3t}2}(\int_{\frac t2}^{\frac{3t}2}|P'_sx|^2sds)
 \kl
  2P_{\frac{3t}2}(\int_0^{\frac{3t}2}|P'_sx|^2sds).
\end{align*}
This implies in particular that
 \[ \sup_t \| P_tx-P_{\frac{3t}2}x\|_\infty \kl
 \sqrt{2} \|x\|_{BMO_c(\hat \Gamma)}.
 \]
For  $1<a\leq \frac 32$, choose $b\geq 0$ such that
$\frac{a-b}{1-b}=\frac 32$. Then we obtain
  \begin{align}
  \|\p  |P_tx-P_{at}x|^2 \|_\infty  &\le \| P_{bt}|P_{(1-b)t}(x)-P_{\frac
 32(1-b)t}(x)|^2 \|_\infty   \nonumber \\
   &\le  \| \p |P_{(1-b)t}(x)-P_{\frac 32(1-b)t}(x)|^2 \|_\infty
   \kl 2 \pl \|x\|_{BMO_c(\hat \Gamma)}^2 \pl .
 \end{align}
We deduce
 \begin{equation}
 \|P_t(x)-P_{at}(x) \|_\infty \kl  \sqrt{2} \pl \|x\|_{BMO_c(\hat \Gamma_A)} \label{sqrt2}
 \end{equation}
for any $1<a\le \frac32$. Consider now $a>\frac 32$. Let  $n$ be the
integer part of $\log _{\frac 32}a$.  We may use a telescopic sum
 \[ P_tx-P_{at}x \lel (P_tx-P_{\frac{3t}2}x)+(P_{\frac{3t}2}x-P_{_{\frac 32\frac{3t}2}}x)+\cdots (P_{(\frac
 32)^nt}x-P_{at}x)\pl .
 \]
We apply \eqref{sqrt2} for every summand. Then the triangle
inequality implies the assertion. \qd

The careful reader will have observed that we need one extra
estimate to complete the cycle.

\begin{prop}\label{gstar} Let $(T_t)$ be a semigroup satisfying
$\Gamma^2\gl 0$. Then
 \[ \|x\|_{BMO_c(\hat{\Gamma})}\kl  \frac{18}{3-\sqrt{7}} \pl \| x\|_{BMO_c^*(P)} \pl .\]
\end{prop}

\begin{proof} We fix $x,t$  and split $x=(x-P_{4t}x)+P_{4t}x$.
Then we have
 \begin{align*}
 & \| \int_0^\infty P_{s+t}\hat{\Gamma }(P_sx,P_sx) \frac{st}{s+t}
 ds \|_\infty ^{\frac 12} \kl
  \|\int_0^\infty P_{s+t}\hat{\Gamma }(P_s(x-P_{4t}x),P_s(x-P_{4t}x))\frac{st}{s+t} ds\|_\infty ^{\frac
  12}\\
 &\quad \quad\quad\quad + \|\int_0^\infty P_{s+t}\hat{\Gamma }(P_{s+4t}x,P_{s+4t}x) \frac{st}{s+t}ds\|_\infty ^{\frac 12}
 \end{align*}
For the first term we may apply Proposition \ref{bmo3}i)  for
$x'=x-P_{4t}x$ and obtain with $\frac{st}{s+t}\kl \min(s,t)$ that
 \[ \|\int_0^\infty P_{s+t}\hat{\Gamma }(P_s(x-P_{4t}x),P_s(x-P_{4t}x))\min(s,t) ds\|_\infty ^{\frac  12}
 \kl 2 \pl \|P_t|x-P_{4t}x|^2\|^{1/2} \pl. \]
The last term can be estimated by the $BMO_c^*$-norm using the
triangle inequality from Lemma \ref{bmo1} as follows
 \begin{align}
 &\|P_t|x-P_{4t}x|^2\|^{1/2} \kl \|P_t|x-P_{t}x|^2\|^{1/2}+
 \|P_t|P_tx-P_{2t}x|^2\|^{1/2}+\|P_t|P_{2t}x-P_{4t}x|^2\|^{1/2}  \nonumber \\
  &\le \|P_t|x-P_{t}x|^2\|^{1/2} + \|P_{2t}|x-P_{t}x|^2\|^{1/2}+\|P_{3t}|x-P_{t}x|^2\|^{1/2}
  \kl 3\|x\|_{BMO_c^*(P)} \label{finish}
 \end{align}
For the tail estimate we apply again the $\Gamma^2$ condition:
 \begin{align*}
 &\|\int_0^\infty P_{s+t}\hat{\Gamma }(P_{s+4t}x,P_{s+4t}x)
  \frac{st}{s+t} ds\|_\infty ^{\frac 12} \kl
 \|\int_0^\infty P_{s+t+3t} \hat{\Gamma}(P_{s+t}x,P_{s+t}x)
 \frac{st}{s+t}  ds \|_\infty^{\frac 12}\\
 &\lel  \|\int_t^\infty P_{s+3t}\hat{\Gamma }(P_sx,P_sx)
  \frac{(s-t)t}{s} ds\|_\infty^{\frac 12} \pl.
 \end{align*}
For $s \gl t$  we  consider the function $f(s)\lel
\frac{(s-t)(s+3t)}{3s^2}=\frac13+\frac{t}{s}[\frac23-\frac ts]$.
Note that $f(t)=0$ and for $t\le s\le \frac32 t$ we have $f(s)\le
\frac13$. For $s\gl \frac32 t$ we $f(s)\le \frac13+\frac49=\frac79$.
Thus in any case
 \[ \frac{(s-t)t}s\leq \frac 79\frac{3st}{s+3t}. \]
Taking supremum over $t,$ we obtain
 \begin{align*}
 &\sup_t \|\int_0^{\infty}
 \int_0^\infty P_{s+t}\hat{\Gamma }(P_sx,P_sx)\frac{st}{s+t} ds\|_\infty^{\frac
 12}\\
 &\kl 2 \sup_t  \|P_t|x-P_{4t}x|^2\|^{1/2}
  +
 \frac{\sqrt{7}}{3} \sup_t \|\int_0^{\infty}
 \int_0^\infty P_{s+3t}\hat{\Gamma }(P_sx,P_sx)\frac{s3t}{s+3t} ds\|_\infty^{\frac
 12} \pl .
 \end{align*}
This implies
 \begin{align*}
 &(1-\frac{\sqrt{7}}{3})
 \sup_t \|\int_0^{\infty} \int_0^\infty P_{s+t}\hat{\Gamma }(P_sx,P_sx)\frac{st}{s+t} ds\|_\infty^{\frac  12}
 \kl 2 \pl \sup_t  \|P_t|x-P_{4t}x|^2\|^{\frac12}\\
 &\kl 6 \pl \|x\|_{BMO_c^*(P)} \pl .
 \end{align*}
We may replace $\frac{st}{s+t}$ by $\min(s,t)$ with an additional
factor $2$. Hence the assertion follows from Lemma \ref{lem1}.\qd

\begin{theorem} Let $(T_t)$ be semigroup of completely positive
maps satisfying $\Gamma^2\gl 0$. Then the norms $\|\pl
\|_{BMO_c(P)}$, $\|\pl \|_{BMO_c^*(P)}$ and $\|\pl
\|_{BMO_c(\hat{\Gamma})}$ are all equivalent on $\A$.
\end{theorem}

\begin{proof} According to Proposition \ref{bmo3} we know that
 \[  \sup_t \|\int_0^\infty P_{s+t}\hat{\Gamma }(P_sx,P_sx)\min(s,t)ds\|_\infty^{\frac 12}
 \sim_{180} \|x\|_{BMO_c(P)} \pl .\]
Then Lemma \ref{lem1} implies that $\|\pl \|_{BMO_c(P)}$ and
$\|\pl\|_{BMO_c(\hat{\Gamma})}$ are equivalent. Proposition
\eqref{gstar} provides the upper estimate of $\|\pl
\|_{BMO_c(\hat{\Gamma})}$ against $\|\pl\|_{BMO_c^*(P)}$.
Conversely, we deduce from Proposition \ref{gp}, Lemma \ref{aa},
Lemma \ref{lem1} and  Proposition \ref{bmo3} i) that
 \begin{align*}
 &\|x\|_{BMO_c^*(P)}\kl 2 \|x\|_{BMO_c(P)}+ \sup_t
 \|P_tx-P_{2t}x\| \\
 &\le 2 \|x\|_{BMO_c(P)} + \sqrt{2}(1+\log_{\frac32}2)
 \|x\|_{BMO_c(\hat{\Gamma})} \\
 &\le 2 \|x\|_{BMO_c(P)} + 2\sqrt{2}\pl 2\sqrt{6}\pl \|x\|_{BMO_c(P)}
 \lel (2+8\sqrt{3}) \|x\|_{BMO_c(P)} \pl .
 \end{align*}
Thus all the norms are equivalent on $\A$. \qd

We conclude this section with two results on interpolation which
show that the BMO spaces are indeed a good endpoints. Let us define
 \[ \|x\|_{BMO(T)} \lel \max\{\|x\|_{BMO_c(T)},\|x^*\|_{BMO_c(T)}\}
 \pl \]
and the space $BMO(T)$ as the completion of $N$ with respect to that
norm.

\begin{theorem}\label{intpol} Let $(T_t)$ be a semigroup of completely positive
maps with a Markov dilation and $\Gamma^2\gl 0$. Then
 \begin{enumerate}
  \item[i)] $[BMO(T),L_p(N)]_{\frac1q} \lel L_{pq}(N)$;
  \item[ii)] $[BMO(P),L_p(N)]_{\frac1q} \lel L_{pq}(N)$.
 \end{enumerate}
holds for $1\le p<\infty$, $1<q<\infty$.
\end{theorem}

\begin{proof} For both proofs we note that the trivial inclusion
$N\subset BMO$ implies
 \[ L_{pq}(N)\subset [BMO(N),L_p(N)]_{\frac{1}{q}} \pl .\]
For the converse in i) we  consider the norm
  \begin{eqnarray} \label{pes}
  \|\pi_0(x)\|_{L_p^cmo(N)} \lel \|\sup_t \pi_t(T_t|x|^2-|T_tx|^2)\|_{L_{p/2}}^{1/2} \pl .
  \end{eqnarray}
Since $\Gamma^2\gl 0$, we know that the reversed martingale
$m_t(x)=\pi_t(T_t(x))$ has continuous path and therefore
 \[ \|\pi_0(x)\|_{L_p^cMO}
 \lel \lim_{|\si|,\U} \|\sup_j
 E_{t_j]}(|\pi_0(x)-m_{t_j}(x)|^2)\|_{p/2}
 \lel \|x\|_{L_p^cmo} \]
holds for every ultrafilter on the set of partitions. For a fixed
partition  and a reversed martingale $(y_s)$ we introduce
 \[ \|y\|_{L_p^cMO(\si)}
 \lel \max\{\|\sup_j
 E_{t_j]}(|y_0-y_{t_j}|^2)\|_{p/2}^{1/2},
 \| \sup_j
 E_{t_j]}(|y_0^*-y_{t_j}^*|^2)\|_{p/2}^{1/2}\} \pl .\]
It was shown in  \cite{JM} that
 \[ [L_pMO(\si),L_q(M)]_{\theta} \subset L_s(M)  \quad \mbox{with}\quad  \frac{1}{s}=\frac{1-\theta}{p}+\frac{\theta}{q} \]
and the constant $c(s)$  is uniformly on compact intervals of
$(1,\infty)$ and independent of $\si$. Thus for the norm
 \[ \|y\|_{L_pMO} \lel \lim_{\si,\U} \|x\|_{L_pMO(\si)} \]
we still have
 \[ [L_pMO,L_q(M)]_{\theta} \subset L_s(M)  \quad \mbox{with}\quad   \frac{1}{s}=\frac{1-\theta}{p}+\frac{\theta}{q}
 \pl .\]
However, on $\pi_0(A)$ we know that the $L_pMO$ and the $L_pmo$ norm
coincide and hence
 \[
 [\overline{\pi_0(\A)}^{\|\pl\|_{L_pmo}},\overline{\pi_0(\A)}^{\|\pl\|_{q}}]_{\theta}
 \subset L_s(M)\cap \overline{\pi_0(\A)}^{\|\pl\|_s} \pl .\]
Let us briefly indicate that $\pi_0(\A)$ is dense in $\pi_0(N)$
viewed as a subspace of $L_pMO$. Indeed, since the continuous
martingales are closed with respect to the $L_p$ norm,  and the
$L_p(M)$ majorizes the $L_p^cMO$ norm up to a constant $c(p)$, the
density follows from the density if $\A$ in $L_p(N)$. Finally, the
embedding $\pi_0$ of $L_s(N)$ in $L_s(\hat{M})$ is isometric and
therefore we have
 \[ [L_pmo(N),L_q(N)]_{\theta} \subset L_s(N) \]
with a constant $c(s)$ bounded on compact intervals. Sending $p\to
\infty$ and $\theta$ to $\frac{1}{v}$, we deduce
 \[ [BMO(T),L_q(N)]_{\frac{1}{v}} \subset L_{qv}(N) \pl .\]
For $P$ we can either work with the Markov dilation $\rho_a$ and
perform a similar argument, or we can observe that a Markov dilation
for $(T_t)$ produces a Markov dilation for $(P_t)$.\qd

In some applications it is worth while to work with the column space
$BMO_c(T)$ as completion of $N$ with respect to the $BMO_c(T)$ norm.

\begin{theorem}\label{intpol2} Let $(T_t)$ be a semigroup of completely positive
maps with a Markov dilation and $\Gamma^2\gl 0$. Then
 \begin{enumerate}
 \item[i)] $[BMO_c(P),L_2(N)]_{\frac{2}{p}}\subset H_p^c(P)$;
 \item[ii)] $[BMO_c(T),L_2(N)]_{\frac{2}{p}}\subset H_p^c(T)$.
 \end{enumerate}
\end{theorem}

\begin{proof} For the first assertion, we just recall that as in
the proof of Theorem \ref{intpol} we have
 \[ [L_p^cMO(M),L_2(M)]_{\theta}
 \subset  H_s^c(M)  \quad ,\quad  \frac{1}{s}\lel
 \frac{1-\theta}{p}+\frac{\theta}{2} \pl .\]
The same remark on the uniformity of constants applies. Note also
that
 \[ \| \pi_0(x)\|_{L_q^cMO(M)}\kl \| \pi_0(x)\|_{BMO_c(M)}
 \kl \|x\|_{BMO_c(P)} \pl .\]
Sending $p\to \infty$ and applying Theorem \ref{stein3}ii) we obtain
assertion ii). Assertion i) can be derived from ii) or adapting the
argument for $\rho_a$ instead of $\pi_0$.\qd

\section{Abstract semigroup theory}

In the theory of semigroups certain  tools from classical
Hardy-Littlewood theory are still available. In this section we will
recall the so-called Hardy-Littlewood-Sobolev theory which is
beautifully presented in \cite{VST}. Almost all (but not all) the
methods from the commutative theory apply in our setting. In
particular, we will prove the von Neumann algebra version of
\cite[Theorem II.5.2]{VST}. We refer to \cite{VST} for history and
credits. We are interested in the space
 \[ L_p^0(N) \lel ({\rm I}-\Pr)L_p(N) \]
of mean $0$ elements. Recall that $\Pr=\lim_{t\to \infty} T_t$ is
the orthogonal projection onto the kernel of $A$ and hence $({\rm
I}-\Pr)$ is a complete bounded (with cb-norm $\le 2$) on all
$L_p(N)$.

\begin{theorem}\label{coul} Let $(T_t)$ be a semigroup of completely positive
selfadjoint contractions on a von Neumann algebra $N$ with negative
generator $A$ and $n>2$. Let $L_p^0(N)$ be the space of mean $0$
elements. The following are equivalent
\begin{enumerate}
 \item[i)] $\|x\|_{2n/(n-2)}^2 \kl  C_1  \pl (x,Ax)  $ for all
 mean $0$ elements $x$,
 \item[ii)] $\|x\|_2^{2+4/n}\kl  C_2 \pl (x,Ax) \pl
 \|x\|_1^{4/n}$ for all mean $0$ elements $x$,
 \item[iii)] $\|T_t:L_1^0(N)\to L_{\infty}(N)\|\kl C_3 t^{-n/2}$.
 \end{enumerate}
\end{theorem}

\renewcommand{\theequation}{R$_n^{pq}$}
An important tool is  the family of conditions
 \begin{equation}\label{Rpq}
  \|T_t:L_p^0(N)\to L_q^0(N)\|\kl C t^{-\frac{n}{2}(1/p-1/q)}
  \quad,\quad 1\le p\le q\le \infty
\end{equation}
\renewcommand{\theequation}{\arabic{section}.\arabic{equation}}

The proof of the following result is verbatim  the same as in the
commutative case.
\begin{lemma}\label{RR} Let $T_t$ be a selfadjoint family of operators,
uniformly bounded on $L_p(N)$. Then $(R_n^{pq})$ holds for one pair
$1\le p<q\le \infty$ if and only if for all $1\le p\le q\le \infty$.
\end{lemma}

\begin{proof}[{\it Sketch of proof.}]  Let $p_1\le p<q$ and
$\frac{1}{p}=\frac{1-\theta}{p}+\frac{\theta}{q}$. Assume
$(R_n^{pq})$. Then we deduce from interpolation that
 \[ \|T_t(x)\|_p\kl C t^{-\frac{n}{2}(1/p-1/q)}\|x\|_p\kl
 C t^{-\frac{n}{2}(1/p-1/q)}
 \|x\|_{p_1}^{1-\theta}\|x\|_q^{\theta} \pl .\]
Thus $(R_n^{p_1q})$ holds with constant $C^{1/1-\theta}$. In
particular, $(R_n^{1q})$ holds. By duality we find
$(R_n^{q'\infty})$. Applying the argument again we get
$(R_n^{1\infty})$. Now, we show that $(R_n^{1,\infty})$ implies
$(R_n^{pq})$. Indeed, by complementation and interpolation we have
 \[ \|T_t:L_p^0(N)\to L_{\infty}\|\kl 2 \|T_t:L_{\infty}^0(N)\to
 L_{\infty}\|^{1-1/p}\|T_t:L_1^0(N)\to L_{\infty}(N)\|^{1/p} \pl. \]
This yields $(R_n^{p\infty})$. The same interpolation argument
implies $(R_n^{pq})$.\qd

In the following we will simply refer to the condition
\renewcommand{\theequation}{R$_n$}
 \begin{equation}\label{R0}
  \|T_t:L_1^0(N)\to L_\infty(N)\|\kl C t^{-\frac{n}{2}}
\end{equation}

Our next result requires a little bit more interpolation theory. We
recall  that two Banach spaces $A_0,A_1\subset V$ are injectively
embedded in a common topological vector space such that $A_0\cap
A_1$ is dense in $A_0$ and $A_1$. The unit ball of the space
$[A_0,A_1]_{\theta,1}$ is the convex hull of element $x$ in the
intersection satisfying
 \[ \|x\|_{A_0}^{1-\theta}\|x\|_{A_1}^{\theta} \le 1 \pl .\]
This implies that a linear operator $T:[A_0,A_1]_{\theta,1}\to X$
with values in a Banach space is continuous if
\renewcommand{\theequation}{$\theta,1$}
 \begin{equation}
  \|T(x)\|\kl C \|x\|_0^{1-\theta}\|x\|_1^{\theta} \quad ,\quad x\in X \pl .
  \end{equation}
The corresponding ``dual'' observation holds for the interpolation
space $[A_0,A_1]_{\theta,\infty}$. We recall that the norm of $x$ in
$[A_0,A_1]_{\theta,\infty}$ is less than $C$ if for every $t>0$ we
can decompose $x=x_0+x_1$ such that
\renewcommand{\theequation}{$\theta,\infty$}
 \begin{equation}
  \|x_0\|_0+t\|x_1\|_1\kl C t^{\theta} \pl .
 \end{equation}
\renewcommand{\theequation}{\arabic{section}.\arabic{equation}}
We will apply this for the scale of noncommutative Lorentz spaces
 \begin{equation}\label{mmar}
   L_{r,s}(N)\lel [L_{p,s_1}(N),L_{q,s_2}(N)]_{\theta,s} \quad
 ,\quad \frac{1}{r}\lel \frac{1-\theta}{p}+\frac{\theta}{q} \pl
 \end{equation}
which holds for all $1\le s_1,s_2\le \infty$ and $0<\theta<1$. We
refer to \cite{BL} for general information on interpolation theory
and to \cite{PX-hand} for  the translation to the noncommutative
setting. Since the space $L_p^0(N)$ is completely  complemented and
$L_{p,p}(N)=L_p(N)$, we may define
$L_{r,s}^0(N)=[L_p^0(N),L_q^0(N)]_{\theta,s}$ and then \eqref{mmar}
remains true for the spaces $L_{r,s}^0(N)$. The next argument is
adapted from \cite{Vo}. The conclusion is slightly weaker than in
the commutative situation. The key ingredient is the resolvent
formula
 \begin{equation}\label{inversf}
  A^{-z} \lel   \Gamma(z)^{-1} \int_0^{\infty} T_t \pl t^{z-1} dt
  \pl\quad \mbox{for $Re(z)>0$.}
  \end{equation}

\begin{lemma}\label{weak} Let $(T_t)$ be a semigroup of normal selfadjoint
contractions such that {\rm (}R$_{n}${\rm )} holds. Let $z\in \cz$
and $\al=Re(z)$.
\begin{enumerate}
\item[i)] Let $1\le p<s<q\le \infty$ and $z\in \cz$ with
$\al=\frac{n}{2}(\frac1s-\frac1q)$. Then
 \[ \| A^{-z}:L_{s,1}^0(N) \to L_q(N)\|\kl  C(\al,n) \pl
 . \]
 \item[ii)] Let $1\le p<r<\infty$ such that $\al\lel \frac{n}{2}(\frac1p-\frac1r)$.
Then
  \[ \|A^{-z}:L_p^0(N)\to L_{r,\infty}(N)\|\kl C(\al,n) \pl .\]
\end{enumerate}
\end{lemma}

\begin{proof} Ad i): We define $\al=Re(z)$,
$\frac1r=\frac1p-\frac1q$ and  $\theta=\frac{2r\al}{n}$. Let $x$ be
an element in $L_q^0(N)$ and $b>0$. In combination with $(R_n)$, we
deduce from \eqref{inversf} that
 \begin{align*}
 &|\Gamma(z)| \|A^{-z}(x)\|_{q} \lel  \|\int_0^{\infty} T_t(x) t^{z-1}
 dt
 \|_q \kl
  \int_0^b \|T_t(x)\|_{q} \pl t^{\al-1} dt +
  \int_b^{\infty} C t^{-n/2r} \|x\|_p  \pl t^{\al-1} dt \\
  &\le \al^{-1} b^{\al} \pl \|x\|_q + C
  (\frac{n}{2r}-\al)^{-1}
 b^{\al-\frac{n}{2r}} \pl \|x\|_p\pl .
 \end{align*}
We choose $b^{\frac{n}{2r}}\lel \frac{\|x\|_p}{\|x\|_q}$. (For
$\|x\|_q=0$ there is nothing to show.) This yields
 \[ \|A^{-z}(x)\|_{q}\kl |\Gamma(z)|^{-1} K' \frac{n}{\al(n- 2r\al)}
  \|x\|_q^{1-\frac{2r\al}{n}}  \|x\|_p^{\frac{2r\al}{n}} \pl
 .\]
The assertion follows from equation $(\theta,1)$ and
$L_{s,1}^0(N)=[L_q^0(N),L_p^0(N)]_{\theta,1}$. Note also that
$1/s=(1-\theta)/q+\theta/p=1/q+\theta/r=1/q+2\al/n$. For the proof
of ii) we define $\theta=1-\frac{2p\al}{n}$. Assume that $x\in
L_p^0(N)$ and decompose $\Gamma(z) A^{-z}x=x_0+x_1$ where
 \[ x_0 \lel \int_b^{\infty}  T_t(x) t^{z-1} dt
  \quad ,\quad
  x_1 \lel \int_0^b T_t(x) t^{z-1} dt
  \pl .\]
As above we deduce from $R_n^{p,\infty}$  and the assumption
$2p\al<n$ that
 \[ \|x_0\|_{\infty} \kl \frac{2Cp}{n-2p\al}
 b^{\al-\frac{n}{2p}}  \pl \|x\|_p
   \pl. \]
On the other hand, we have $\|x_1\|_p\le \frac{b^{\al}}{\al}
\|x\|_p$. For fixed $t>0$ we choose $b$ such that $b^{-n/2p}=t$.
Then
 \[ \|x_0\|_{\infty}+t\|x_1\|_{p} \kl K' \|x\|_p \pl b^{\al-\frac{n}{2p}}
 \lel K' \|x\|_p \pl  t^{1-\frac{2p\al}{n}}\pl .\]
Thus we have verified condition $(\theta,\infty)$ and the assertion
follows from the real interpolation method
$L_{r,\infty}(N)=[L_{\infty}(N),L_p(N)]_{\theta,\infty}$.\qd

As an immediate application of the Marcinkiewicz interpolation
theorem (in the form of \eqref{mmar}), we can remove the Lorentz
spaces from the conclusion. 

\begin{cor}\label{strong} Let $(T_t)$ be a semigroup of normal selfadjoint
contractions such that {\rm (}R$_{n}${\rm )} holds. In the type
$III$ case we assume in addition that $T_t$ commutes with the
modular group of a normal faithful state. Let $z\in \cz$ and
$\al=Re(z)$. Then
 \[ \|A^{-z}:L_p^0(N)\to L_q^0(N)\|\kl C(\al) \]
holds for all $1<p<q<\infty$ such that $\al\lel
\frac{n}{2}(\frac1p-\frac1q)$.
\end{cor}

\begin{proof}[{\it Proof of Theorem \ref{coul}.}] For the proof of the implication iii) $\Rightarrow$ i) we choose
$z=\frac12$ and $p=2$. Note that $(x,Ax)=\|A^{1/2}x\|_2^2$. Since
$(R_n)$ is satisfied, we obtain $1=n(1/2-1/q)$, i.e. $q=2n/(n-2)$.
The implication i) $\Rightarrow$ ii) follows from
 \[ \|x\|_2\kl \|x\|_{_{\frac{2n}{n-2}}}^{^{\frac{n}{n+2}}}
 \|x\|_1^{\frac{2}{n+2}} \pl .\]
The implication ii) $\Rightarrow$ iii) follows verbatim as in
\cite[Theorem III.3.2]{VST}. One first shows $(R_n^{1,2})$ by
differentiation for a selfadjoint mean $0$ element $x$  using
$\frac{d}{dt}\|T_tx\|_2^2=-2Re(AT_tx,T_tx)$. Remark \ref{RR} implies
the assertion.\qd

For our applications we need compactness results of the operator
$A^{-\al}$ on $L_p^0(N)$. Our aim is to derive them from compactness
on $L_2^0(N)$. Let us consider the following conditions
\begin{enumerate}
 \item[gap$_c$)] The spectrum of $(1-\Pr)A$ on $L_2(N)$ is
 contained in $[c,\infty)$,
 \item[com$_{\p}$)] $A^{-1}$ is compact on $L_2^0(N)$.
 \end{enumerate}

\begin{prop}\label{gap} Let $A$ be a generator which satisfies
$gap_c$. Then
 \begin{enumerate}
 \item[i)] Let $Re(z)>0$. Then $A^{-z}$ is (completely) bounded on
 $L_p^0(N)$ for $1<p<\infty$,
 \item[ii)] $\|T_t:L_p^0(N)\to L_p^0(N)\|_{cb}\kl 2 e^{-\frac{2tc}{p}}$ for
 all $2\le p<\infty$.
 \end{enumerate}
\end{prop}

\begin{proof} First we note that gap$_c$) means $A\gl c$ on the
Hilbert space $L_2^0(N)$ and hence
 \[ \|e^{-tA}:L_2^0(N)\to L_2^0(N)\|\kl e^{-tc} \pl .\]
Let $2\le p<\infty$. Since $L_p^0(N)$ forms an interpolation scale,
we deduce
 \[ \|T_t:L_p^0(N)\to L_p^0(N)\|\kl 2 \|T_t:L_2^0(N)\to
 L_2^0(N)\|^{\frac{2}{p}} \|T_t:N\to N\|^{1-2/p} \kl 2e^{-2tc/p} \pl .\]
This shows ii) and the same argument for $T_t\ten id$ provides the
cb-estimate. For the proof of i) we assume $2<p<\infty$. Then
\eqref{inversf} implies with $\al=Re(z)$ that
 \begin{align*}
 \|A^{-z}:L_p^0(N)\to L_p^0(N)\|_{cb}
 &\le 2|\Gamma(z)|^{-1}  \int_0^{\infty} e^{-2ct/p} s^{\al-1}
 ds<\infty \pl .
 \end{align*}
Since $A$ is selfadjoint the same estimate holds on
$L_{p'}^0(N)$.\qd



The next Lemma allows to interpolate compactness (see \cite{Pie}).

\begin{lemma}\label{comint} Let $(A_0,A_1)$ be an interpolation
couple as above, $T:X\to A_0\cap A_1$ a linear map such that $T:A\to
A_0$ is bounded and $T:X\to A_1$ is compact. Then $T:X\to A_\theta$
is compact.
\end{lemma}

\begin{proof} Let us recall that $T:X\to Y$ is compact if and only
if the entropy numbers
 \[ e_k(T)\lel \inf\{ \eps : T(B_X)\subset \bigcup_{j=0}^{2^{k-1}}
 y_j+\eps B_Y\} \]
satisfy $\lim_k e_k(T)=0$. Here $B_X$, $B_Y$, is the unit ball of
$X$, $Y$, respectively. The infimum is taken over arbitrary points
in $y$. We recall from \cite{Pie} that
 \[ e_{k+j-1}(T:X\to A_\theta)
 \kl 2 e_k(T:X\to A_0)^{1-\theta} e_j(T:X\to A_1)^{\theta} \pl .\]
In particular, $e_k(T:X\to A_{\theta})\kl 2 \|T:X\to
A_0\|^{1-\theta} e_k(T:X\to A_1)^{\theta}$ still converges to $0$.
\qd

\begin{theorem}\label{commm}  Let $(T_t)$ be a semigroup of
selfadjoint, positive contractions on a finite von Neumann algebra
satisfying
 \[ \|T_t:L_1^0(N)\to L_{\infty}(N)\|\kl C t^{-n/2}\]
and such that $A^{-1}$ is compact on $L_2^0(N)$. Then
$A^{-z}:L_p^0(N)\to L_q^0(N)$ is compact for all $1\le p<q\le
\infty$ such that $\frac{2Re(z)}{n}>\frac1p-\frac1q$.
\end{theorem}

\begin{proof} By assumption $A^{-1}$ is bounded on $L_2^0(N)$ and hence we have a spectral gap.
Now, we consider $2<p<\infty$ and want to show that
$A^{-\al}:L_p^0(N)\to L_p^0(N)$ is compact for all $\al>0$.
According to Proposition \ref{gap} it suffices to consider
$\al<n/2p$. Define $1/q=1/p-2\al/n$. According to Corollary
\ref{strong} we know that $A^{-\al}:L_p^0(N)\to L_q^0(N)$ is
bounded. Since $A^{-\al}$ is compact on $L_2^0(N)$ we also know that
$A^{-\al}:L_p^0(N)\to L_2^0(N)$ is compact. We may write
$L_p^0(N)=[L_q^0(N),L_2^0(N)]_{\theta}$ where
$1/p=(1-\theta)/q+\theta/2$ and $0<\theta<1$. Hence Lemma
\ref{comint} implies that $A^{-\al}:L_p^0(N)\to L_p^0(N)$ is
compact. By duality we conclude that $A^{-\al}:L_p^0(N)\to L_p^0(N)$
is compact for all $\al>0$ and $1<p<\infty$.

Now, we consider $z=\al+is$ and assume that $1<p<q$. By our
assumption $2\al/n>1/p-1/q$. This allows us to find $1<s<p$, and
$\al_1>0$ such that $2(\al-\al_1)/n=1/s-1/q$. According to Lemma
\ref{weak}i)  we know that $A^{\al_1-z}:L_{s,1}^0(N)\to L_q^0(N)$ is
bounded. On the other hand $A^{-\al_1}:L_p^0(N)\to L_p^0(N)$ is
compact and the inclusion $L_p^0(N)\subset L_{s,1}^0(N)$ continuous.
Then
 \[ A^{-z}=A^{\al_1-z}A^{-\al_1}:
 L_p^0(N)\stackrel{A^{-\al_1}\mbox{\scriptsize compact}}{\longrightarrow} L_p^0(N)
 \subset L_{s,1}^0(N) \stackrel{A^{\al_1-z}}{\longrightarrow}
 L_q^0(N)\]
is the composition of a bounded operator and a compact operator,
hence itself compact.

In the case $1\le p<q<\infty$ we use the same argument  and find
$q<r<\infty$, a decomposition $A^{-z}=A^{-\al_1}A^{\al_1-z}$ such
that $A^{\al_1-z}:L_p^0(N)\to L_{r,\infty}^0(N)$ is continuous,
$A^{-\al_1}:L_q^0(N)\to L_q^0(N)$ is compact, and the inclusion
$L_{r,\infty}^0(N)\subset L_q^0(N)$ is continuous. Thus $A^{-z}$ is
compact. Finally for $p=1$ and $q=\infty$ we write
$A^{-z}=A^{-z/2}A^{-z/2}$ and make a pit stop at  $L_2^0(N)$.\qd

\section{Applications}

\subsection{Quantum metric spaces}

We recall from \cite{Rie} that a quantum metric space is given by a
$C^*$-algebra $C$,  a $^*$-subalgebra $\A$ and a norm $\tnorm\pl$ on
$\A$ such that
 \[ d_{\left\vert\xyspace\left\vert\xyspace\left\vert \pl
 \right\vert\xyspace\right\vert\xyspace\right\vert
}(\phi,\psi) \lel \sup\{ |\phi(a)-\psi(a)|: a\in \A, \tnorm a \le
1\} \] induces the weak$^*$ topology on the state space $S(C)$.  A
norm $\tnorm\pl $ is a Lipschitz norm if in addition
 \begin{equation}\label{Lip}
  \tnorm{ab} \kl \tnorm a\|b\|+\|a\|\tnorm b \pl .
  \end{equation}
In \cite{OR} a Lipschitz norm generating the weak$^*$ topology is
called a Lip-norm, i.e. quantum metric space require a Lip-norm
instead of Lipschitz norm.

\begin{lemma}\label{Lip2} Let $T_t$ be a unital completely positive semigroup on a von Neumann algebra $N$. Let
$-A$ be the generator  and  $\A$ be a (non-complete) $^*$-algebra
contained in the domain of $A$. Then
 \[ \|a\|_{\Gamma}\lel \max\{\|\Gamma(a,a)\|^{1/2},
 \|\Gamma(a^*,a^*)\|^{1/2}\} \]
and $\tnorm a =\|\Gamma(a,a)\|^{1/2}$ satisfy \eqref{Lip}.
\end{lemma}

\begin{proof} We recall from \cite{Prig}  that $H_N=\{\sum_i a_i\ten y_i: \sum_i
x_iy_i=0\}$ equipped wit the $N$-valued  inner product
 \[ \langle a_1\ten x_1,a_2\ten x_2\rangle
 \lel x_1^*\Gamma(a_1,a_2)x_2\]
defines a $N$-valued Hilbert module. Then $\delta(a)=a\ten 1-1\ten
a$ is a derivation, i.e.
 \[ \delta(ab)\lel ab\ten 1-1\ten ab
 \lel (a\ten 1) (b\ten 1-1\ten b)+(a\ten 1-1\ten a)(1\ten b)
 \lel (a\ten 1)\delta(b)+\delta(a)(1\ten b) \pl. \]
Since $T_t(1)=1$ we have $A(1)=0$ and
  \[ \Gamma(1,a)\lel 1A(a)+a^*A(1)-A(1a)\lel 0 \pl .\]
Hence $\langle \delta(a),\delta(a)\rangle= \Gamma(a,a)$. This
implies
 \begin{align*}
  \|\Gamma(ab,ab)\|^{1/2} \lel
  \|\delta(ab)\| \kl \|(1\ten a)\delta(b)\| +
  \|\delta(a)(1\ten b)\|
  \kl \|a\|\|\delta(b)\|+\|\delta(a)\|\|b\| \pl .
  \end{align*}
Recall that $\|\delta(a)\|=\|\Gamma(a,a)\|^{1/2}$. This also shows
 \[ \|\Gamma((ab)^*,(ab)^*)\|^{1/2}
  \lel \|\Gamma(b^*a^*,b^*a^*)\|^{1/2}
 \kl \|b^*\|\|\Gamma(a^*,a^*)\|^{1/2} +  \|\Gamma(b^*,b^*)\|^{1/2}
 \|a\| \pl .\]
Taking the maximum yields \eqref{Lip}.  \qd

We also need the following observation from \cite[Proposition
1.3]{OR}

\begin{lemma}\label{ro} Let $\tnorm \pl$ be a Lipschitz norm and $\si$ be a state. Then $(C,\A,\tnorm\pl
)$ is a quantum metric space  iff
 \[ \{ x\in \A\pl:\pl \tnorm a\le 1, \si(a)=0\} \]
is relatively compact in $C$.
\end{lemma}

\begin{theorem}\label{quant} Let $(T_t)$ be a completely positive semigroup of
selfadjoint maps on a finite von Neumann algebra $N$ such that
$\A\subset N$ is weakly dense and with a Markov dilation. Assume in
addition
 \begin{enumerate}
 \item[i)] $\ker(A)=\cz1 $ and
 $A^{-1}$ is compact on $L_2^0(N)=(I-\Pr)L_2(N)$,
 \item[ii)] $\|T_t:L_2^0(N)\to N\|\kl C t^{-n/4}$ for some $n>0$.
\end{enumerate}
Then
 \[ \|x\|_{\Gamma}\lel \max\{\|\Gamma(x,x)^{1/2}\|,
 \|\Gamma(x^*,x^*)^{1/2}\|\} \]
and $\tnorm{x}=\|\Gamma(x,x)\|^{1/2}$ define a quantum metric spaces
for the norm closure $C$ of $\A\subset N$.
\end{theorem}

\begin{proof} Let us recall that $\Pr$ is the projection onto
the kernel of the selfadjoint operator $A$. Thus i) implies in
particular that $\ker(A)=\cz 1$ and that $A$ has a spectral gap
 \[ c=\|A^{-1}:L_2^0(N)\to L_2^0(N)\| \pl .\]
Moreover, $\lim_{t\to \infty} T_t(x)=\tau(x)1$ and hence $L_p^0(N)$
is the closure of elements $x\in N$ such that $\tau(x)=0$. Let
$1<s<p<\infty$ such that $2<p$ and $\frac{2\al}{n}>\frac 1 p$.
According to Corollary \ref{commm} we know that
 \[ \{x\in L_p^0: \|A^{\al}x\|_p\le 1\} \subset L_{\infty}^0(N) \]
is relatively compact in $N$. Let $\delta>0$. Then we deduce from
\cite{JR1} and  Theorem \ref{hpp} that
 \begin{align*}
   \|A^{\frac12-\delta}x\|_p &\le  c(\delta)
 \|A^{\frac12}x\|_{H_{p}^c(T)}
  \kl c(\delta)c(p) \|\Gamma(x,x)^{1/2}\|_p \\
  &\le c(\delta)c(p) \|\Gamma(x,x)^{1/2}\|_\infty
  \lel c(\delta)c(p) \|\Gamma(x,x)\|_\infty^{1/2}
  \pl .
  \end{align*}
Hence we need $\frac12-\delta> \frac{n}{2p}$ which is satisfied for
$p>n$. Lemma \ref{ro} implies the assertion.\qd

\begin{rem} Let $M$ be a compact Riemannian manifold. Then
 \[ d(p,q) \lel \sup\{|f(x)-f(y)|\pl :\pl \| \p |\nabla f| \p \|_{\infty}\le 1\}
 \pl. \]
Moreover, $\Gamma(f,f)=|\nabla f|^2$. The condition ii) corresponds
to a Sobolov embedding theorem and Theorem \ref{quant} provides an
appropriate gradient norm in this context.\end{rem}

\subsection{Rapid decay and quantum metric spaces}

Let us recall that a finitely generated discrete group has
\emph{rapid decay (RD) of order $s$} if there exists an $s<\infty$
such that
 \[ \|x\|_{\infty}\kl C(s) k^s \|x\|_2  \]
holds for all linear combinations $x=\sum_{|g|=k} a_g\la(g)$. Here
$|\pl|$ is the word length function with respect to fixed number of
generators. The notion is, however, independent of that choice. We
refer to \cite{Jo1} for more information. The following observation
is closely related to the work of Rieffel and Ozawa \cite{OR}.

\begin{lemma}\label{dim} Let $G$ be a discrete, finitely generated
group with word length function $|\pl|$ and rapid decay of order
$s$. Let $\psi:G\to \rz$ be a conditionally negative function such
that
 \begin{equation}\label{kalp}
  \inf_{l(g)=k} \psi(g) \gl c_{\al} k^{\al}
  \end{equation}
for some $\al>0$. Then the operator
$T_t(\la(g))=e^{-t\psi(g)}\la(g)$ satisfies
 \[ \|T_t:L_2^0(N)\to N\|\kl C(s,\al) \pl t^{-\frac{2s+1}{2\al}} \pl .\]
\end{lemma}

\begin{proof} We consider a decomposition $x=\sum_k x_k$ such that
$x_k=\sum_{|g|=k} a_{g}\la(g)$ is supported by words of length $k$.
Note that $T_t(x_k)$ is still supported on words of length $k$ and
for such $g$ we have $e^{-t\psi(g)}\le e^{-tc_\al k^{\al}}$. Hence
we get
\begin{align*}
 \|T_tx\|&\le \sum_k \|T_tx\|_{\infty}
 \kl C(s) \sum_k k^s \|T_tx\|_2 \\
 &\le C(s) \sum_k k^s e^{-tc_{\al}k^{\al}} \|x_k\|_2 \kl
 C(s) (\sum_k k^{2s} e^{-2tc_{\al}k^{\al}})^{1/2} (\sum_k
 \|x_k\|_2^2)^{1/2}\pl.
 \end{align*}
Now it remains to estimate the sum via some calculus (i.e.
$y=2tc_{\al}x^{\al}$, $dy/y=\al dx/x$)
 \begin{align*}
 \sum_{k\gl 1} k^{2s} e^{-2tc_{\al}k^{\al}}
 &= e^{-2tc_{\al}} + 2^{2s} \int_{2c_{\al}}^{\infty} x^{2s+1}
 e^{-t2c_{\al}x^{\al}}\frac{dx}{x} \\
 &=
   e^{-2tc_{\al}}+ 2^{2s}\al^{-1} (2tc_{\al})^{-(2s+1)/\al}
 \int_1^{\infty} y^{\frac{2s+1}{\al}} e^{-y} \frac{dy}{y}  \pl.
 \end{align*}
Thus for $0<t\le 2$ we obtain
 \[ \|T_tx\|_{\infty}\kl C(s,\al) \pl  t^{-\frac{2s+1}{2\al}}
 \|x\|_2 \pl .\]
We recall that on $L_2^0(N)=\cz 1^{\perp}$ we have a spectral gap
$\psi(w)\gl c(\al) |w|^{\al}\gl c(\al)$ for all $w\neq 1$. Hence
$\|T_t:L_2^0(N)\to L_2^0(N)\|\le e^{-tc_{\al}}$. Hence for $t\gl 2$
we have
 \[ \|T_tx\|_{\infty} \lel \|T_1(T_{t-1}x)\|_{\infty}
 \kl C(s,\al) \pl  \|T_{t-1}x\|_2
 \kl C(s,\al)e^{c(\al)} e^{-tc_{\al}} \pl \|x\|_2 \pl .\]
The assertion follows.\qd

\begin{rem}{\rm In case of the free group and $\psi(g)=|g|$ we
have $\al=1$ and $s=1$. This yields the order  $t^{-3/2}$ and hence
property (R$_6$). According to Varopoulos' definition \cite{Vo} this
means dimension $d=n/2=3$, as predicted by P. Biane \cite{Bia}.}
\end{rem}

\begin{rem}\label{oorr}{\rm According to the work of Rieffel and Ozawa\cite{OR} hyperbolic groups satisfies
rapid decay with $s=1$ and $d=3$. }
\end{rem}

\begin{proof}[Proof of Corollary  \ref{hyp}] According to Lemma \ref{dim} the assumption ii) of Theorem
\ref{quant} are satisfied for $\A=\cz[G]$. Since $G$ is finitely
generated we know that that the span  $F_k$  of work of length $k$
are finite dimensional. By assumption the inverse of the operator
$A(\la(g))=\psi(g)\la(g)$ satisfies $\|A^{-1}:F_k\to F_k\|\le
c_{\al}^{-1}k^{-\al}$ and hence $A^{-1}$ is compact on $L_2^0(N)$.
This provides assumption i) and  Theorem \ref{quant} implies the
assertion. \qd

\begin{exam}{\rm  1) The most natural examples are cocompact lattices
$\Gamma\subset G$, where
 \[ G\in\{SO_0(n,1), SU(n,1)\} \pl .\]
Let us indicate that the assumptions are verified for $\al=1$.
Indeed, we first recall that $G$ acts on a hyperbolic space $X$ and
isometrically on the virtual boundary $\partial X$. Moreover, there
exists a quadratic form $Q$ on the boundary such that
 \[ \phi(d(x,y))\lel Q(\mu_x-\mu_y) \]
holds for all $x,y\in \partial x$. Here $d$ is the hyperbolic
distance and $\phi(r)$ behaves like $2\log\cosh(r)$ for large $r$.
This means $c_1r\le \phi(r)\le c_2r$. By the Milnor-Swarc Lemma (see
e.g. \cite{Roe}), we also know that for cocompact discrete lattice
the word length is quasi isometric to hyperbolic distance
 \[ c_1^{-1}l(g)\kl  d(gx_0,x_0) \kl c_2 l(g) \]
given by a fixed base point. This yields $s=1$. Hence we find
dimension $3$ in all of these cases.

2) The assumptions are satisfied for the free group in finitely many
 generators by the work of Haagerup \cite{Haa} (see also
\cite{many}).

3) Let $G_1$ and $G_2$ be two groups with rapid decay and
conditionally negative functions $\psi_1$, $\psi_2$ satisfying
\eqref{kalp} with $\al=\min(\al_1,\al_2)\le 1$. Then
$\psi(g,h)=\psi_1(g)+\psi_2(h)$ also satisfies \eqref{kalp}.
According to Jolissaint's  work \cite[Lemma 2.1.2]{Jo1}, the product
also has rapid decay. Thus $T_t(\la((g,h))=e^{-t\psi(g,h)}\la(g,h)$
defines a completely positive semigroup for which the assumptions of
Theorem \ref{quant} are also satisfied.

4) Let $(G_i,l_i,\psi_i)$ be groups with rapid decay and
conditionally length functions $\psi$ satisfying \eqref{kalp} with
parameter $k_{\al}$. According to \cite[Theorem 2.2.2]{Jo1} we know
that $(\ast_i G_i,\ast l_i)$ has property RD where
 \[ \ast l_i(w_1\cdots w_n)\lel \sum_i |w_i|_{l_i} \pl \]
here $w_j\in G_{i_j}$. Bo\'{z}ejko proved that $\psi_t(w_1\cdots
w_n)=e^{-t\sum_j\psi_{i_j}(w_j)}$ are still positive definite and
hence the free sum $\psi(w_1\cdots w_n)=\sum_j \psi_{i_j}(w_j)$ is a
conditionally negative definite function on $\ast_i G_i$ such that
 \begin{align*}
 \psi(w_1\cdots w_n) \lel \sum_j \psi_{i_j}(w_j) \gl  c(\al)
 \sum_j |w_j|^{\al} \gl c(\al) (\sum_j |w_j|)^{\al}
 \end{align*}
holds for $\al \le \min\{1,\al_j\}$. Hence the free product is again
a quantum metric space.} \end{exam}

\subsection{Torsion free ordered groups}

In this section we show that multipliers on $\zz$ can be used to
obtain result for torsion free ordered groups. Our main application
is the well-known Hilbert transform in this context of sub-diagonal
von Neumann algebras. Let us consider a discrete group $G$ with
normal divisors
 \[ G=G_0 \unrhd G_1 \unrhd \cdots \]
such that $\bigcap_i G_i=\{1\}$ and
 \begin{equation}\label{fil}
  G_i/G_{i+1} \lel \zz \pl .
  \end{equation}
It is very easy to see that if we were to have
$G_{i}/G_{i+1}=\zz^{n_i}$ that the sequence can be further refined
to satisfy \eqref{fil}. Our aim is to use Riesz transforms to show
the boundedness of the Hilbert transform for ordered groups. Let us
recall that in the situation above the cone of positive group
elements $P$ is given by
 \[ P \lel \{g\in G: g\in G_i\setminus G_{i+1} \mbox{ and }
 gG_{i+1} \gl 0 \} \pl .\]
Clearly, the integer $i$ is uniquely determined by $g$. Here $``\gl
0$'' is the usual  relation in $\zz$. We have $P\cup
P^{-1}=G\setminus \{1\}$. In the following we denote by
$(VN(G_i))_{i\gl 0}$ the reversed martingale filtration given by the
conditional expectation $E_i(\la(g))=1_{g\in G_i}\la(g)$.

\begin{defi} Let $\bigcup_k N_k\subset N$ be a weakly dense  martingale filtration. A \emph{tangent
dilation} is given by a von Neumann algebra $M$ and trace (state and
modular group) preserving homomorphisms $\pi_k:N_k\to M$, $\rho:N\to
M$ such that
 \begin{enumerate}
 \item[i)] The conditional expectation $E_{\rho}:M\to \rho(N)$ satisfies
 \[  \rho E_{k-1}\lel E_{\rho}\pi_k \]
 for all $k$.
 \item[ii)] The von Neumann algebras $M_k=\pi_k(N_k)$ are
 successively independent over $\rho(N)$.
 \end{enumerate}
\end{defi}

\begin{lemma}\label{rosen} Let $(N_k)$ be a martingale filtration and
$(\rho,(\pi_k)_k)$ a tangent dilation. Let $1<p<\infty$. Then the
map $d:L_p(N)\to L_p(M)$ given by
 \[ dx\lel \sum_k \pi_k(d_k(x)) \]
gives a  linear isomorphic embedding. In the limit cases $d$ is
bounded between the corresponding martingale $BMO$ and $H_1$ spaces.
\end{lemma}

\begin{proof} Let us first consider $p\gl 2$. We apply the
Rosenthal inequality and deduce that
 \begin{align*}
  \|dx\|_p &\kl cp \kla \kla \sum_k \|\pi_k(d_k(x))\|_p^p\mer^{\frac1p}
 +\|(\sum_k E_\rho(\pi_k(|d_k(x)|^2+|d_k(x)^*|^2))^{1/2}\|_p \mer \\
 &\le  cp \kla \kla \sum_k \|d_k(x)\|_p^p\mer^{\frac1p}
 +\|(\sum_k E_{k-1}(|d_k(x)|^2+|d_k(x)^*|^2))^{1/2}\|_p \mer
 \end{align*}
Therefore the Burkholder/Rosenthal inequality implies
 \[ \|dx\|_p \kl cp \pl c(p) \|x\|_p  \pl .\]
A similar  argument applies for the BMO norms. Indeed, we denote by
$\hat{E}_k$ the conditional expectation onto the von Neumann algebra
$\hat{M}_k$ generated by $\rho(N_k)$ and
$\pi_1(N_1),....,\pi_k(N_k)$. Then we deduce from being successively
independent that
 \begin{align*}
 \|dx\|_{BMO_c}&= \sup_n \|\hat{E}_n(\sum_{k\gl n}
 \pi_k(d_k(x)^2))\| \sim_2 \|\pi_k(d_k(x)^2)\|+
 \|\sum_{k>n} \hat{E}_n\hat{E}_{k-1}(\pi_k(d_k(x)^2))\| \\
 &= \|\pi_k(d_k(x)^2)\|+
 \|\sum_{k>n} \hat{E}_n\hat{E}_{\rho}(\pi_k(d_k(x)^2))\| \\
 &= \|\pi_k(d_k(x)^2)\|+
 \|\sum_{k>n} \rho(\sum_k E_{k-1}(d_k(x)^2))\| \sim_2 \|x\|_{BMO_c} \pl .
 \end{align*}
Therefore $d$ yields an isomorphic embedding $d:BMO_{c/r}(N_k)\to
BMO_{c/r}(\hat{M}_k)$. For $1\le p\le 2$, we see that similarly $d$
is bounded on $h_{p}^d$, $h_{p}^c$ and $h_{p}^r$. However, for $1\le
p <2$, we know that $H_p^c=h_p^d+h_p^c$ holds with equivalent norms
and hence $d$ is also continuous on $H_p^c$. Using
$tr(dx^*dy)=tr(x^*y)$ it then follows that $d$ is an isomorphism on
$H_p^c$ and $H_p^r$ for all $1\le p<\infty$ and on $BMO_{c}$,
$BMO_r$. The assertion follows.\qd

\begin{lemma}\label{filtra} For a group $G=G_0\unrhd G_1 \cdots $
there is a canonical tangent dilation.
\end{lemma}

\begin{proof}  Let $\tilde{G}=G\times \zz$. Then we define
 \[ \pi_k:G_k\to G\times \zz \pl ,\pl \pi(g) \lel (g,gG_{k+1}) \pl
 .\]
Let $E:VN(\tilde{G})\to VN(\tilde{G})$ be the conditional
expectation onto $VN(G)$ and $\rho$ the canonical embedding . Then
clearly,
 \[ E(\pi_k(g))\lel \begin{cases} g & g\in G_{k+1} \\
                                  0& \mbox{else}
                                  \end{cases} \pl .\]
This means $E(\pi_k(g))\lel \rho(E_{G_{k+1}}(g))$.  Note here that
for a reversed filtration the definition of tangent filtration has
to be suitably modified. Finally, let $\tilde{G}_k\subset \tilde{G}$
be the subgroup generated by $\rho(G_k)$ and $\zz$. Let $g\in
G_{k-1}$. Then we have
 \[ E_{\tilde{G}_{k}}(\pi_{k-1}(g)) \lel 1_{g\in G_k}(g) (g,gG_{k})
 \lel E(\pi_{k-1}(g)) \pl, \]
because only for $g\in G_k$ we have a non-trivial term. \qd

In the following we consider the Hilbert transform
 \[ H(g) \lel i (1_{P}(g)-1_{P}(g^{-1}))g \pl \]
induced by the order. Note that
 \[ H(g)^* \lel H(g^{-1}) \pl .\]

\begin{lemma} Let $P_t^{\zz}$ be the Poisson semigroup on
$VN(\zz)$ and $P_t=id\ten P_t^{\zz}$ the Poisson semigroup with
generator $A$ and gradient form $\Gamma$. Then
 \[ \Gamma(dHx,dHx) \lel \Gamma(dx,dx) \]
and
 \[  P_t|dHx|^2-|P_tdHx|^2 \lel  P_t|dx|^2-|P_tdx|^2 \]
holds for all $x\in \cz[G]$.
\end{lemma}

\begin{proof} For the first assertion we consider $g,h \in G$ and
$i,j$ such that $g_1\in G_i\setminus G_{i+1}$, $g_2\in G_j\setminus
G_{j+1}$. If $g_1G_{i+1}\gl 0$ and $g_2G_{j+1}\gl 0$ or
$g_1G_{i+1}\le 0$ and $g_2G_{j+1}\gl 0$ we have
 \[ \Gamma(dHg_1,dHg_2) \lel \Gamma(dg_1,dg_2) \pl .\]
The interesting case is given by $k=g_1G_{i+1}\gl 0$ and
$j=g_2G_{i+1}\le 0$. Then we note that
 \[ \Gamma(\la(k),\la(j))
 \lel \frac{|k|+|j|-|k-j|}{2} \la(k)^*\la(j) \lel 0 \pl .\]
The second assertion follows similarly in this case
 \[ P_t^{\zz}(\la(k)^*\la(j))-P_t^{\zz}(\la(k))^*P_t^{\zz}(\la(j))
 \lel (e^{-t|j-k|}-e^{-t|k|}e^{-t|j|})\la(k)^*\la(j) \lel 0 \pl
 .\]
Thus the sign change only occurs when we have a $0$-coefficient. \qd

\begin{theorem} $H$ is bounded on $L_p(VN(G))$ for all
$1<p<\infty$.
\end{theorem}

\begin{proof} Let $x$ be selfadjoint. Theorem
\ref{stein3}iii)  and $H^{\infty}$-calculus implies
 \[ \|dx\|_p \sim \|(\int \Gamma(P_sdx,P_sdx) ds)^{1/2}\|
 \lel \|(\int \Gamma(P_sdHx,P_sdHx) ds)^{1/2}\| \sim \|Hx\|_p \pl
 .\]
Here we use that for selfadjoint $x$ the element $Hx$ is also
selfadjoint. For an alternative proof, we can use interpolation and
note that $dH:VN(G)\to BMO(\tilde{G})$ is bounded and
$dH:L_2(VN(G))\to L_2(VN(\tilde{G}))$ and hence
 \[ dH:L_p(VN(G))\to L_p(VN(\tilde{G})) \]
is bounded. Finally, the assertion follows from Lemma \ref{rosen}
 \begin{align*}
  \|dHx\|_p \sim_{c(p)} \|Hx\|_p
  \end{align*}
This yields the assertion for $2\le p<\infty$. Duality  implies the
result for $1<p\le 2$ as well. \qd

\begin{rem}{\rm With the help of the tangent dilation, we can show that every completely bounded Fourier multiplier on $\zz$ induces a Fourier multiplier on $VN(G)$, by applying it to
$d(x)$. We can even use different multipliers on $\zz^{\infty}$ by
modifying $d(g)\lel (g,...,gG_{k+1},...)$ for $g\in G_k\setminus
G_{k+1}$.
   }
\end{rem}



\end{document}